\documentclass[10pt]{article}
\usepackage[french]{babel}
\usepackage{amssymb}
\usepackage{graphicx}
\usepackage{psfrag}
\usepackage{graphics}
\usepackage{amsthm}

\input pictex.tex

\newcommand{\clo}{\mathrm{S}^1}

\newcommand{\+}{\tiny{+}}
\newcommand{\esp}{\thinspace}
\theoremstyle{definition}
 
\newtheorem{thm}{Th\'eor\`eme}[section]

\newtheorem{prop}[thm]{Proposition}

\newtheorem{lem}[thm]{Lemme}

\newtheorem{rem}[thm]{Remarque}

\newtheorem{defn}[thm]{D\'efinition}

\input epsf

\begin{document} 

\date{}
\author{Deroin, B. \qquad Kleptsyn, V. \qquad Navas, A.}

\title{Sur la dynamique unidimensionnelle en r\'egularit\'e interm\'ediaire}
\maketitle

\vspace{-0.5cm}

\noindent{\bf R\'esum\'e.} Par des m\'ethodes 
de nature probabiliste, nous d\'emontrons de nouveaux 
r\'esultats de rigidit\'e pour des groupes et des pseudo-groupes 
de diff\'eomorphismes de vari\'et\'es unidimensionnelles dont la 
classe de diff\'erentiabilit\'e est interm\'ediaire ({\em i.e.} entre $C^1$ 
et $C^2$). En particulier, nous prouvons des g\'en\'eralisations du 
th\'eor\`eme de Denjoy et d'un lemme classique de Kopell pour des groupes 
ab\'eliens. Ensuite, nous appliquons les techniques introduites \`a l'\'etude 
des feuilletages de codimension 1 dont la r\'egularit\'e transverse 
est interm\'ediaire. Nous obtenons notamment des versions g\'en\'eralis\'ees 
du th\'eor\`eme de Sacksteder en classe $C^1$. Nous finissons par 
quelques remarques \`a propos de la mesure stationnaire.

\vspace{0.15cm}

\noindent{\bf Abstract.} Using probabilistic methods, we prove 
new rigidity results for groups and pseudo-groups of diffeomorphisms 
of one dimensional manifolds with intermediate regularity class ({\em i.e.} 
between $C^1$ and $C^2$). In particular, we demonstrate some generalizations 
of Denjoy's Theorem and the classical Kopell's Lemma for Abelian groups. These 
techniques are also applied to the study of codimension 1 foliations. 
For instance, we obtain several generalized versions of Sacksteder's Theorem 
in class $C^1$. We conclude with some remarks about the stationary measure.

\vspace{0.15cm}

\noindent{\bf Mots-cl\'es :} diff\'eomorphismes, centralisateurs, 
rigidit\'e, feuilletages, mesure stationnaire.

\vspace{0.15cm}


\vspace{0.1cm}


\markright{\sc Dynamique unidimensionnelle en r\'egularit\'e interm\'ediaire}

\thispagestyle{empty}

\tableofcontents

\newpage

\section*{Introduction}
\addcontentsline{toc}{section}{Introduction}

\subsection*{Quelques motivations}
\addcontentsline{toc}{subsection}{Quelques motivations}
 
Dans ce travail, nous nous int\'eressons \`a la dynamique 
des groupes de diff\'eomorphismes de vari\'et\'es unidimensionnelles dont 
la classe de diff\'erentiabilit\'e est {\em interm\'ediaire}, c'est-\`a-dire 
sup\'erieure ou \'egale \`a $C^1$ mais strictement 
inf\'erieure \`a $C^2$. \`A premi\`ere 
vue, cela peut para\^{\i}tre un sujet d'\'etude \`a la fois artificiel et 
tr\`es technique. Il y a n\'eanmoins plusieurs raisons qui le justifient.

La premi\`ere est d'ordre dynamique. On sait depuis le c\'el\`ebre travail 
de Denjoy \cite{denjoy} que l'\'etude de la dynamique des diff\'eomorphismes 
de classe $C^2$ du cercle dont le nombre de rotation est irrationnel 
se ram\`ene (\`a conjugaison topologique pr\`es) \`a celle des rotations. 
Ce n'est plus le cas en classe $C^1$, ce qui a \'et\'e 
illustr\'e aussi dans \cite{denjoy} par des 
\begin{tiny}$^{_\ll}$\end{tiny}contre-exemples\hspace{0.05cm}\begin{tiny}$^{_\gg}$\end{tiny} 
(voir \'egalement \cite{bol}). Au \S 3 du chapitre X de \cite{herman}, et en suivant 
une id\'ee de Milnor et Sergeraert, Herman a am\'elior\'e ces contre-exemples en 
les rendant de classe $C^{1+\tau}$ pour tout $\tau \!<\! 1$. Au \S 5 de ce m\^eme 
chapitre, on trouve des contre-exemples de Denjoy continus avec des 
centralisateurs non triviaux. Il est bien possible que Herman ait d\'ej\`a 
r\'efl\'echi au probl\`eme d'am\'eliorer la classe de diff\'erentiabilit\'e 
de ces derniers. Ce probl\`eme est en fait au c\oe ur de notre travail.

Une autre raison pour entreprendre cette \'etude, de nature cohomologique, 
est li\'ee \`a l'invariant de Godbillon-Vey. Rappelons que pour des groupes 
agissant sur des vari\'et\'es unidimensionnelles, cet invariant correspond 
\`a une classe de 2-cohomologie (\`a valeurs r\'eelles) dont le 
repr\'esentant (non homog\`ene) est donn\'e par
$$\mathcal{GV}(f,g) = \int \log(f') \thinspace d\log(g' \!\circ\! f).$$
Cette d\'efinition a un sens \'evident lorsque les diff\'eomorphismes 
$f$ et $g$ sont de classe $C^2$. Or, on s'aper\c{c}oit facilement 
que le domaine de d\'efinition le plus naturel est celui des 
diff\'eomorphismes dont le logarithme de la d\'eriv\'ee appartient \`a 
l'espace de Sobolev $H^{1/2}$. En particulier, l'invariant peut \^etre 
d\'efini (et il est continu) dans le groupe des diff\'eomorphismes 
de classe $C^{3/2 + \varepsilon}$ pour tout $\varepsilon > 0$. Cette 
\begin{tiny}$^{_\ll}$\end{tiny}extension\esp\begin{tiny}$^{_\gg}$\end{tiny} 
propos\'ee par Hurder et Katok dans \cite{HK} 
(voir aussi \cite{Ts-area}) est \`a comparer avec le travail 
remarquable \cite{Ts2} de Tsuboi, o\`u il est d\'emontr\'e (parmi bien 
d'autres choses) qu'aucune extension n'est envisageable en classe $C^1$ 
(voir \'egalement \cite{Ts3,TsP}). La preuve de Tsuboi s'appuie sur 
l'existence d'une famille assez large de contre-exemples $C^1$ (d\^us 
\`a Pixton \cite{Pix}) \`a un lemme classique de Kopell \cite{Ko}, selon 
lequel le centralisateur de tout diff\'eomorphisme de classe $C^2$ de 
l'intervalle $[0,1[$ sans point fixe \`a l'int\'erieur agit librement sur $]0,1[$. 
Si l'on tient compte de l'extension de Hurder et Katok, on devrait s'attendre 
\`a ce que de tels contre-exemples ne puissent pas \^etre construits en classes 
interm\'ediaires trop \'el\'ev\'ees. Cette motivation a amen\'e 
Tsuboi \`a construire des contre-exemples au lemme de Kopell (et aussi 
au th\'eor\`eme de Denjoy avec des centralisateurs non triviaux) dont la r\'egularit\'e 
est sup\'erieure \`a $C^1$~\cite{TsP}. Il a \'egalement conjectur\'e des ph\'enom\`enes de 
rigidit\'e en classe interm\'ediaire, mais il n'a pas donn\'e de r\'esultats 
dans cette direction. Notre travail confirme l'intuition de Tsuboi, et de 
plus il montre que ses constructions sont en fait optimales en ce qui 
concerne la r\'egularit\'e atteinte.

Signalons en passant que la classe de diff\'erentiabilit\'e $C^{3/2}$ appara\^{i}t de mani\`ere 
naturelle dans d'autres situations. Par exemple, dans \cite{Na-kazh,Na-superrig} elle se 
trouve \^etre la classe critique pour la d\'efinition d'un autre invariant cohomologique : 
le \begin{tiny}$^{_\ll}$\end{tiny}cocycle de Liouville\hspace{0.05cm}\begin{tiny}$^{_\gg}$\end{tiny}. 
Ce cocycle est un outil important dans l'\'etude de ph\'enom\`enes de rigidit\'e (dans l'esprit 
des th\'eor\`emes de Margulis et Zimmer) pour les actions de groupes 
\begin{tiny}$^{_\ll}$\end{tiny}de rang sup\'erieur\esp\begin{tiny}$^{_\gg}$\end{tiny} 
par diff\'eomorphismes du cercle (voir \cite{mori} pour un lien entre le cocycle de 
Liouville et la classe de Godbillon-Vey ; voir aussi \cite{NS}). 
D'autre part, la r\'esolution d'\'equations diff\'erentielles 
stochastiques dans le groupe des diff\'eomorphismes du cercle donne lieu 
\`a des analogues de la mesure de Wiener sur $\mathrm{Diff}_+^{\infty}(\clo)$ 
dont le support s'av\`ere \^etre pr\'ecis\'ement l'espace des diff\'eomorphismes 
de classe $C^{3/2}$ (voir \cite{mall,mal}). \`A notre connaissance, le cas 
$d=2$ des th\'eor\`emes A et B que nous verrons plus loin correspondent aux 
premiers r\'esultats de nature essentiellement dynamique et li\'es \`a la 
classe $C^{3/2}$.

La derni\`ere justification de notre \'etude concerne la th\'eorie 
les feuilletages de codimension 1. En effet, en s'inspirant du 
th\'eor\`eme classique de Denjoy, cette th\'eorie 
s'est developp\'ee en admettant la plupart des cas 
une hypoth\`ese de r\'egularit\'e transverse $C^2$. Or, si l'on tient 
compte de sa nature dynamique, et du fait que de nombreux r\'esultats 
fondamentaux en syst\`emes dynamiques, notamment la th\'eorie de Pesin, 
sont encore valables en classe $C^{1+\tau}$ (et parfois m\^eme en classe $C^1$), 
on peut esp\'erer que plusieurs propri\'et\'es des feuilletages de classe $C^2$ soient 
encore valables dans ce contexte (quelques progr\`es importants dans cette 
direction sont d\'ej\`a connus \cite{DK,hurdito,Hu,hurder1,hurder2,HL}). 
Si tel est le cas, on devrait s'attendre \'egalement \`a voir 
appara\^{\i}tre des obstructions en r\'egularit\'e interm\'ediaire.


\subsection*{Pr\'esentation des r\'esultats}
\addcontentsline{toc}{subsection}{Pr\'esentation des r\'esultats}

Suivant un principe exprim\'e par Herman dans l'introduction 
de \cite{herman}, \begin{tiny}$^{_\ll}$\end{tiny}tout chercheur d\'esireux de 
travailler sur les diff\'eomorphismes du cercle doit s'habituer \`a construire 
et \'etudier des exemples\hspace{0.05cm}\begin{tiny}$^{_\gg}$\end{tiny}. Ce point 
de vue est illustr\'e par le fait que Denjoy lui-m\^eme ait abouti 
\`a son c\'el\`ebre th\'eor\`eme parce que \begin{tiny}$^{_\ll}$\end{tiny}il 
ne pouvait pas en construire des contre-exemples de classe 
$C^2$\hspace{0.05cm}\begin{tiny}$^{_\gg}$\end{tiny} 
(des contre-exemples qui, n\'eanmoins, \'etaient vraisemblables \`a Poincar\'e). 
Dans ce m\^eme esprit, 
si l'on essaie de construire des contre-exemples en r\'egularit\'e 
inf\'erieure \`a $C^2$ \`a plusieurs r\'esultats de dynamique unidimensionnelle, 
on voit appara\^{\i}tre des obstructions pour des classes de diff\'erentiabilit\'e 
interm\'ediaires tr\`es pr\'ecises. La premi\`ere est li\'ee \`a des actions de 
groupes ab\'eliens sur le cercle.

\vspace{0.35cm}

\noindent{\bf Th\'eor\`eme A.} {\em Si $d$ est un entier sup\'erieur ou \'egal 
\`a $2$ et $\varepsilon\! >\! 0$, alors toute action libre de $\mathbb{Z}^d$ par 
diff\'eomorphismes de classe $C^{1+1\!/\!d+\varepsilon}$ du cercle est minimale.}

\vspace{0.35cm}

Ce r\'esultat peut \^etre consid\'er\'e comme une g\'en\'eralisation du th\'eor\`eme de 
Denjoy (le cas $d = 1$) avec une petite hypoth\`ese de r\'egularit\'e suppl\'ementaire 
(donn\'ee par le $\varepsilon > 0$). Il est bien possible qu'il soit encore valable 
en classe $C^{1+1\!/\!d}$. Nous en donnons une preuve simple dans ce cas sous une hypoth\`ese 
dynamique assez naturelle (voir la proposition~\ref{tangente}). Signalons que le 
th\'eor\`eme A a \'et\'e conjectur\'e 
dans \cite{TsP} par Tsuboi, qui l'a \'egalement illustr\'e par des contre-exemples de 
classe $C^{1+1\!/\!d-\varepsilon}$ pour tout $\epsilon > 0$ (nous discuterons de telles 
constructions au \S \ref{exemples}). Il serait sans doute int\'eressant d'obtenir, sous les 
hypoth\`eses du th\'eor\`eme, un r\'esultat d'ergodicit\'e (par rapport \`a la mesure 
de Lebesgue) similaire \`a celui valable pour les diff\'eomorphismes de classe $C^2$ 
du cercle dont le nombre de rotation est irrationnel (voir \cite{HasK,herman}). 
Quant \`a la r\'egularit\'e de la lin\'earisation, le lecteur trouvera dans le 
\S 3 de \cite{kra} des r\'esultats optimaux (qui g\'en\'eralisent ceux de Moser 
\cite{moser}) pour des groupes commutatifs engendr\'es par des petites perturbations 
de rotations\footnote{Malheureusement, nous devons signaler que les exemples 
du \S 2 de \cite{kra} sont erron\'es. Ceci est une cons\'equence directe 
de la proposition 1.1, et r\'esulte \'egalement du th\'eor\`eme A.}.

D\^u \`a l'absence de r\'ecurrence pour la dynamique, le cas o\`u il existe des points fixes 
globaux, {\em i.e.} le cas de l'intervalle, est l\'eg\`erement diff\'erent. Dans ce contexte, 
on a un analogue du th\'eor\`eme de Denjoy~; c'est le th\'eor\`eme de Kopell-Szekeres 
concernant la rigidit\'e des centralisateurs de diff\'eomorphismes de classe $C^2$ (voir par 
exemple \cite{yoccoz}). Le r\'esultat suivant peut \^etre 
vu comme une g\'en\'eralisation du fameux {\em lemme de Kopell} 
(le cas $d=1$) pour les actions de $\mathbb{Z}^{d+1}$ (sous une l\'eg\`ere hypoth\`ese 
suppl\'ementaire de r\'egularit\'e). En effet, l'hypoth\`ese combinatoire (\ref{comb}) 
ci-dessus est exactement celle que l'on trouve lors de la d\'emonstration du lemme 
de Kopell (voir le \S \ref{ex2}). De plus, c'est en utilisant des diff\'eomorphismes 
(de petite r\'egularit\'e) qui v\'erifient cette propri\'et\'e que Tsuboi a 
d\'emontr\'e son th\'eor\`eme d'acyclicit\'e cohomologique dans \cite{Ts2}.

\vspace{0.4cm}

\noindent{\bf Th\'eor\`eme B.} {\em Soient $d\!\geq\!2$ un entier et 
$\varepsilon \!>\! 0$. Soient $f_1,\ldots,f_{d+1}$ des diff\'eomorphismes de 
classe $C^1$ de l'intervalle $[0,1]$ qui commutent entre eux. Supposons qu'il 
existe des intervalles ouverts disjoints $I_{n_1,\ldots,n_{d}}$ qui 
sont dispos\'es dans $]0,1[$ en respectant l'ordre lexicographique et tels 
que, pour tout $(n_1,\ldots,n_d) \! \in \! \mathbb{Z}^d$,}
\begin{equation}
f_i(I_{n_1,\ldots,n_i,\ldots,n_{d}}) = I_{n_1,\ldots,n_i-1,\ldots,n_{d}}  
\quad \mbox{pour tout} \quad i\! \in\! \{1,\ldots,d\} \quad 
\mbox {et} \quad f_{d+1}(I_{n_1,\ldots,n_d}) =I_{n_1,\ldots,n_d}.
\label{comb}
\end{equation}
{\em Si $f_1,\ldots,f_d$ sont de classe $C^{1+1\!/\!d+\varepsilon}$, alors 
la restriction de $f_{d+1}$ \`a la r\'eunion des $I_{n_1,\ldots,n_d}$ est l'identit\'e.}

\vspace{0.4cm}

Une nouvelle fois, la classe de diff\'erentiabilit\'e $C^{1\!/\!d}$ pour la d\'eriv\'ee est 
optimale. La construction de contre-exemples, bien plus d\'elicate que dans le cas 
du cercle, appara\^{\i}t aussi dans \cite{TsP}. 

Dans les deux r\'esultats pr\'ec\'edents, la classe interm\'ediaire 
optimale est li\'ee \`a la  
\begin{tiny}$^{_\ll}$\end{tiny}croissance\hspace{0.05cm}\begin{tiny}$^{_\gg}$\end{tiny} 
du groupe (ou plut\^ot des orbites des actions respectives). Il serait 
int\'eressant (et pas tr\`es difficile) d'obtenir des r\'esultats de ce type 
pour des actions de groupes nilpotents, \`a croissance sous-exponentielle ou 
r\'esolubles, qui soient des extensions en classe interm\'ediaire de ceux obtenus 
dans \cite{PT}, \cite{Na-subexp,ceuno} et \cite{Na-solv} respectivement. En 
pr\'esence d'une dynamique \begin{tiny}$^{_\ll}$\end{tiny}\`a croissance 
exponentielle\hspace{0.05cm}\begin{tiny}$^{_\gg}$\end{tiny}, 
on voit appara\^{\i}tre des ph\'enom\`enes de rigidit\'e 
analogues en classe $C^{1+\tau}$ pour tout $\tau \!>\! 0$. De plus, en 
utilisant des m\'ethodes un peu techniques mais tout \`a fait standard 
en syst\`emes dynamiques, on peut d\'emontrer que certains de ces 
ph\'enom\`enes ont encore lieu en classe $C^1$. Un exemple est donn\'e 
par le r\'esultat ci-dessous, lequel peut \^etre pens\'e comme une 
version g\'en\'eralis\'ee du th\'eor\`eme de Sacksteder.

\vspace{0.35cm}

\noindent{\bf Th\'eor\`eme C.} {\em Soit $\mathcal{F}$ un feuilletage 
de codimension 1 et transversalement de classe $C^1$. Si $\mathcal{F}$ 
n'admet pas de mesure transverse invariante (au sens de Plante} 
\cite{Pl-anals}), {\em alors il poss\`ede des feuilles ressort hyperboliques.}

\vspace{0.35cm}

Ce r\'esultat n'est absolument pas nouveau. Dans toute sa g\'en\'eralit\'e, 
il doit \^etre atribu\'e \`a Hurder \cite{hurder2}, si bien que des r\'esultats 
reli\'es (mais plus faibles) se trouvent dans \cite{candel,DK,GLW}. La d\'emonstration 
que nous proposons est \'el\'ementaire et s'applique en g\'eneral \`a des pseudo-groupes 
de diff\'eomorphismes en dimension 1. En effet, d'apr\`es la th\'eorie bien connue 
des feuilletages de codimension 1, la preuve du th\'eor\`eme C se r\'eduit 
\`a d\'emontrer que, en pr\'esence d'une feuille ressort topologique, il 
existe des feuilles ressort hyperboliques (voir l'appendice \ref{ap0}). 
Or, c'est exactement cette affirmation que nous d\'emontrons~; la m\'ethode 
est simple et directe, et elle pourrait \^etre utile dans d'autres circonstances. 

Le probl\`eme de la validit\'e du th\'eor\`eme de Sacksteder en classe inf\'erieure \`a $C^2$ 
a \'et\'e d'abord abord\'e par Hurder pour les groupes de diff\'eomorphismes du 
cercle~\cite{hurdito,Hu}. En pr\'esence d'un minimal exceptionnel avec une dynamique 
\begin{tiny}$^{_\ll}$\end{tiny}suffisamment riche\hspace{0.05cm}\begin{tiny}$^{_\gg}$\end{tiny}, 
il d\'emontre l'existence d'\'el\'ements avec des points fixes hyperboliques. 
Nous proposons ici une version \begin{tiny}$^{_\ll}$\end{tiny}globale\begin{tiny}$^{_\gg}$\end{tiny} 
de ce fait.

\vspace{0.35cm}

\noindent{\bf Th\'eor\`eme D.} {\em Si $\Gamma$ est un sous-groupe d\'enombrable de 
$\mathrm{Diff}^1_+(\clo)$ qui ne pr\'eserve aucune mesure de probabilit\'e du cercle, 
alors $\Gamma$ contient des \'el\'ements n'ayant que des points fixes hyperboliques.}

\vspace{0.35cm}

Remarquons que l'hypoth\`ese de non existence d'une mesure de probabilit\'e invariante 
ci-dessus \'equivaut \`a ce que le groupe $\Gamma$ ne soit pas semi-conjugu\'e \`a 
un groupe de rotations et qu'il n'ait pas d'orbite finie. 

Tel qu'il est \'enonc\'e, le th\'eor\`eme D g\'en\'eralise les r\'esultats obtenus par 
Hurder dans \cite{hurdito,Hu} dans trois directions : l'hypoth\`ese faite sur 
la dynamique (nous admettons des 
groupes non ab\'eliens avec toutes ses orbites denses), le nombre de points 
fixes hyperboliques, et la non existence d'autres points fixes. La premi\`ere 
direction n'est pas nouvelle, car il avait \'et\'e d\'ej\`a 
remarqu\'e par Ghys que, en classe $C^2$, le th\'eor\`eme de Sacksteder 
(qui en g\'en\'eral est pr\'esent\'e comme un r\'esultat qui n'est 
valable qu'en pr\'esence d'un minimal exceptionnel) reste valide 
pour des groupes non commutatifs qui agissent de mani\`ere minimale 
(voir la page 11 de \cite{ET}). De plus, Hurder a r\'ecemment d\'emontr\'e 
que tout sous-groupe de $\mathrm{Diff}_+^1(\clo)$ sans probabilit\'e 
invariante poss\`ede des \'el\'ements avec des points fixes hyperboliques 
\cite{hurder1}. La deuxi\`eme direction, suivant laquelle il existe des 
\'el\'ements avec des points fixes hyperboliquement dilatants et contractants, 
est une am\'elioration non banale d\'ej\`a en classe $C^2$. Nous en donnons une 
preuve simple dans ce contexte au \S \ref{sac-circ}. Finalement, l'obtention d'un 
\'el\'ement n'ayant que des points fixes hyperboliques est un raffinement bien 
plus subtil. En effet, m\^eme dans le cas analytique r\'eel, on ne dispose 
d'aucune d\'emonstration de ce fait par des m\'ethodes 
\begin{tiny}$^{_\ll}$\end{tiny}standard\hspace{0.05cm}\begin{tiny}$^{_\gg}$\end{tiny}. 
Pour la preuve nous nous appuyons fortement sur des r\'esultats concernant 
la mesure stationnaire et les exposants de Lyapunov d'une action, 
lesquels portent un inter\^et en soi et sont inclus dans l'appendice de ce travail.


\subsection*{Sur la m\'ethode de d\'emonstration}
\addcontentsline{toc}{subsection}{Sur la m\'ethode de d\'emonstration}

L'un des aspects techniques les plus importants de la 
dynamique unidimensionnelle diff\'erentiable est li\'e \`a la possibilit\'e 
de contr\^oler la distorsion (lin\'eaire ou projective) des it\'er\'ees d'une 
application (ou de compos\'ees de diff\'erentes applications). Dans le cas d'une dynamique 
\begin{tiny}$^{_\ll}$\end{tiny}inversible\hspace{0.05cm}\begin{tiny}$^{_\gg}$\end{tiny}, 
{\em i.e.} lorsqu'il n'y a pas de point critique, deux principes sont plus au moins 
canoniques pour contr\^oler la distorsion lin\'eaire. Le premier vient du th\'eor\`eme 
de Denjoy : en classe $C^2$, on peut contr\^oler la distorsion des it\'er\'ees sur des 
intervalles disjoints en la comparant avec la somme des longueurs de 
ces intervalles. Le deuxi\`eme principe, dit du 
\begin{tiny}$^{_\ll}$\end{tiny}folklore\hspace{0.05cm}\begin{tiny}$^{_\gg}$\end{tiny} 
dans \cite{SS}, est valable en classe $C^{1+\tau}$ lorsqu'on sait {\em a priori} 
qu'il y a suffisamment d'hyperbolicit\'e. 

Dans ce travail nous introduisons 
des nouvelles m\'ethodes pour contr\^oler la distorsion en classe interm\'ediaire 
lorsque la dynamique a une structure topologique et combinatoire bien pr\'ecise 
et qui appara\^{\i}t de mani\`ere naturelle dans de nombreuses situations. 
Ces m\'ethodes s'appuient sur des arguments de nature probabiliste 
(inspir\'es en partie de \cite{DK}) : nous cherchons \`a assurer 
un contr\^ole uniforme pour la distorsion d'une suite al\'eatoire 
\begin{tiny}$^{_\ll}$\end{tiny}typique\hspace{0.05cm}\begin{tiny}$^{_\gg}$\end{tiny} 
de compositions (si bien qu'en g\'en\'eral il est impossible de d\'eterminer 
pr\'ecis\'ement quelle suite conviendra !). Nous croyons que, 
au del\`a m\^eme des th\'eor\`emes pr\'esent\'es, 
ces m\'ethodes soient celles qui donnent le plus d'inter\^et \`a cet 
article. D'ailleurs, il est tr\`es raisonnable d'essayer de les utiliser afin 
de raffiner d'autres r\'esultats classiques de la th\'eorie des feuilletages 
de codimension 1, comme par exemple ceux de Duminy \cite{CC-dum,yo} (voir \`a ce propos 
\cite{HL}). Il pourrait s'av\'erer \'egalement int\'eressant de repenser en r\'egularit\'e 
interm\'ediaire la {\em th\'eorie des niveaux} \cite{CClibro,hector}, tout 
en s'appuyant sur les techniques introduites dans ce travail.

\vspace{0.45cm}

\noindent{\bf Quelques recommandations pour la lecture.} La structure de cet article 
n'est pas toujours lin\'eaire. Le lecteur press\'e pour la d\'emonstration du th\'eor\`eme 
de Denjoy g\'en\'eralis\'e peut aller directement au \S 2 (si bien que la lecture des 
exemples du \S 1.1 peut aider \`a \'eclaircir un peu le panorama). Celui qui est 
plut\^ot int\'eress\'e par le th\'eor\`eme de Sacksteder peut passer imm\'ediatement au 
\S 2.1 et puis aux \S 4.1 et \S 4.2 (o\`u l'on d\'emontre les versions g\'en\'eralis\'ees 
pour des pseudo-groupes). Le cas des diff\'eomorphismes du cercle n\'ecessite 
cependant les r\'esultats probabilistes de l'appendice (qui peuvent \^etre 
\'etudi\'es sans la n\'ecessit\'e du reste de l'article). Nous proposons ci-dessous un 
sch\'ema d'interd\'ependance logique entre les diff\'erents paragraphes. Une fl\`eche 
continue repr\'esente une lecture indispensable, alors qu'une fl\`eche pointill\'ee 
indique une lecture non indispensable mais tout \`a fait convenable pour la compr\'ehension.

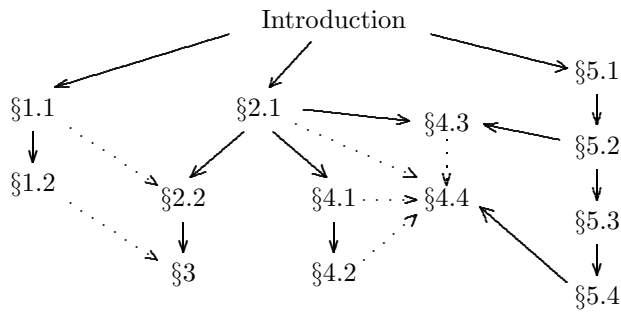
\vspace{0.45cm}


\beginpicture

\setcoordinatesystem units <1cm,1cm>

\put{$\mbox{Introduction}$} at 0 3.2 
\put{$\S 1.1$} at -4 2 
\put{$\S 1.2$} at -4 1 
\put{$\S 2.1$} at -1 2 
\put{$\S 2.2$} at -2 0.8 
\put{$\S 3$} at -2 -0.2 
\put{$\S 4.1$} at 0 0.8 
\put{$\S 4.2$} at 0 -0.2 
\put{$\S 4.3$} at 1.5 1.8 
\put{$\S 4.4$} at 1.5 0.8 
\put{$\S 5.1$} at 3.5 2.5 
\put{$\S 5.2$} at 3.5 1.5 
\put{$\S 5.3$} at 3.5 0.5 
\put{$\S 5.4$} at 3.5 -0.5 

\putrule from -4 1.3 to -4 1.7 
\putrule from -2 0.1 to -2 0.5 
\putrule from 0 0.1 to 0 0.5 
\putrule from 3.5 1.8 to 3.5 2.2 
\putrule from 3.5 0.8 to 3.5 1.2 
\putrule from 3.5 -0.2 to 3.5 0.2 

\plot 
-4 1.3 
-4.05 1.45 /

\plot 
-4 1.3 
-3.95 1.45 /

\plot 
-2 0.1 
-2.05 0.25 /

\plot 
-2 0.1 
-1.95 0.25 /

\plot 
0 0.1 
0.05 0.25 /

\plot 
0 0.1 
-0.05 0.25 /

\plot 
3.5 1.8 
3.55 1.95 /

\plot 
3.5 1.8 
3.45 1.95 /

\plot 
3.5 0.8 
3.55 0.95 /

\plot 
3.5 0.8 
3.45 0.95 /

\plot 
3.5 -0.2 
3.55 -0.05 /

\plot 
3.5 -0.2 
3.45 -0.05 /


\plot 3.1 -0.3 
1.95 0.7 /

\plot 1.95 0.7 2.035 0.51 /
\plot 1.95 0.7 2.13 0.65 /

\plot 
-1.4 3 
-3.7 2.3 /

\plot -3.7 2.3 -3.52 2.305 /
\plot -3.7 2.3 -3.55 2.4 /

\plot 
-0.3 2.9 
-0.85 2.3 /

\plot -0.85 2.3 -0.82 2.42 /
\plot -0.85 2.3 -0.7 2.36 /

\plot 
1.3 3 
3.1 2.55 /

\plot 3.1 2.55 2.95 2.5 /
\plot 3.1 2.55 2.98 2.65 /

\plot 
-1.2 1.7 
-1.9 1.1 /

\plot -1.9 1.1 -1.705 1.17 /
\plot -1.9 1.1 -1.83 1.25 /

\plot 
-0.8 1.7
-0.1 1.1 /

\plot -0.1 1.1 -0.28 1.145 /
\plot -0.1 1.1 -0.16 1.26 /

\plot 1.1 1.1 1 1.22 /
\plot 1.1 1.1 0.93 1.1 /

\plot 
-0.4 2 
1 1.85 /

\plot 1 1.85 0.86 1.93 /
\plot 1 1.85 0.84 1.8 /

\plot
3 1.6
2 1.8 /

\plot 2 1.8 2.18 1.82 /
\plot 2 1.8 2.155 1.7 /


\plot 1.5 1 1.55 1.15 /
\plot 1.5 1 1.45 1.15 /

\plot 1.1 0.8 0.95 0.85 /
\plot 1.1 0.8 0.95 0.75 /

\plot 1.1 0.6 0.95 0.55 /
\plot 1.1 0.6 1.04 0.455 /

\plot 
-2.35 1 
-2.43 1.12 /

\plot 
-2.35 1 
-2.5 1.03 /

\plot 
-2.35 0 
-2.43 0.12 /

\plot 
-2.35 0 
-2.5 0.03 /


\setdots

\putrule from 1.5 1.6 to 1.5 1 

\plot 
-0.5 1.8 
1.1 1.1 /

\plot
-3.5 1.75 
-2.35 1 /

\plot
-3.5 0.75 
-2.35 0 /

\plot 
0.4 0.8 
1.1 0.8 /

\plot 
0.4 0 
1.1 0.6 /

\put{} at -8.1 0

\endpicture


\vspace{0.35cm}

\noindent{\bf Remerciemments.} Nous sommes tr\`es reconnaissants envers \'Etienne Ghys pour 
avoir partag\'e avec nous ses id\'ees et connaissances sur les groupes de diff\'eomorphismes 
du cercle, ainsi qu'envers Sylvain Crovisier, sans l'aide de qui nous n'aurions jamais r\'ealis\'e 
la validit\'e des versions $C^1$ (et non seulement $C^{1+\tau}$) du th\'eor\`eme de Sacksteder. Nous 
exprimons aussi notre gratitude \`a Jean Christophe Yoccoz, Danijela Damjanovic et Anna Erschler 
pour l'int\^eret qu'ils ont port\'e \`a cet article. Ce travail s'est d\'eroul\'e \`a l'IH\'ES et 
\`a l'UMPA de l'ENS-Lyon, et nous voudrions remercier ces deux institutions par leur hospitalit\'e.

\section{Quelques exemples}
\label{exemples}

Tout au long de ce travail, nous ne consid\'ererons que des transformations qui respectent 
l'orientation. En nous inspirant du travail de Denjoy \cite{denjoy}, o\`u la construction 
des contre-exemples $C^1$ pr\'ec\`ede au th\'eor\`eme (et en constitue la partie la 
plus longue), nous ferons une r\'evision rapide des constructions d'exemples de 
diff\'eomorphismes du cercle et de l'intervalle \`a des centralisateurs non 
triviaux. La pr\'esentation est inspir\'ee de \cite{FF,ceuno}, bien que des 
constructions similaires, bas\'ees sur les exemples de Pixton \cite{Pix}, 
apparaissent d\'ej\`a dans le travail de Tsuboi \cite{TsP}.

Rappelons qu'\'etant donn\'e un hom\'eomorphisme $\eta: [0,\infty[ \rightarrow [0,\infty[$, 
on dit qu'une fonction $\psi$ d\'efinie sur le cercle ou un intervalle et \`a valeurs 
r\'eelles est $\eta$-continue s'il existe $C<\infty$ tel que pour tout $x,y$ on ait
$$|\psi(x) - \psi(y)| \leq C \thinspace \eta \big( |x-y| \big).$$
Lorsque $\tau$ appartient \`a $]0,1[$, un diff\'eomorphisme $f$ 
est dit de classe $C^{1+\tau}$ si sa d\'eriv\'ee est $\eta$-continue 
par rapport \`a $\eta(s) = s^{\tau}$, c'est-\`a-dire s'il existe une constante 
$C < \infty$ telle que pour tout $x,y$ on ait $|f'(x) - f'(y)| \leq C |x-y|^{\tau}$ 
(on dit \'egalement que $f'$ est $\tau$-H\"older continue). Nous dirons que 
$f$ est de classe $C^{1+lip}$ si sa d\'eriv\'ee est lipschitzienne. 


\subsection{Les contre-exemples de Denjoy-Herman}
\label{ex1}

Pour construire nos exemples d'actions libres et non minimales 
de $\mathbb{Z}^d$ sur le cercle, nous utiliserons la fa- mille d'applications 
$\varphi_{a,b} \!:\! [0,a] \rightarrow [0,b]$ introduite par Yoccoz et 
donn\'ee sur $]0,a[$ par $\varphi_{a,b}(x) \!=\! \varphi_b \circ (\varphi_a)^{-1}(x)$, 
o\`u $\varphi_a : \mathbb{R} \rightarrow ]0,a[$ est d\'efini par
$$\varphi_a(u) = \frac{1}{\pi} \int_{-\infty}^u
\frac{ds}{s^2+(1/a)^2} = \frac{a}{2} + \frac{a}{\pi} \arctan(au).$$
En faisant  $u = \varphi_a^{-1}(x)$ on voit que
$$\varphi_{a,b}'(x) = \varphi_b'(u) \thinspace 
(\varphi_a^{-1})'(x) = \frac{\varphi_b'(u)}{\varphi_a'(u)} 
= \frac{u^2 + 1/a^2}{u^2 + 1/b^2},$$
ce qui montre que $\varphi_{a,b}$ est un diff\'eomorphisme de classe $C^1$ tangent \`a 
l'identit\'e aux extr\'emit\'es. De plus, pour la deuxi\`eme d\'eriv\'ee on trouve 
ais\'ement la majoration
\begin{equation}
\left| \varphi_{a,b}''(x) \right| \leq \frac{6 \pi}{a} \left| \frac{b}{a} - 1 \right|.
\label{tres}
\end{equation} 
Si pour des intervalles non d\'eg\'en\'er\'es $I=[x_0,y_0]$ et $J=[x_1,y_1]$ on d\'esigne 
par $\varphi(I,J): I \rightarrow J$ le diff\'eomorphisme donn\'e par 
$$\varphi(I,J)(x) = x_1 + \varphi_{y_0-x_0,y_1-x_1}(x-x_0),$$
alors on constate imm\'ediatemment que la famille des $\varphi(I,J)$ est \'equivariante, 
dans le sens que
$$\varphi(J,K) \circ \varphi(I,J) = \varphi(I,K).$$
Gr\^ace aux propri\'et\'es pr\'ec\'edentes, cette famille permet de 
construire des contre-exemples de Denjoy de classe $C^{1+1\!/\!d-\varepsilon}$ 
dont le centralisateur contient un sous-groupe isomorphe \`a $\mathbb{Z}^d$ 
qui agit librement. Plus g\'en\'eralement, nous construirons de tels 
contre-exemples en classe $C^{1+\eta}$ par rapport au module de continuit\'e 
$\eta(s) = s^{1\!/\!d} [\log(1/s)]^{1\!/\!d + \varepsilon}$. Pour cela, on commence avec 
$d$ rotations $R_{\theta_1}, \ldots, R_{\theta_d}$ d'angles lin\'eair\'ement 
ind\'ependants sur les rationnels. On fixe un entier $m \geq 2$ et un point 
$p \in \clo$, et pour chaque $(i_1,\ldots,i_d) \in \mathbb{Z}^d$ on remplace 
le point $p_{i_1,\ldots,i_d} = R_{\theta_1}^{i_1} \cdots R_{\theta_d}^{i_d}(p)$ 
par un intervalle $I_{i_1,\ldots,i_d}$ de longueur  
$$\ell_{(i_1,\ldots,i_d)} =
\frac{1}{\big(|i_1|+\cdots+|i_d|+m\big)^d \thinspace 
\big[\log(|i_1|+\cdots+|i_d|+m)\big]^{1+\varepsilon}}.$$ 
On obtient ainsi un nouveau cercle (de longueur finie 
$T_m \leq 2^d \varepsilon / [\log(m)]^{\varepsilon} (d\!-\!1)!$), sur lequel 
les $R_{\theta_j}$ induisent de mani\`ere unique des hom\'eomorphismes $f_j$ 
v\'erifiant, pour tout $x \in I_{i_1,\ldots,i_j,\ldots,i_d}$,
$$f_j(x) = \varphi(I_{i_1,\ldots,i_j,\ldots,i_d},I_{i_1,\ldots,1+i_j,\ldots,i_d})(x).$$
En vertu des propri\'et\'es d'\'equivariance des $\varphi(I,J)$, ces 
hom\'eomorphismes $f_j$ commutent entre eux. V\'erifions maintenant que $f_1$ est 
de classe $C^{1+\eta}$, le cas des autres $f_i$ \'etant analogue. Pour cela, 
on remarque d'abord que si $x$ et $y$ appartiennent au 
m\^eme intervalle $I_{i_1,\ldots,i_d}$, alors $|f_1'(x)-f_1'(y)| = f_1''(p)|x-y|$ 
pour un point $p \in I_{i_1,\ldots,i_d}$. Ceci donne, gr\^ace \`a (\ref{tres}) 
et au fait que la fonction $s \mapsto s / \eta(s)$ est croissante,
$$\frac{|f_1'(x)-f_1'(y)|}{\eta\big(|x-y|\big)} 
\leq \frac{|x-y|}{\eta\big(|x-y|\big)}\cdot \frac{6\pi}{|I_{i_1,\ldots,i_d}|} 
\left| \frac{|I_{1+i_1,\ldots,i_d}|}{|I_{i_1,\ldots,i_d}|}-1 \right|
\leq \frac{6\pi}{\eta\big(|I_{i_1,\ldots,i_d}|\big)} 
\left| \frac{|I_{1+i_1,\ldots,i_d}|}{|I_{i_1,\ldots,i_d}|}-1 \right|.$$
Or, l'expression \`a droite est \'egale \`a 
\begin{small}
$$\left| \frac{(|i_1|\!\+\!\cdots\!\+\!|i_d|\!\+\!m)^d 
[\log(|i_1|\!\+\!\cdots\!\+\!|i_d|\!\+\!m)]^{1\+\varepsilon}}
{(|1\!\+\!i_1|\!\+\!\cdots\!\+\!|i_d|\!\+\! m)^d 
[\log(|1\!\+\!i_1|\!\+\!\cdots\!\+\!|i_d|\!\+\!m)]^{1\+\varepsilon}} 
- 1  \right|  \frac{6 \pi \thinspace (|i_1|\!\+\!\cdots\!\+\!|i_d|\!\+\!m) 
[\log(|i_1|\!\+\!\cdots\!\+\!|i_d|\!\+\!m)]^{(1\+\varepsilon)\!/\!d}}
{[d\log(|i_1|\!\+\!\cdots\!\+\!|i_d|\!\+\!m) \+ (1\+\varepsilon) 
\log (\log (|i_1|\!\+\!\cdots\!\+\!|i_d|\!\+\!m))]^{1\!/\!d \+ \varepsilon}}.$$
\end{small}En appliquant le th\'eor\`eme des accroissements finis \`a la fonction 
$s \mapsto s^d [\log(s)]^{1+\varepsilon}$, on constate ais\'ement que cette expression 
est majorable par $(d+2) / [\log(m)]^{\varepsilon(1-1\!/\!d)}$. 
Lorsque $x$ et $y$ appartiennent \`a des intervalles $|I_{i_1,\ldots,i_d}|$ distincts, 
la m\^eme majoration modulo un facteur $2$ est valable, car les diff\'eomorphismes 
sont tangents \`a l'identit\'e aux extr\'emit\'es de ces intervalles. Cela montre 
qu'en faisant tendre $m$ vers l'infini, les diff\'eomorphismes $f_j$ deviennent 
aussi proches que l'on veut de l'identit\'e en topologie $C^{1 + \eta}$. Puisque 
$T_m \leq 1$ pour tout $m \in \mathbb{N}$ assez grand, ceci reste encore valable apr\`es 
renormalisation du cercle (de fa\c{c}on \`a ce que sa longueur soit \'egale \`a $1$).

\vspace{0.15cm}

Les exemples que nous venons de construire v\'erifient $f_j'(x) = 1$ pour tout 
point $x$ qui est une extr\'emit\'e de l'un des intervalles $I_{i_1,\ldots,i_d}$. 
Nous verrons ci-dessous que, sous une telle hypoth\`ese, il est impossible de 
fabriquer des exemples analogues de classe $C^{1+1\!/\!d}$. L'argument de la 
preuve pr\'esent\'ee \`a continuation est \`a rapprocher avec \cite{norton}.

\vspace{0.1cm}

\begin{prop} {\em Consid\'erons une action de $\mathbb{Z}^d$ par diff\'eomorphismes du cercle 
qui est libre et semiconjugu\'ee \`a une action par des rotations sans y \^etre conjugu\'ee. 
Supposons que pour toute composante connexe $I$ du compl\'ementaire du Cantor invariant 
minimal, et pour tout \'el\'ement $h \in \mathbb{Z}^d$, la d\'eriv\'ee de (l'image de) 
$h$ aux extr\'emit\'es de $I$ soit \'egale \`a $1$. Alors la r\'egularit\'e de 
l'action est strictement inf\'erieure \`a $C^{1+1\!/\!d}$.}
\label{tangente}
\end{prop}

\noindent{\bf D\'emonstration.} Supposons le contraire et soient $f_1,\ldots,f_d$ les 
g\'en\'erateurs de $\mathbb{Z}^d$ en tant que sous-groupe de $\mathrm{Diff}_+^{1+1\!/\!d}(\clo)$. 
Notons $C_i$ la constante de H\"older d'exposant $1/d$ pour $f_i'$, et posons 
$C = \max \{C_1,\ldots,C_d\}$. Fixons une composante connexe $I$ du 
compl\'ementaire du Cantor invariant, et pour chaque entier $n \geq 0$ 
notons $\ell_n$ la longueur minimale d'un intervalle de la forme 
$f_{1}^{i_1} \cdots f_{d}^{i_d}(I)$, avec $i_j \geq 0$ et $\sum_{j} i_j = n$. 
Sans perdre en g\'en\'eralit\'e, nous pouvons supposer que tous ces 
intervalles ont une longueur inf\'erieure ou \'egale \`a 
$(\frac{1}{C(1+1\!/\!d)})^{1\!/\!d}$. Nous affirmons que 
\begin{equation}
\ell_{n+1} \geq \ell_n (1 - C \thinspace \ell_n^{1\!/\!d}).
\label{eles}
\end{equation}
En effet, si $\ell_{n+1}$ est r\'ealis\'e comme la longueur correspondante \`a un 
intervalle $J$, alors $J$ est l'image par l'un des g\'en\'erateurs $f_i$ d'un intervalle 
$K = [a,b]$ dont la longueur est par d\'efinition sup\'erieure ou \'egale \`a $\ell_n$. 
Puisqu'il existe $q \in K$ tel que $f_i' (q) = |J|/|K|$ et que $f_i'(a)=f_i'(b)=1$, on a 
$$\left| \frac{|J|}{|K|} -1 \right| \leq C_i |q - a|^{1\!/\!d} \leq C |K|^{1\!/\!d},$$
d'o\`u on obtient 
$$\big| \ell_{n+1} - |K| \big| = \big| |J|-|K| \big| \leq C|K|^{1+1\!/\!d}.$$
Donc 
$$\ell_{n+1} \geq |K| \thinspace \big( 1 - C |K|^{1\!/\!d} \big),$$
et ceci implique (\ref{eles}), car $|K| \geq \ell_n$ et la fonction 
$s \mapsto s (1 - C s^{1\!/\!d})$ est croissante sur l'intervalle 
$\big[ 0,(\frac{1}{C(1+1\!/\!d)})^{1\!/\!d} \big]$.

Nous affirmons maintenant que pour $A = \min \{ \ell_1, d^d / 2^{d^2} C^d \}$ 
et pour tout $n\! \in\! \mathbb{N}$,
\begin{equation}
\ell_n \geq \frac{A}{n^{d}}.
\label{pata2}
\end{equation}
Cette in\'egalit\'e est v\'erifiable ais\'ement par r\'ecurrence. Pour $n=1$ elle 
a lieu \`a cause de la condition $A \leq \ell_1$. D'autre part, si elle est valable 
pour un entier $n \geq 1$ alors, la fonction $s \mapsto s (1 - C s^{1\!/\!d})$ \'etant 
croissante sur l'intervalle $\big[ 0,(\frac{1}{C(1+1\!/\!d)})^{1\!/\!d} \big]$, 
\`a partir de l'in\'egalit\'e (\ref{eles}) et de la condition 
$A^{1\!/\!d} \leq d / 2^{d} C$ on obtient 
$$\ell_{n+1} \geq \frac{A}{n^d} \left( 1 - C \thinspace \frac{A^{1\!/\!d}}{n} \right) 
\geq \frac{A}{(n+1)^d}.$$

Nous sommes maintenant en mesure d'achever la preuve de la proposition. 
Pour cela, remarquons que la quantit\'e d'intervalles de la forme 
$I_{i_1,\ldots,i_d} = f_1^{i_1} \cdots f_d^{i_d}(I)$, avec $i_j \geq 0$ et 
$\sum_j i_j = n$, est sup\'erieure ou \'egale \`a $B n^{d-1}$ pour une certaine 
constante universelle $B > 0$. Ceci donne, gr\^ace \`a (\ref{pata2}), 
$$\sum_{(i_1,\ldots,i_d) \in \mathbb{N}_0^d} |I_{i_1,\ldots,i_d}| = 
|I| + \sum_{n \geq 1} \sum_{\sum {i_j}=n} |I_{i_1,\ldots,i_d}| \geq 
B \sum_{n \geq 1} n^{d-1} \ell_n \geq AB \sum_{n \geq 1} \frac{1}{n} = \infty,$$
ce qui est absurde, car la somme des longueurs des intervalles (disjoints) 
$I_{i_1,\ldots,i_d}$ est finie. \esp $\square$


\subsection{Les contre-exemples de Pixton-Tsuboi}
\label{ex2}

Dans le contexte des groupes ab\'eliens de diff\'eomorphismes 
de l'intervalle on dispose du fameux lemme de Kopell \cite{Ko}, dont la 
version la plus connue n'est d'habitude \'enonc\'ee qu'en classe $C^2$. Pour 
la commodit\'e du lecteur, nous donnons ci-dessous la version g\'en\'erale 
avec une preuve simple dont l'id\'ee sous-jacente sera utilis\'ee \`a 
plusieurs reprises. Rappelons qu'un diff\'eomorphisme $f$ de $[0,1[$ est 
de classe $C^{1+vb}$ si sa d\'eriv\'ee est \`a variation born\'ee sur tout 
intervalle compact contenu dans $[0,1[$.

\vspace{0.25cm}

\noindent{\bf Lemme [Kopell].} 
{\em Soient $f$ et $g$ deux diff\'eomorphismes de l'intervalle $[0,1[$ 
qui commutent entre eux. Supposons que $f$ soit de classe $C^{1+vb}$ et 
$g$ de classe $C^1$. Si $f$ n'a pas de point fixe sur $]0,1[$ et $g$ 
poss\`ede de tels points, alors $g$ est l'identit\'e.}

\vspace{0.15cm}

\noindent{\bf D\'emonstration.} Quitte \`a changer $f$ par son inverse, 
on peut supposer que $f(x)<x$ pour tout $x \!\in ]0,1[$. Soit $b\!\in ]0,1[$ 
l'un des points fixes de $g$. Pour chaque $n \in \mathbb{Z}$ notons 
$b_n = f^n(b)$, et soit $a = b_1 = f(b)$\footnote{Remarquons que pour $f_1 = f, 
f_2 = g$ et $I_n = f^n(]a,b[)$, nous sommes exactement sous une hypoth\`ese 
combinatoire du type (\ref{comb}), avec $d = 1$.}. 
Puisque $g$ fixe l'intervalle $[b_{n+1},b_n]$, 
pour chaque $n \in \mathbb{N}$ il existe $c_n \in [b_n,b_{n+1}]$ tel que $g'(c_n) = 1$. 
\'Etant donn\'e que la suite $(c_n)$ tend vers l'origine et que $g$ est de classe 
$C^1$, on a n\'ecessairement $g'(0)=1$. Notons $M = M(f)$ la variation 
du logarithme de la d\'eriv\'ee de $f$ sur l'intervalle $[0,b]$. Si $u$ 
et $v$ appartiennent \`a $[a,b]$ alors
\begin{equation}
\left| \log \Big( \frac{(f^n)'(v)}{(f^n)'(u)} \Big) \right| \leq
\sum\limits_{i=1}^{n} \big| \log \big( f'(f^{i-1}(v)) \big) - 
\log \big(f' (f^{i-1}(u)) \big) \big| \leq M. 
\label{clasiquita}
\end{equation}
En posant $u = x \in [a,b]$ et $v = f^{-n} g f^n (x) = g(x) \in [a,b]$, 
en utilisant l'\'egalit\'e
$$g'(x) = \frac{(f^{n})'(x)}{(f^{n})'(f^{-n} g f^n(x))} 
\hspace{0.1cm} g'(f^n(x)) = \frac{(f^{n})'(x)}{(f^{n})'(g(x))} 
\hspace{0.1cm} g'(f^n(x)),$$
et en passant \`a la limite lorsque $n$ tend vers l'infini, on obtient l'in\'egalit\'e 
$\sup_{x \in [a,b]} g'(x) \leq \exp(M)$. Or, ceci reste valable lorsqu'on remplace 
$g$ par $g^j$ pour n'importe quel $j \in \mathbb{N}$ (car $M$ ne d\'epend que de 
$f$). On en d\'eduit que $\sup_{x \in [a,b]} (g^j)'(x) \leq \exp(M)$. 
Comme $g$ fixe $a$ et $b$, ceci implique que la restriction de $g$ \`a 
l'intervalle $[a,b]$ est l'identit\'e. Finalement, en conjugant successivement par 
$f$, on conclut que $g$ est l'identit\'e sur tout l'intervalle $[0,1[$. \esp $\square$

\vspace{0.1cm}

\begin{rem} Notons que dans la d\'emonstration pr\'ec\'edente, tout le contr\^ole de 
distorsion est r\'ealis\'e par le diff\'eomorphisme $f$, qui doit donc \^etre de 
classe au moins $C^{1+vb}$, tandis que l'autre diff\'eomorphisme $g$ n'a besoin 
que d'\^etre de classe $C^1$. Le lemme n'est cependant plus valable 
lorsque $g$ est seulement lipschitzien. En effet, si l'on fixe 
n'importe quel hom\'eomorphisme lipschitzien $h$ de $[a,b]$ et on 
l'\'etend \`a $[0,1[$ en commutant avec $f$, alors pour tout $c \in ]0,1[$ 
l'application qu'on obtient est encore lipschitzienne sur $[0,c]$, 
et sa constante de Lipschitz y diff\`ere de celle de $h$ par un 
facteur au plus \'egal \`a $\exp(M_c)$, o\`u $M_c$ d\'esigne la variation 
du logarithme de la d\'eriv\'ee de $f$ sur $[0,c]$. L'inter\^et de cette observation vient 
du fait que, dans \cite{GT}, il a \'et\'e remarqu\'e que si $f:[0,1[ \rightarrow [0,1[$ est 
un diff\'eomorphisme de classe $C^2$ sans point fixe \`a l'int\'erieur et $g$ est un 
hom\'eomorphisme de $[0,1[$ qui commute avec $f$, alors la restriction de $g$ \`a 
$]0,1[$ est automatiquement de classe $C^2$ lorsque $g$ est un diff\'eomorphisme 
de classe $C^1$. De plus, si $f$ est hyperbolique, {\em i.e.} si $f'(0) \neq 1$, 
alors $g$ est de classe $C^2$ sur $[0,1[$ d\`es qu'il est 
diff\'erentiable \`a l'origine \cite{St}.
\end{rem}

On peut donner des contre-exemples au lemme de Kopell analogues
\`a ceux du th\'eor\`eme de Denjoy, tout en respectant la propri\'et\'e 
combinatoire (\ref{comb}).  Pour cela, fixons un entier $m \geq 2$, 
et pour $(i_1,\ldots,i_d) \in \mathbb{Z}^d$ posons une nouvelle fois
$$\ell_{(i_1,\ldots,i_d)} = \frac{1}{\big(|i_1|+\cdots+|i_d|+m\big)^d \thinspace 
\big[\log(|i_1|+\cdots+|i_d|+m)\big]^{1+\varepsilon}}.$$
Par r\'ecurrence sur \thinspace $j$ \thinspace d\'efinisons \thinspace 
$\ell_{(i_1,\ldots,i_{j-1})} = \sum_{i_j \in \mathbb{Z}} 
\ell_{(i_1,\ldots,i_{j-1},i_j)}$. \thinspace Notons \thinspace 
$[x_{(i_1,\ldots,i_j,\ldots,i_d)},y_{(i_1,\ldots,i_j,\ldots,i_d)}]$ 
\thinspace l'intervalle d'extr\'emit\'es 
$$x_{(i_1,\ldots,i_j,\ldots,i_d)} = \sum_{i_1' < i_1} \ell_{(i_1')} + 
\sum_{i_2' < i_2} \ell_{(i_1,i_2')} + \ldots + 
\sum_{i_d' < i_d} \ell_{(i_1,\ldots,i_{d-1},i_d')}, 
\qquad y_{(i_1,\ldots,i_j,\ldots,i_d)} = 
x_{(i_1,\ldots,i_j,\ldots,i_d)} + \ell_{(i_1,\ldots,i_d)}.$$
Consid\'erons le diff\'eomorphisme $f_j$ de $[0,1]$ dont la restriction aux intervalles 
$[x_{(i_1,\ldots,i_j,\ldots,i_d)},y_{(i_1,\ldots,i_j,\ldots,i_d)}]$ co\"{\i}ncide avec 
$$\varphi \big( [x_{(i_1,\ldots,i_j,\ldots,i_d)},y_{(i_1,\ldots,i_j,\ldots,i_d)}], 
[x_{(i_1,\ldots,i_j-1,\ldots,i_d)},y_{(i_1,\ldots,i_j-1,\ldots,i_d)}] \big).$$ 
Les estim\'ees qui pr\'ec\`edent la proposition (\ref{tangente}) permettent 
de d\'emontrer une nouvelle fois que les $f_i$ ainsi obtenus sont de classe 
$C^{1+\eta}$ par rapport au module de continuit\'e 
$\eta(s) = s^{1\!/\!d} [\log(1/s)]^{1\!/\!d+\varepsilon}$. 
On obtient donc des contre-exemples au lemme de Kopell donn\'es par des actions de 
$\mathbb{Z}^d$ par diff\'eomorphismes de classe $C^{1+1\!/\!d-\varepsilon}$ de l'intervalle. 
Cependant, nous devons imp\'erativement souligner que la r\'egularit\'e atteinte 
par cette m\'ethode n'est pas optimale. Pour aboutir \`a la classe optimale 
$C^{1+1\!/\!(d-1)-\varepsilon}$, il est n\'ec\'essaire de modifier la construction 
pr\'ec\'edente afin de suprimer les tangences \`a l'identit\'e 
\begin{tiny}$^{_\ll}$\end{tiny}excessives\hspace{0.05cm}\begin{tiny}$^{_\gg}$\end{tiny}. 
En effet, une modification simple de la proposition \ref{tangente} montre que 
par cette m\'ethode on n'atteindra m\^eme pas la classe $C^{1+1\!/\!d}$. Le 
lecteur int\'eress\'e trouvera dans \cite{TsP} la construction des 
contre-exemples optimaux au lemme de Kopell ; nous reviendrons sur 
ce point dans le \S \ref{kop}.

Pour finir ce paragraphe, nous voudrions pr\'esenter une preuve simple d'une version du 
lemme de Kopell pour des actions de ${\mathbb Z}^3$ avec une hypoth\`ese de r\'egularit\'e 
interm\'ediaire (qui n\'eanmoins n'est pas optimale). Pour cela, commen\c{c}ons par un lemme 
\'el\'ementaire mais tr\`es utile.

\begin{lem} {\em Soit $g: [a,b] \rightarrow [a,b]$ un diff\'eomorphisme 
de classe $C^{1+\tau}$. Si $C = C(g)$ d\'esigne la constante de 
$\tau$-h\"olderianit\'e de $g'$, alors pour tout $x \!\in ]a,b[$ on a}
\begin{equation}
|x - g(x)| \leq C \thinspace |b - a|^{1+\tau}.
\label{grosero}
\end{equation}
\end{lem}

\noindent{\bf D\'emonstration.} Fixons un point arbitraire $x$ dans $]a,b[$, 
et prenons deux points $p \in [a,x]$ et $q \in [a,b]$ tels que 
\thinspace $g'(p) = \big( g(x) - g(a) \big) / (x-a)$ \thinspace 
and \thinspace $g'(q) = \big( g(b)-g(a) \big)/(b-a) = 1$. 
\thinspace Alors on a
\begin{equation}
|x - g(x)| = \big| x - \big[ (x-a)g'(p) + g(a) \big] \big| = 
\big| (x-a) \big( 1 - g'(p) \big) \big|.
\label{unila}
\end{equation}
Puisque $g'$ est $\tau$-H\"older continue,
\begin{equation}
|1-g'(p)| = |g'(q) - g'(p)| \leq C |q-p|^{\tau} \leq C |b-a|^{\tau}.
\label{dorila}
\end{equation}
De (\ref{unila}) et (\ref{dorila}) on d\'eduit 
$$|x - g(x)| \leq |x-a| \cdot C|b-a|^{\tau} \leq C|b-a|^{1+\tau},$$
ce qui conclut la preuve du lemme. \esp $\square$

\vspace{0.25cm}

La proposition suivante n'est qu'une cons\'equence du lemme de Kopell g\'en\'eralis\'e 
(th\'eor\`eme B). Cependant, l'id\'ee de la preuve que nous pr\'esentons, dans un esprit 
dynamique \begin{tiny}$^{_\ll}$\end{tiny}classique\hspace{0.05cm}\begin{tiny}$^{_\gg}$\end{tiny}, 
pourrait \^etre utile dans d'autres situations.

\vspace{0.1cm}

\begin{prop} {\em Soient $f$, $g$ et $h$ des diff\'eomorphismes de l'intervalle $[0,1]$ 
qui commutent entre eux, avec $f$ et $g$ de classe $C^{1+\tau}$ et $h$ de classe $C^1$. 
Supposons qu'il existe des intervalles ouverts disjoints de la forme $I_{i,j}$ qui 
soient dispos\'es en respectant l'ordre lexicografique sur $]0,1[$, et de sorte que 
pour tout $(i,j) \in \mathbb{Z}^2$ on ait}
$$f(I_{i,j}) = I_{i-1,j}, \qquad g(I_{i,j}) = I_{i,j-1} \qquad 
h(I_{i,j}) = I_{i,j}.$$
{\em Si $\tau \geq (\sqrt{5}-1)/2$ alors la restriction 
de $h$ \`a la r\'eunion des $I_{i,j}$ est l'identit\'e.}
\label{cano}
\end{prop}

\noindent{\bf D\'emonstration.} Reprenons l'argument de la preuve du lemme de Kopell. 
Pour tout $x \in I = I_{0,0}$ et tout $(n,m) \in \mathbb{N}^2$ on a 
$$\log(h^n)'(x) = \log \big((h^n)'(f^m(x)) \big) + \log \big( (f^m)'(x) \big) 
- \log \big( (f^m)'(h^n(x)) \big).$$
Puisque $f^m(x)$ converge vers un points fixe (parabolique) de $h$, 
en prenant la limite lorsque $m$ tend vers l'infini on obtient
\begin{eqnarray*}
|\log \big( (h^n)'(x) \big) | &\leq& \sum_{k \geq 0} \big|\log \big( f'(f^k(x)) \big) 
- \log \big( f' (h^n(f^k(x)))\big) \big|\\
&\leq& \bar{C} \sum_{k \geq 0} \big| f^k(x) - h^n(f^k(x)) \big|^{\tau}\\
&\leq& \bar{C} \sum_{k \geq 0} |f^k(I)|^{\tau},
\end{eqnarray*}
o\`u $\bar{C} = \bar{C}(f)$ d\'esigne la constante de $\tau$-h\"olderianit\'e de $\log(f')$. 
Pour obtenir la convergence de la derni\`ere s\'erie, l'id\'ee consiste \`a remarquer que 
$|f^k(I)|$ est tr\`es petit par rapport \`a $|f^k(J)|$, o\`u $J$ d\'esigne le plus petit 
intervalle contenant tous les $g^m(I)$, avec $m \in \mathbb{Z}$ (notons que $g$ fixe 
l'intervalle $J$). En effet, si $C=C(g)$ est la constante de $\tau$-h\"olderianit\'e 
de $g'$, alors l'in\'egalit\'e (\ref{grosero}) donne \thinspace 
$|f^k(I)| \leq C \thinspace |f^k(J)|^{1+\tau}$. \thinspace Par suite, 
$$\big| \log \big( (h^n)'(x) \big) \big| \leq 
\bar{C} C \sum_{k \geq 0} |f^k(J)|^{\tau(1+\tau)}.$$
Si \thinspace $\tau \geq (\sqrt{5}-1)/2$ \thinspace alors \thinspace 
$\tau(1+\tau) \geq 1$, \thinspace et l'expression pr\'ec\'edente est 
major\'ee par \esp $\bar{C} C \sum_{k \geq 0} |f^k(J)| \leq \bar{C} C$. 
\esp Autrement dit, pour tout $x \in I$ et tout $n \in \mathbb{N}$ on 
a l'in\'egalit\'e \thinspace $(h^n)'(x) \leq \exp(\bar{C}C)$. \thinspace Comme dans la 
fin de la preuve du lemme (classique) de Kopell, ceci implique que la restriction 
de $h$ \`a $I$ (et donc \`a tous les $I_{i,j}$) est l'identit\'e. \esp $\square$

\vspace{0.28cm}

Si pour chaque entier $d \geq 3$ on note $\tau_d$ l'unique r\'eel positif qui v\'erifie 
l'\'egalit\'e $\tau_d (1+\tau_d)^{d-2} = 1$, alors la m\'ethode pr\'ec\'edente permet 
de d\'emontrer une proposition analogue pour des actions de $\mathbb{Z}^{d}$ par  
diff\'eomorphismes de classe $C^{1+\tau_d}$ (et qui v\'erifient une condition 
combinatoire du type (\ref{comb})). Or, cet exposant $\tau_d$ n'est pas du 
tout optimal. En effet, bien que la suite $(\tau_d)$ converge vers z\'ero lorsque 
$d$ tend vers l'infini, le nombre $\tau_d$ est sup\'erieur \`a $1/(d-1)$ (par exemple, 
$\tau_3$ est \'egal au nombre d'or $(\sqrt{5}-1)/2 > 1/2$).


\section{Le th\'eor\`eme de Denjoy g\'en\'eralis\'e}
\label{denj}

\subsection{Le principe g\'en\'eral}

Nous rappelons dans la suite l'id\'ee de la preuve du th\'eor\`eme de Denjoy donn\'ee 
par Schwartz dans \cite{Sc}. C'est un principe qui est devenu classique gr\^ace  \`a 
la formulation et les applications aux feuilletages de codimension 1 donn\'ees par 
Sacksteder dans \cite{Sa} (voir \cite{CClibro} pour une discussion plus d\'etaill\'ee).

\vspace{0.05cm}

\begin{lem} {\em Soit $\Gamma$ un groupe de diff\'eomorphismes de classe 
$C^{1+lip}$ d'une vari\'et\'e unidimensionnelle compacte. Supposons qu'il existe 
un sous-ensemble fini $\mathcal{G}$ de $\Gamma$ et un intervalle ouvert $I$ 
satisfaisant la propri\'et\'e suivante : pour chaque $n \in \mathbb{N}$ il 
existe un \'el\'ement $h_n = g_{i_{n,n}} g_{i_{n-1,n}} \cdots g_{i_{1,n}}$ dans 
$\Gamma$ tel que tous les $g_{i_{m,n}}$ appartiennent \`a $\mathcal{G}$, 
les intervalles $I, g_{i_{1,n}}(I), g_{i_{2,n}} g_{i_{1,n}}(I), 
\ldots, g_{i_{n-1,n}} \cdots g_{i_{1,n}}(I)$ sont deux \`a deux disjoints, 
et $h_n(I)$ s'accumule sur l'une des extr\'emit\'es de $I$. 
Alors pour $n \in \mathbb{N}$ assez grand, l'application $h_n$ poss\`ede un point 
fixe hyperbolique (qui est proche de l'extr\'emit\'e correspondante de $I$).}
\label{!}
\end{lem}

\vspace{0.05cm}

La d\'emonstration de ce lemme utilise le fait que les intervalles $I,\ldots, 
g_{i_{n-1,n}} \cdots g_{i_{1,n}}(I)$ sont disjoints, ainsi que l'hypoth\`ese de 
r\'egularit\'e $C^{1+lip}$, pour contr\^oler la distorsion de $h_n$ sur $I$ gr\^ace 
\`a l'argument classique de Denjoy. Ensuite, l'id\'ee consiste \`a contr\^oler la 
distorsion de $h_n$ sur un intervalle plus large $J$ contenant $I$ et ind\'ependant 
de $n \in \mathbb{N}$. On aboutit \`a ceci par un argument de r\'ecurrence assez 
subtil (remarquons que, en g\'en\'eral, 
les intervalles $J,\ldots, g_{i_{n-1,n}} \cdots g_{i_{1,n}}(J)$ 
ne sont pas disjoints!). La contraction (topologique) devient alors 
\'evidente (voir la figure 1), et cette contraction doit \^etre hyperbolique 
\`a cause du contr\^ole (uniforme) de la distorsion.

\vspace{0.2cm}

\beginpicture

\setcoordinatesystem units <1cm,1cm>

\putrule from -5.5 0 to 5.5 0 

\putrule from -2.8 -0.015 to 2.8 -0.015

\putrule from -2.8 0.015 to 2.8 0.015

\putrule from 3 -0.015 to 3.5 -0.015 

\putrule from 3 0.015 to 3.5 0.015

\circulararc -45 degrees from 0 0.5 center at 1.6 -3.3

\circulararc 132 degrees from 3.8 -0.5 center at 4.3 -0.3 

\circulararc -55 degrees from 2.2 -0.5 center at -1.25 6

\plot 3.8 -0.5 
3.9 -0.56 /

\plot 3.8 -0.5 
3.8 -0.615 /

\plot 2.2 -0.5 
2.12 -0.63 /

\plot 2.2 -0.5 
2.03 -0.51 /

\plot 3.14 0.52 
3 0.52 /

\plot 3.14 0.52 
3.05 0.62 /

\put{$\Big|$} at -4.8 0  
\put{$\Big|$} at 4.8 0 
\put{$|$} at 3.722 0  
\put{$|$} at 2.5 0
\put{$\Big($} at -2.8 0
\put{$\Big)$} at 2.8 0 
\put{$($} at 3 0
\put{$)$} at 3.5 0 
\put{$h_n$} at -3.8 -1.2
\put{$h_n$} at 4.4 -1.2
\put{$h_n$} at 1.3 1.1
\put{$a$} at -5.051 0.23 
\put{$b$} at 5.051 0.23

\plot 3.615 1 
3.615 0.3 /

\plot 3.615 0.3
3.55 0.45 /

\plot 3.615 0.3 
3.68 0.45 /

\put{$J \!=\! [a,b]$} at -3.8 1 
\put{Figure 1} at 0 -1.8 
\put{} at -8.1 0

\small
\put{$h_n(I)$} at 3.25 -0.6
\put{$I$} at 0 -0.6 
\put{$\bullet$} at 3.62 0 
\put{point fixe} at 3.615 1.6
\put{hyperbolique} at 3.615 1.25

\endpicture


\vspace{0.4cm}

Si l'on veut contr\^oler les distorsions en classe $C^{1+\tau}$, 
on est amen\'e \`a estimer des sommes du type
\begin{equation}
|I|^{\tau} + \sum_{k=1}^{n-1} |g_{i_{k,n}} \dots g_{i_{1,n}}(I)|^{\tau}
\label{ento}
\end{equation}
Or, m\^eme si les intervalles $I,\ldots, g_{i_{n-1,n}} \cdots g_{i_{1,n}}(I)$ 
sont deux \`a deux disjoints, cette somme peut devenir trop grande avec 
$n$. Pour r\'esoudre ce probl\`eme, notre id\'ee consiste \`a 
penser les compositions de mani\`ere al\'eatoire (contrairement 
au cas de la preuve de Schwartz et Sacksteder, o\`u la suite 
des compositions est fix\'ee de mani\`ere d\'eterministe). Plus 
pr\'ecis\'ement, nous consid\'ererons des suites d'applications $h_n$ 
de la forme $h_n = g_n h_{n-1}$ de fa\c{c}on \`a ce que les intervalles 
$\thinspace h_0(I)\!=\!I, h_1(I),\ldots \thinspace$ soient deux 
\`a deux disjoints, et par des arguments d'ordre probabiliste 
nous chercherons \`a montrer que pour 
un \begin{tiny}$^{_\ll}$\end{tiny}chemin 
typique\hspace{0.05cm}\begin{tiny}$^{_\gg}$\end{tiny} la valeur 
de la somme (\ref{ento}) est uniform\'ement born\'ee. Cela nous 
permettra de trouver des \'el\'ements avec des points fixes hyperboliques 
gr\^ace au lemme g\'en\'eral suivant.

\vspace{0.15cm}

\begin{lem} {\em Soit $\Gamma$ un pseudo-groupe de diff\'eomorphismes 
de classe $C^{1+\tau}$ d'une vari\'et\'e unidimensionnelle compacte. 
Supposons qu'il existe un sous-ensemble fini $\mathcal{G}$ de $\Gamma$, 
un intervalle ouvert $I$ et une constante $M\!<\!\infty$ tels que la 
propri\'et\'e suivante soit satisfaite : \`a chaque \'el\'ement 
$g \in \mathcal{G}$ on peut associer un intervalle compact 
$\mathrm{C}_g$ contenu dans un domaine ouvert de d\'efinition de $g$ 
de sorte que, pour tout $n \in \mathbb{N}$, il existe un 
\'el\'ement $h_n = g_{i_{n,n}} \cdots g_{i_{1,n}}$ dans $\Gamma$ 
v\'erifiant les propri\'et\'es suivantes :

\vspace{0.1cm}

\noindent{-- \thinspace tous les $g_{i_{m,n}}$ appartiennent \`a $\mathcal{G}$ ;}

\vspace{0.1cm}

\noindent{-- \thinspace si $g_{i_{k,n}} = g$ alors l'intervalle 
$g_{i_{k-1,n}} \cdots g_{i_{1,n}}(I)$ est contenu dans $\mathrm{C}_g$ \esp (o\`u nous 
convenons que \esp $g_{i_{k-1,n}} \cdots g_{i_{1,n}}$ est l'identit\'e lorsque $k=1$) ;}

\vspace{0.1cm}
  
\noindent{-- \thinspace on a l'in\'egalit\'e}
$$\sum_{k = 0}^{n-1} |g_{i_{k,n}} \cdots g_{i_{1,n}}(I)|^{\tau} \leq M.$$
Alors il existe une constante strictement positive $L = L(\tau,M,|I|;\mathcal{G})$ telle 
que si $n \in \mathbb{N}$ est tel que $h_n(I)$ est contenu dans un $L$-voisinage de 
l'intervalle $I$, alors l'application $h_n$ poss\`ede un point fixe hyperbolique 
(qui est proche de l'extr\'emit\'e correspondante de $I$).}
\label{?}
\end{lem}

\noindent{\bf D\'emonstration.} La preuve du lemme \'etant bien connue, 
nous ne la r\'ep\'etons que pour le cas d'un groupe de diff\'eomorphismes. 
Fixons une constante $C > 0$ telle que pour tout $g \in \mathcal{G}$ et tout 
$x,y$ on ait
$$\big| \log(g'(x)) - \log(g'(y)) \big| \leq C \esp |x-y|^{\tau}.$$
Nous montrerons alors que pour
$$L = L(\tau,M,|I|;\mathcal{G}) = \frac{|I|}{2 \exp(2^{\tau} C M )},$$ 
l'affirmation du lemme est satisfaite.

D\'esignons par $J$ le $2L$-voisinage de $I$, et notons $I'$ (resp. $I''$) 
la composante connexe de $J \setminus I$ \`a droite (resp. \`a gauche) de 
$I$. Pour $n \in \mathbb{N}$ fix\'e nous montrerons par r\'ecurrence sur 
$k \!\in\! \{ 0,\ldots,n \}$ que les deux propri\'et\'es suivantes sont 
simultan\'ement v\'erifi\'ees :

\vspace{0.1cm}

\noindent{$(i)_k \hspace{0.25cm} |g_{i_{k,n}} \cdots g_{i_{1,n}}(I')| \leq  
|g_{i_{k,n}} \cdots g_{i_{1,n}}(I)|$ ;}

\vspace{0.1cm}

\noindent{$(ii)_k \hspace{0.25cm} \sup_{x,y \in I \cup I'} 
\frac{(g_{i_{k,n}} \cdots g_{i_{1,n}})'(x)}{(g_{i_{k,n}} \cdots g_{i_{1,n}})'(y)} 
\leq \exp(2^{\tau} \thinspace C M )$.}

\vspace{0.1cm}

La condition $(ii)_0$ est trivialement v\'erifi\'ee, tandis que $(i)_0$ 
est satisfaite gr\^ace \`a l'hypoth\`ese $|I'| \!=\! 2L \!\leq\! |I|$. 
Supposons que $(i)_j$ et $(ii)_j$ soient valables pour tout $j \in \{0,\ldots,k-1\}$. 
Dans ce cas, pour tout $x,y$ dans $I \cup I'$  nous avons
\begin{eqnarray*}
\left|\log\left(\frac{(g_{i_{k,n}} \cdots g_{i_{1,n}})'(x)}
{(g_{i_{k,n}} \cdots g_{i_{1,n}})'(y)}\right)\right| 
&\leq& \sum_{j=0}^{k-1} \big| \log(g_{i_{j+1,n}}'(g_{i_{j,n}} \cdots g_{i_{1,n}}(x))) - 
\log(g_{i_{j+1,n}}'(g_{i_{j,n}} \cdots g_{i_{1,n}}(y))) \big|\\ 
&\leq& C \sum_{j=0}^{k-1} \big| g_{i_{j,n}} \cdots g_{i_{1,n}}(x) - 
g_{i_{j,n}} \cdots g_{i_{1,n}}(y) \big|^{\tau}\\
&\leq& C \sum_{j=0}^{k-1} \big( |g_{i_{j,n}} \cdots g_{i_{1,n}}(I)| + 
|g_{i_{j,n}} \cdots g_{i_{1,n}}(I')| \big)^{\tau}\\ 
&\leq& C \thinspace 2^{\tau} M.
\end{eqnarray*}
Ceci montre $(ii)_k$. Quant \`a $(i)_k$, remarquons 
qu'il existe $x \in I$ et $y \in I'$ tels que
$$|g_{i_{k,n}} \cdots g_{i_{1,n}} (I)| = |I| \cdot (g_{i_{k,n}} \cdots g_{i_{1,n}})'(x) 
\qquad \mbox { et } \qquad |g_{i_{k,n}} \cdots g_{i_{1,n}}(I')| 
= |I'| \cdot (g_{i_{k,n}} \cdots g_{i_{1,n}})'(y).$$
Donc, par $(ii)_k$,
$$\frac{|g_{i_{k,n}} \cdots g_{i_{1,n}}(I')|}{|g_{i_{k,n}} \cdots g_{i_{1,n}}(I)|}= 
\frac{(g_{i_{k,n}} \cdots g_{i_{1,n}})'(x)}{(g_{i_{k,n}} \cdots g_{i_{1,n}})'(y)} \cdot 
\frac{|I'|}{|I|} \leq \exp(2^{\tau} C M) \frac{|I'|}{|I|} \leq 1,$$
ce qui montre $(i)_k$. Bien s\^ur, un argument analogue 
montre que $(i)_k$ et $(ii)_k$ sont v\'erifi\'ees pour tout 
$k$ dans $\{0,\ldots,n\}$ lorsqu'on remplace $I'$ par $I''$.

Supposons maintenant que $h_n(I)$ soit contenu dans le $L$-voisinage de 
l'intervalle $I$. La propri\'et\'e $(i)_n$ donne alors $h_n(J) \subset J$. De 
plus, si $h_n(I) \subset J$ se trouve \`a droite (resp. \`a gauche) de $I$, alors 
$h_n(I \cup I') \subset I'$ (resp. $h_n(I'' \cup I) \subset I''$). Les deux 
cas \'etant analogues, consid\'erons seulement le premier. Bien \'evidemment, 
$h_n$ poss\`ede au moins un point fixe $x$ dans $I'$. Il nous reste 
donc \`a v\'erifier qu'il s'agit d'un point fixe hyperbolique 
contractant. Or, il existe $y \in I$ tel que 
$$h_n'(y) = \frac{|h_n(I)|}{|I|} \leq \frac{L/2}{|I|}.$$
Par cons\'equent, si $h_n'(x) \geq 1$ alors on aurait \esp
$h_n'(x) / h_n'(y) \geq  2 \thinspace |I| / L$, \esp
et donc, d'apr\`es $(ii)_n$,
$$\exp(2^{\tau} C M) \geq \frac{2 \thinspace |I|}{L},$$
ce qui contredirait la d\'efinition de $L$. \esp $\square$


\subsection{Preuve du th\'eor\`eme de Denjoy g\'en\'eralis\'e}
\label{denjoy-gen}

Avant de rentrer dans les d\'etails de la d\'emonstration du th\'eor\`eme A, 
nous en donnons l'id\'ee essentielle. Supposons que $I$ soit un intervalle ouvert 
errant pour la dynamique de deux diff\'eomorphismes $g_1$ et $g_2$ du cercle qui 
sont de classe $C^{1+\tau}$ et qui commutent entre eux. Remarquons que l'ensemble 
$$\{ g_1^m g_2^n(I): (m,n) \in \mathbb{N}_0 \times \mathbb{N}_0, m+n \leq k \}$$ 
contient exactemment $(k+1)(k+2)/2$ intervalles. Puisqu'ils sont deux 
\`a deux disjoints, leur \begin{tiny}$^{_\ll}$\end{tiny}longueur 
typique~\begin{tiny}$^{_\gg}$\end{tiny} est de l'ordre de $1/k^2$. 
Donc, pour une \begin{tiny}$^{_\ll}$\end{tiny}suite al\'eatoire 
typique\hspace{0.05cm}\begin{tiny}$^{_\gg}$\end{tiny} $I,h_1(I),h_2(I) \ldots$, 
o\`u $h_{n+1} = g_1 h_n$ ou $h_{n+1} = g_2  h_n$, on s'attend \`a ce que, pour $\tau>1/2$, 
$$\sum_{k \geq 1} |h_k(I)|^{\tau} \leq C \sum_{k \geq 1} \frac{1}{k^{2\tau}} < \infty.$$
Or, la s\'erie \`a gauche est exactement celle dont la convergence permet 
de contr\^oler les distorsions, et donc de trouver des \'el\'ements avec des points 
fixes (hyperboliques), contredisant ainsi la libert\'e de l'action.

Afin de \begin{tiny}$^{_\ll}$\end{tiny}modeler\hspace{0.05cm}\begin{tiny}$^{_\gg}$\end{tiny} 
une preuve dans l'esprit de l'id\'ee ci-dessus, 
nous devons pr\'eciser quelles sont nos \begin{tiny}$^{_\ll}$\end{tiny}suites 
al\'eatoires typiques\hspace{0.05cm}\begin{tiny}$^{_\gg}$\end{tiny}. Pour cela, 
consid\'erons le processus de Markov sur $\mathbb{N}_0 \times \mathbb{N}_0$ 
dont les probabilit\'es de transition sont
\begin{equation}
p \big( (m,n) \rightarrow (m+1,n) \big) = \frac{m+1}{m+n+2} \qquad \mbox{et} 
\qquad p \big( (m,n) \rightarrow (m,n+1) \big) = \frac{n+1}{m+n+2}.
\label{paso}
\end{equation}
Ce processus markovien induit une mesure de probabilit\'e 
$\mathbb{P}$ sur l'espace de chemins correspondant $\Omega$.
On v\'erifie ais\'ement que, en partant de l'origine, la probabilit\'e 
d'arriver en $k$ pas au point $(m,n)$ est \'egale \`a $1/(k+1)$ (resp. 0) 
si \thinspace $m+n=k$ \thinspace (resp. \thinspace $m+n \neq k$). 

\begin{rem} Il est int\'eressant de constater que les probabilit\'es de passage (\ref{paso}) 
ci-dessus co\"{\i}ncident avec celles qui apparaissent dans le mod\`ele d'urne de Polya.
\end{rem}

\begin{rem} On peut canoniquement identifier l'espace $(\Omega,\mathbb{P})$ \`a 
l'intervalle unit\'e (muni de la mesure de Lebesgue). Pour cela, \`a chaque 
\thinspace $x \in [0,1]$ \thinspace on associe le chemin dont la position au $k$-i\`eme 
pas est \'egale \`a \thinspace $\big( [kx],k\!-\![kx] \big)$. \thinspace Avec 
cette identification, la propri\'et\'e d'\'equidistribution pr\'ec\'edente 
devient compl\`etement naturelle.
\end{rem}


\beginpicture

\setcoordinatesystem units <1cm,1cm>

\putrule from 1.5 0 to 1.5 6 
\putrule from 3 0 to 3 6 
\putrule from 4.5 0 to 4.5 6

\putrule from 0 1.5 to 6 1.5    
\putrule from 0 3 to 6 3  
\putrule from 0 4.5 to 6 4.5  

\putrule from 0 0 to 0 6  
\putrule from 0 0 to 6 0

\putrule from 3 -0.6 to 5 -0.6 
\putrule from -0.9 2 to -0.9 4 

\plot 5 -0.6 
4.85 -0.65 /

\plot 5 -0.6 
4.85 -0.55 /

\plot -0.9 4
-0.95 3.85 /

\plot -0.9 4 
-0.85 3.85 /

\put{Figure 2} at 4 -1.3


\plot 1 0 
0.85 0.05 /

\plot 1 0 
0.85 -0.05 /

\plot 1 1.5 
0.85 1.55 /

\plot 1 1.5 
0.85 1.45 /

\plot 1 3 
0.85 2.95 /

\plot 1 3 
0.85 3.05 /

\plot 1 4.5 
0.85 4.55 /

\plot 1 4.5 
0.85 4.45 /



\plot 2.5 0 
2.35 0.05 /

\plot 2.5 0 
2.35 -0.05 /

\plot 2.5 1.5 
2.35 1.55 /

\plot 2.5 1.5 
2.35 1.45 /

\plot 2.5 3 
2.35 2.95 /

\plot 2.5 3 
2.35 3.05 /

\plot 2.5 4.5 
2.35 4.55 /

\plot 2.5 4.5 
2.35 4.45 /



\plot 4 0 
3.85 0.05 /

\plot 4 0 
3.85 -0.05 /

\plot 4 1.5 
3.85 1.55 /

\plot 4 1.5 
3.85 1.45 /

\plot 4 3 
3.85 2.95 /

\plot 4 3 
3.85 3.05 /

\plot 4 4.5 
3.85 4.55 /

\plot 4 4.5 
3.85 4.45 /


\plot 5.5 0 
5.35 0.05 /

\plot 5.5 0 
5.35 -0.05 /

\plot 5.5 1.5 
5.35 1.55 /

\plot 5.5 1.5 
5.35 1.45 /

\plot 5.5 3 
5.35 2.95 /

\plot 5.5 3 
5.35 3.05 /

\plot 5.5 4.5 
5.35 4.55 /

\plot 5.5 4.5 
5.35 4.45 /



\plot 0 1
-0.05 0.85 /

\plot 0 1 
0.05 0.85 /

\plot 1.5 1
1.45 0.85 /

\plot 1.5 1 
1.55 0.85 /

\plot 3 1
2.95 0.85 /

\plot 3 1 
3.05 0.85 /

\plot 4.5 1
4.55 0.85 /

\plot 4.5 1 
4.45 0.85 /



\plot 0 2.5
-0.05 2.35 /

\plot 0 2.5 
0.05 2.35 /

\plot 1.5 2.5
1.45 2.35 /

\plot 1.5 2.5 
1.55 2.35 /

\plot 3 2.5
2.95 2.35 /

\plot 3 2.5 
3.05 2.35 /

\plot 4.5 2.5
4.45 2.35 /

\plot 4.5 2.5 
4.55 2.35 /



\plot 0 4
-0.05 3.85 /

\plot 0 4 
0.05 3.85 /

\plot 1.5 4
1.45 3.85 /

\plot 1.5 4 
1.55 3.85 /

\plot 3 4
2.95 3.85 /

\plot 3 4 
3.05 3.85 /

\plot 4.5 4
4.45 3.85 /

\plot 4.5 4 
4.55 3.85 /



\plot 0 5.5
-0.05 5.35 /

\plot 0 5.5 
0.05 5.35 /

\plot 1.5 5.5
1.45 5.35 /

\plot 1.5 5.5 
1.55 5.35 /

\plot 3 5.5
2.95 5.35 /

\plot 3 5.5
3.05 5.35 /

\plot 4.5 5.5
4.45 5.35 /

\plot 4.5 5.5 
4.55 5.35 /


\setdots

\plot 
0 0 
6 6 /

\small
 
\put{$(0,0)$} at -0.4 -0.4

\put{$p \big( (m,n) \rightarrow (m,n+1) \big) \geq \frac{1}{2}$} at 3 6.5  

\put{$n \geq m \hspace{0.3cm} \Longrightarrow $} at 3 6.9 

\put{$p \big( (m,n) \rightarrow (m+1,n) \big) \geq \frac{1}{2}$} at 7 3.58  

\put{$m \geq n \hspace{0.3cm} \Longrightarrow $} at 7 3.98 

\put{$1/2$} at 0.6 -0.3 
\put{$1/3$} at 0.6 1.2 
\put{$1/4$} at 0.6 2.7
\put{$2/3$} at 2.1 -0.3 
\put{$2/4$} at 2.1 1.2 
\put{$3/4$} at 3.6 -0.3 

\put{$1/2$} at -0.4 0.55
\put{$1/3$} at 1.1 0.55 
\put{$1/4$} at 2.6 0.55
\put{$2/3$} at -0.4 2.05
\put{$2/4$} at 1.1 2.05 
\put{$3/4$} at -0.4 3.55

\put{$g_1$} at 4 -0.8
\put{$g_2$} at -1.1 3 

\put{$\bullet$} at 0 0 

\put{} at -4.1 0 

\endpicture


\vspace{0.5cm}

Pour d\'emontrer le th\'eor\`eme A dans le cas $d=2$, nous proc\'edons par 
contradiction. Soient $g_1$ et $g_2$ deux diff\'eomorphismes du cercle
de classe $C^{1+\tau}$ qui commutent. Le semigroupe $\Gamma^+$ engendr\'e par 
$g_1$ et $g_2$ s'identifie \`a $\mathbb{N}_0 \times \mathbb{N}_0$. Par 
suite, le processus markovien d\'ecrit pr\'ec\'edemment induit une 
\begin{tiny}$^{_\ll}$\end{tiny}promenade 
al\'eatoire\hspace{0.05cm}\begin{tiny}$^{_\gg}$\end{tiny} 
sur $\Gamma^+$. Dans ce qui suit nous identifierons $\Omega$ \`a l'espace 
des chemins correspondants sur $\Gamma^+$. Pour tout $\omega \in \Omega$ 
et tout $n \in \mathbb{N}$, d\'esignons par $h_n(\omega)\in \Gamma^+$ 
le produit des $n$ premi\`eres cordonn\'ees de $\omega$. 
Autrement dit, pour $\omega = (g_{i_1},g_{i_2},\ldots) \in \Omega$ notons 
$h_n(\omega) = g_{i_n} \cdots g_{i_1}$ (o\`u les $g_{i_j}$ ce sont des 
$g_1$ ou des $g_2$), et convenons que $h_0(\omega) = id$.

Si l'action du groupe $\Gamma \!=\! \langle g_1, g_2 \rangle \!\sim \!\mathbb{Z}^2$ est libre, 
alors les nombres de rotation $\rho(g_1)$ et $\rho(g_2)$ de $g_1$ et $g_2$ respectivement sont 
{\em ind\'ependants} sur les rationnels, dans le sens que pour tout $(r_0,r_1,r_2) \in \mathbb{Q}^3$ 
distinct de $(0,0,0)$ on a \esp $r_1 \rho(g_1) + r_2 \rho(g_2) \neq r_0$. \esp En effet, dans 
le cas contraire on pourrait trouver des \'el\'ements non triviaux (et donc d'ordre infini) 
avec un nombre de rotation rationnel ; de tels \'el\'ements poss\'edant des points 
p\'eriodiques, ceci contredirait la libert\'e de l'action. 

Supposons maintenant que l'action de $\Gamma$ soit (libre et) non minimale. Dans 
ce cas, il existe un ensemble de Cantor invariant et minimal pour l'action. De 
plus, toute composante connexe $I$ du compl\'ementaire de cet ensemble 
est errante pour la dynamique, dans le sens que ses images par des 
\'el\'ements distincts du groupe sont disjointes.

\vspace{0.1cm}

\begin{lem} {\em Si $\tau > 1/2$ alors la valeur de la s\'erie \thinspace 
$\sum_{n \geq 0} |h_n(\omega)(I)|^{\tau}$ \thinspace est finie pour 
$\mathbb{P}$-presque tout $\omega \in \Omega$.}
\label{2d}
\end{lem}

\noindent{\bf D\'emonstration.} Pour tout chemin $\omega \in \Omega$ 
d\'efinissons 
$$\ell_{\tau} (\omega) = \sum_{k \geq 0} |h_k(\omega)(I)|^{\tau}.$$ 
Nous allons v\'erifier que, si $\tau > 1/2$, alors l'esp\'erance (par 
rapport \`a $\mathbb{P}$) de la fonction $\ell_{\tau}$ est finie, ce qui 
implique \'evidemment l'affirmation du lemme. Remarquons d'abord que 
$$\mathbb{E}(\ell_{\tau}) 
= \mathbb{E} \Big( \sum_{k \geq 0} |h_k(\omega)(I)|^{\tau} \Big)
= \sum_{k \geq 0} \mathbb{E} \big( |h_k(\omega)(I)|^{\tau} \big)
= \sum_{k \geq 0} \sum_{m+n=k} \frac{|g_1^m g_2^n(I)|^{\tau}}{k+1}.$$
Or, l'in\'egalit\'e de H\"older montre ais\'ement que 
$$\sum_{m+n=k} \frac{|g_1^m g_2^n(I)|^{\tau}}{k+1} \leq 
\left( \sum_{m+n=k} |g_1^m g_2^n(I)| \right)^{\tau} 
\left( (k+1) \cdot \frac{1}{(k+1)^{1/(1-\tau)}}\right)^{1-\tau},$$
et donc, 
$$\mathbb{E}(\ell_{\tau}) \leq \sum_{k \geq 0} 
\frac{\left(\sum_{m+n=k}|g_1^m g_2^n(I)|\right)^{\tau}}{(k+1)^{\tau}}.$$
En appliquant une nouvelle fois l'in\'egalit\'e de H\"older on obtient
$$\mathbb{E}(\ell_{\tau}) \leq 
\left[ \sum_{(m,n) \in \mathbb{N}_0 \times \mathbb{N}_0} 
|g_1^m g_2^n(I)| \right]^{\tau} \left[ \sum_{k \geq 1} \left( \frac{1}{k^{\tau}} 
\right)^{\frac{1}{1-\tau}} \right]^{1-\tau}.$$
\'Etant donn\'e que $\tau > 1/2$, la s\'erie
$$\sum_{k \geq 1} \left( \frac{1}{k^{\tau}} \right)^{\frac{1}{1-\tau}} = 
\sum_{k \geq 1} \frac{1}{k^{\tau/(1-\tau)}}$$
converge, et puisque les intervalles $g_1^mg_2^n(I)$ sont deux \`a deux 
disjoints, ceci montre la finitude de $\mathbb{E}(\ell_{\tau})$. \esp $\square$

\vspace{0.1cm}

\begin{rem} Dans la d\'emonstration pr\'ec\'edente, la seule propri\'et\'e 
du processus de diffusion sur $\mathbb{N}_0 \times \mathbb{N}_0$ que l'on 
a utilis\'e est le fait que les probabilit\'es d'arriv\'ee en $k$ pas sont 
\'equidistribu\'ees sur l'ensemble des points \`a distance 
(simpliciale) $k$ de l'origine.
\end{rem}

D'apr\`es le lemme pr\'ecedent, si $M$ est suffisamment grand alors l'ensemble 
$\Omega(M) \!=\! \big\{ \omega \in \Omega: \ell_{\tau}(\omega) \leq M \big\}$ 
poss\`ede une probabilit\'e strictement positive (en fait, $\mathbb{P}[\Omega(M)]$ 
converge vers $1$ lorsque $M$ tend vers l'in- fini). Fixons un tel $M$, 
et soit \thinspace 
$L \!=\! L(\tau,M,|I|;\{g_1,g_2\})$ \thinspace la constante du lemme \ref{?}. 
Consid\'erons finalement l'intervalle ouvert \esp $K'$ \esp de taille 
$|K'| = L$ et adjacent \`a droite \`a $I$. Nous affirmons que 
\begin{equation}
\mathbb{P} \big[ \omega \in \Omega: h_n(\omega)(I) \not\subset K' 
\mbox{ pour tout } n \in \mathbb{N} \big] = 0.
\label{total}
\end{equation}

Pour d\'emontrer (\ref{total}) remarquons d'abord que l'action du groupe 
engendr\'e par les diff\'eomorphismes $g_1$ et $g_2$ est semiconjugu\'ee 
\`a une action par des rotations. Par suite, si l'on 
\begin{tiny}$^{_\ll}$\end{tiny}\'ecrase\hspace{0.05cm}\begin{tiny}$^{_\gg}$\end{tiny} 
chaque composante connexe du compl\'ementaire du Cantor invariant minimal 
$\Lambda$, alors 
on obtient un cercle topologique $\clo_{\Lambda}$ sur lequel $g_1$ et $g_2$ 
induisent des hom\'eomorphismes dont toutes les orbites sont denses. De plus, 
l'intervalle $K'$ devient un intervalle $U$ d'int\'erieur non vide dans 
$\clo_{\Lambda}$. Les nombres de rotation de $g_1$ et $g_2$ \'etant irrationnels, 
il existe $N \in \mathbb{N}$ tel que, apr\`es \'ecrasement, $g_1^{-1}(U), \ldots, 
g_1^{-N}(U)$ recouvrent le cercle topologique $\clo_{\Lambda}$, 
et de m\^eme pour $g_2^{-1}(U),\ldots,g_2^{-N}(U)$. Sur le cercle original 
$\clo$ cela se traduit par le fait que, pour toute composante connexe $I_0$ de 
$\clo \setminus \Lambda$, il existe $n_1$ et $n_2$ dans $\{1,\ldots,N\}$ 
tels que \thinspace $g_1^{n_1}(I_0) \subset K'$ \thinspace et 
\thinspace $g_2^{n_2}(I_0) \subset K'$. \thinspace 

Soulignons maintenant la propri\'et\'e \'el\'ementaire suivante 
et qui d\'ecoule directement de (\ref{paso})~: 
les probabilit\'es de passage \`a droite (resp. vers le haut) du processus 
markovien consid\'er\'e sont $\geq 1/2$ au dessous (resp. au dessus) de la 
diagonale (voir la figure 2). D'apr\`es la d\'efinition de $N$, cette propri\'et\'e 
donne, pour tout entier $k \geq 0$, 
$$\mathbb{P} \big[ g_1^i h_k(\omega)(I) \not\subset K' \thinspace \mbox { et }
\thinspace g_2^i h_k(\omega)(I) \not\subset K' \mbox{ pour tout } i \leq N  
\thinspace \big| \thinspace h_j(\omega)(I) \not\subset K' \mbox{ pour tout } 
j \leq k \big] \leq 1 - \frac{1}{2^N}.$$
Cette derni\`ere in\'egalit\'e implique \'evidemment que 
\begin{equation}
\mathbb{P} \big[ h_{k+i} (\omega)(I) \not\subset K' 
\mbox{ pour tout } i \leq N  \thinspace 
\big| \thinspace h_j(\omega)(I) \not\subset K' \mbox{ pour tout } 
j \leq k \big] \leq 1 - \frac{1}{2^N}.
\label{total2}
\end{equation}
Par cons\'equent, pour tout $r \in \mathbb{N}$,
$$\mathbb{P} \big[ h_n(\omega)(I) \not\subset K' \mbox{ pour tout } n \in \mathbb{N} \big] 
\leq \mathbb{P} \big[ h_n(\omega)(I) \not\subset K' \mbox{ pour tout } 
n \in \{ 1, \ldots, rN \} \big] \leq \Big( 1 - \frac{1}{2^N} \Big)^r,$$
d'o\`u l'on obtient (\ref{total}) en faisant tendre $r$ vers l'infini. 

Pour finir la preuve du th\'eor\`eme A (toujours dans le cas $d=2$), remarquons que  si 
$\omega \in \Omega(M)$ et $n \in \mathbb{N}$ sont tels que $h_n(\omega)(I) \subset K'$, 
alors le lemme \ref{?} permet de trouver un point fixe hyperbolique pour 
$h_n(\omega) \in \Gamma^+$, contredisant ainsi l'hypoth\`ese de libert\'e de l'action.

\vspace{0.35cm}

\noindent{\bf Modifications pour le cas ${\bf d > 2}$.} La d\'emonstration du 
th\'eor\`eme A pour $d > 2$ est tout \`a fait analogue \`a celle donn\'ee 
pour le cas $d = 2$. Elle se fait aussi par contradiction : on suppose 
l'existence d'un intervalle errant et on consid\`ere le processus 
markovien sur $\mathbb{N}_0^d$ \`a probabilit\'es de transition 
$$p \big( (n_1,\ldots,n_i,\ldots,n_d) \longrightarrow 
(n_1,\ldots,1+n_i,\ldots,n_d) \big) = \frac{1+n_i}{n_1 + \cdots + n_d + d}.$$
Les probabilit\'es d'arriv\'ee en $k$ pas pour ce processus sont aussi 
\'equidistribu\'ees sur l'ensemble des points \`a distance (simpliciale) $k$ 
de l'origine. Cela permet \`a nouveau de contr\^oler les distortions pour 
presque toute suite al\'eatoire, c'est-\`a-dire de d\'emontrer un analogue du 
lemme~\ref{2d} lorsque $\tau \!>\! 1/d$. On remarque ensuite que chaque point 
$(n_1,\ldots,n_d)$ de $\mathbb{N}_0^d$ est le point de d\'epart d'au moins 
une ligne droite telle les probabilit\'es de passage entre deux sommets 
cons\'ecutifs est sup\'erieure ou \'egale \`a $1/d$ (il suffit de suivre 
la direction de la coordonn\'ee $i$-\`eme pour laquelle $n_i$ prend la 
valeur la plus grande). Cela permet d'obtenir une in\'egalit\'e analogue 
\`a (\ref{total2}) (dont le membre \`a droite sera \'egal \`a 
\thinspace $1\!-\!1/d^N$ \thinspace  
pour un certain entier $N$ assez grand). Une telle in\'egalit\'e 
entra\^{\i}ne la propri\'et\'e (\ref{total}), qui gr\^ace au 
contr\^ole de distorsion permet d'utiliser le lemme \ref{?}. 
On trouve ainsi des \'el\'ements avec des points fixes hyperboliques, 
contredisant une nouvelle fois la libert\'e de l'action.

\vspace{0.1cm}

\begin{rem} Dans la d\'emonstration pr\'ec\'edente nous n'avons eu besoin 
de la finitude de la fonction $\ell_{\tau}$ que pour un ensemble de suites 
de mesure strictement positive. 
N\'eanmoins, la m\'ethode du lemme \ref{2d} m\`ene \`a une conclusion beaucoup plus forte~: 
l'esp\'erance de la fonction $\ell_{\tau}$ est finie d\`es que $\tau > 1/d$. 
On pourrait donc imaginer qu'en affaiblisant cette derni\`ere affirmation on puisse attaquer 
le cas critique $\tau \!=\! 1/d$ par des m\'ethodes analogues. On constatera cependant que pour 
la preuve du lemme \ref{2d} nous nous sommes appuy\'es uniquemment sur le fait que $I$ \'etait 
un intervalle errant. Or, c'est aussi le cas de l'intervalle $I_{0,\ldots,0}$ de l'exemple du 
\S \ref{ex1}, alors que pour cet exemple la valeur de la fonction $\ell_{1\!/\!d}$ est infinie 
pour toute suite $\omega$ (lorsque $\varepsilon \leq d-1$). Ceci indique que pour d\'emontrer 
le th\'eor\`eme dans le cas critique, il est n\'ecessaire d'introduire des m\'ethodes 
plus fines et qui prennent en compte la nature dynamique des intervalles 
$g_1^m g_2^n(I)$ (et non seulement le fait qu'il soient deux \`a deux disjoints). En effet, 
la preuve pr\'ec\'edente ne permet pas de conclure que (pour $\varepsilon \leq d-1$) 
l'exemple du \S \ref{ex1} n'est pas de classe $C^{1+1\!/\!d}$, alors que ceci est ais\'ement 
v\'erifiable d'apr\`es les d\'efinitions (et r\'esulte \'egalement de la proposition 
\ref{tangente}). 
\end{rem}

Nous voudrions conclure ce paragraphe en donnant le sch\'ema d'une autre preuve du th\'eor\`eme A 
dont l'id\'ee sera essentielle \`a la fin du \S \ref{kop}. Pour simplifier, nous ne consid\'erons 
que le cas $d = 2$, et nous gardons les notations introduites tout au long de ce paragraphe. En 
raffinant l\'eg\`erement les arguments donn\'es pr\'ec\`e- demment, on d\'emontre que 
pour presque toute suite $\omega \in \Omega$ et pour tout point $p$ appartenant au cercle 
topolo- gique $\clo_{\Lambda}$, l'ensemble des points de la forme $h_n(\omega)(p)$ est 
dense dans $\clo_{\Lambda}$. Fixons $M > 0$ suffisamment grand de fa\c{c}on \`a 
ce que $\mathbb{P}[\Omega(M)] > 0$, et comme dans la preuve du lemme \ref{?} notons 
$J \! = \! I' \cup I \cup I''$ le $2L$-voisinage de $I$. Choisissons un \'el\'ement 
$h \in \Gamma$ tel que $h(I)$ et $h^{-1}(I)$ soient des intervalles de taille strictement 
inf\'erieure \`a \thinspace $|I| \exp(-2^{\tau} C M)$ \thinspace contenus dans $I'$ et $I''$ 
respectivement. Fixons un point arbitraire $a \in \Lambda$, et supposons que $h'(a) \geq 1$ 
(le cas o\`u $h'(a) \leq 1$ est r\'egl\'e en appliquant les arguments qui suivent \`a $h^{-1}$ 
au lieu de $h$, tout en remarquant que dans ce cas $(h^{-1})'(h(a)) = 1 / h'(a) \geq 1$ et 
$h(a) \in \Lambda$). La premi\`ere partie de la d\'emonstration du lemme \ref{?} 
montre que pour tout $x$ et $y$ appartenant \`a $\bar{I}' \cup \bar{I}$, tout 
$\omega \in \Omega(M)$ et tout $n \in \mathbb{N}$, 
$$\frac{h_n(\omega)'(x)}{h_n(\omega)'(y)} \leq \exp(2^{\tau} C M).$$
\`A partir de l'\'egalit\'e \thinspace $h = h_n(\omega)^{-1} \circ h \circ h_n(\omega)$ 
\thinspace on en d\'eduit que, pour tout $y \in I$,
$$h'(y) = \frac{h_n(\omega)'(y)}{h_n(\omega)'(h(y))} \thinspace 
h' \big( h_n(\omega)(y) \big) \geq \exp (-2^{\tau} C M ) \thinspace 
\limsup_{n \in \mathbb{N}} h' \big( h_n(\omega)(y) \big).$$
Or, puisque la suite $\big( h_n(\omega)(p) \big)$ est dense 
sur $\clo_{\Lambda}$ pour tout $p \!\in\! \clo_{\Lambda}$, 
il existe une suite croissante et infinie d'en- tiers $n_k$ tels que l'intervalle 
$h_{n_k}(\omega)(I)$ tend vers le point $a$. Par suite, \thinspace 
$\limsup_{n \in \mathbb{N}} h' \big( h_n(\omega)(y) \big) \geq h'(a) \geq 1$, 
\thinspace et donc \thinspace $h'(y) \geq \exp(-2^{\tau}CM)$ \thinspace pour tout $y \in I$. 
Ceci implique que la taille de l'intervalle $h(I)$ est au moins \'egale \`a \thinspace 
$|I| \exp(-2^{\tau} C M)$, \thinspace ce qui contredit notre choix de $h$.


\section{Le lemme de Kopell g\'en\'eralis\'e}
\label{kop}

De mani\`ere analogue \`a ce que nous avons fait pour la g\'en\'eralisation 
du th\'eor\`eme de Denjoy, pour la preuve du th\'eor\`eme B nous 
ne consid\'ererons que le cas $d = 2$, et nous laisserons 
au lecteur le soin d'adapter nos arguments au cas $d \geq 3$. 
Tout en gardant les notations de l'\'enonc\'e du th\'eor\`eme (avec $d=2$), identifions 
le semigroupe $\Gamma^+$ engendr\'e par les \'el\'ements $f_1$ et $f_2$ de 
$\mathrm{Diff}_+^{1+\tau} \left( [0,1] \right)$ avec $\mathbb{N}_0 \times \mathbb{N}_0$, 
et consid\'erons le processus markovien du \S \ref{denjoy-gen}. Si l'on fixe 
l'intervalle $I=I_{0,0}$ et pour chaque $\omega \in \Omega$ on d\'efinit
$$\ell_{\tau}(\omega) = \sum_{i \geq 0} |h_i(\omega)(I)|^{\tau},$$
alors l'argument de la preuve du lemme \ref{2d} montre que, lorsque $\tau > 1/2$, 
la fonction $\ell_{\tau} : \Omega  \rightarrow \mathbb{R}$ est presque s\^urement 
finie (et qu'en fait, son esp\'erance est finie).

Soit $C$ une constante de $\tau$-h\"olderianit\'e pour $\log(f_1')$ et $\log(f_2')$. Pour 
chaque $\omega = (f_{j_1},f_{j_2},\ldots) \in \Omega$, tout $n,k$ dans $\mathbb{N}$ et tout 
$x \in I$, l'\'egalit\'e \thinspace 
$f_3^n = h_k(\omega)^{-1} \circ f_3^n \circ h_k(\omega)$ \thinspace donne
$$\log \big( (f_3^n)'(x)\big) = \log \big( (f_3^n)' (h_k(\omega)(x)) \big) + 
\sum_{i=1}^k \big[ \log \big( f_{j_i}'(h_{i-1}(\omega)(x)) \big) - 
 \log \big( f_{j_i}'(f_3^n \circ h_{i-1}(\omega)(x)) \big) \big],$$
et donc
\vspace{-0.2cm}
\begin{eqnarray*}
\big| \log \big( (f_3^n)'(x)\big) \big| 
&\leq& \big| \log \big( (f_3^n)' (h_k(\omega)(x)) \big) \big| + 
C \sum_{i=1}^k \big| h_{i-1}(\omega)(x) - f_3^n \circ h_{i-1}(\omega)(x) \big|^{\tau} \\
&\leq& \big| \log \big( (f_3^n)' (h_k(\omega)(x)) \big) \big| + C \sum_{i=1}^k 
|h_{i-1}(\omega)(I)|^{\tau}.
\end{eqnarray*}
En prenant $\omega \in \Omega$ tel que $\ell_{\tau}(\omega) = M$ soit fini, l'in\'egalit\'e 
pr\'ec\'edente implique que 
$$\big| \log \big( (f_3^n)'(x)\big) \big| \leq 
\big| \log \big( (f_3^n)' (h_k(\omega)(x)) \big) \big| + CM.$$
Or, le point $h_k(\omega)(x)$ converge n\'ecessairement vers un point fixe (parabolique) de $f_3$. 
En faisant tendre $k$ vers l'infini on en d\'eduit que \thinspace 
$\big| \log \big( (f_3^n)'(x)\big) \big| \leq CM.$ \thinspace 
Par suite, $(f_3^n)'(x) \leq \exp(CM)$ pour tout $x \in I$ et tout $n \in \mathbb{N}$, ce qui 
entra\^{\i}ne \'evidemment que la restriction de $f_3$ \`a l'intervalle $I$ est l'identit\'e. 
Par commutativit\'e, ceci est aussi vrai sur tous les intervalles $I_{n_1,n_2}$, ce qui conclut
la d\'emonstration.

\begin{rem} Une lecture attentive de la preuve pr\'ec\'edente montre que le th\'eor\`eme B 
est encore valable pour des diff\'eomorphismes de classe $C^{1+\tau}$ de l'intervalle $[0,1[$, 
{\em i.e.} pour des diff\'eomorphismes dont on ne dispose d'une borne uniforme pour la 
constante de $\tau$-h\"olderianit\'e que sur chaque intervalle compact contenu dans $[0,1[$. 
\end{rem}

Le lecteur pourrait \^etre incommod\'e \`a cause de 
l'hypoth\`ese combinatoire (\ref{comb}). N\'eanmoins, nous avons 
d\'ej\`a expliqu\'e que les contre-exemples au th\'eor\`eme B qui v\'erifient cette condition 
correspondent aux actions de $\mathbb{Z}^{d+1}$ sur $[0,1]$ \begin{tiny}$^{_\ll}$\end{tiny}les 
plus int\'eressantes\hspace{0.05cm}\begin{tiny}$^{_\gg}$\end{tiny} du point de vue cohomologique. 
Ils portent aussi un inter\^et dynamique, car il n'est pas difficile de construire des actions 
de $\mathbb{Z}^{d+1}$ par diff\'eomorphismes de classe $C^{2-\varepsilon}$ de l'intervalle sans 
point global \`a l'int\'erieur. En effet, consid\'erons par exemple des diff\'eomorphismes 
$f_2,\ldots,f_{d+1}$ d'un intervalle $[a,b]\! \subset ]0,1[$ dont les supports soient disjoints, 
et soit $f_1$ un diff\'eomorphisme de l'intervalle 
$[0,1]$ dans lui m\^eme sans point fixe \`a l'int\'erieur et qui envoie 
$]a,b[$ sur un intervalle disjoint. En \'etendant $f_2,\ldots,f_{d+1}$ \`a tout l'intervalle 
$[0,1]$ de fa\c{c}on \`a ce qu'ils commutent avec $f_1$, on obtient une action (fid\`ele) 
de $\mathbb{Z}^{d+1}$ par hom\'eomorphismes de l'intervalle et sans point fixe global \`a 
l'int\'erieur. Bien s\^ur, les m\'ethodes du \S \ref{ex2} (resp. de \cite{TsP}) permettent 
de rendre cette action de classe $C^{3/2 - \varepsilon}$ (resp. $C^{2-\varepsilon}$) pour 
tout $\varepsilon > 0$.

L'exemple na\"{\i}f ci-dessus rend naturel le probl\`eme d'obtenir un r\'esultat qui 
g\'en\'eralise simultan\'ement les th\'eor\`emes A et B, et qui permette de d\'ecrire (toutes) 
les actions de $\mathbb{Z}^{d+1}$ par diff\'eomorphismes du cercle 
de r\'egularit\'e interm\'ediaire~; pour cela il faudrait \'etudier les actions 
non libres mais sans point fixe global. 
Remarquons \`a ce propos l'existence d'un exemple int\'eressant : il s'agit d'une action 
de $\mathbb{Z}^2$ par diff\'eomorphismes de classe $C^{2-\varepsilon}$ du cercle. Dans 
cet exemple, l'un des g\'en\'erateurs correspond \`a un contre-exemple de Denjoy, 
alors que l'autre g\'en\'erateur laisse invariant le Cantor minimal du premier et 
agit non trivialement sur les composantes connexes de son compl\'ementaire. Pour 
rendre cette action de classe $C^{2-\varepsilon}$, on doit utiliser les m\'ethodes 
de \cite{TsP}.

Un autre aspect int\'eressant et important \`a remarquer est l'existence d'actions 
sur l'intervalle qui \`a l'int\'erieur sont libres mais non conjugu\'ees \`a des actions 
par des translations. Ce ph\'enom\`ene ne peut pas se produire en classe $C^{1+vb}$ 
(voir le lemme 3.2 de \cite{Na-subexp}). Cependant, en classe $C^{3/2 - \varepsilon}$ 
on peut construire des actions de $\mathbb{Z}^2$ sur $[0,1]$ qui sont libres sur $]0,1[$ 
mais qui admettent des intervalles errants, tout en utilisant les techniques du \S 
\ref{exemples}. Une nouvelle fois, la classe de diff\'erentiabilit\'e \thinspace 
$3/2\! -\! \varepsilon$ \thinspace est optimale pour de tels exemples. Le r\'esultat 
suivant est \`a rapprocher avec celui qui stipule l'existence de champs de vecteurs 
associ\'es aux diff\'eomorphismes de classe $C^2$ de l'intervalle sans point fixe 
\`a l'int\'erieur \cite{serg}.

\vspace{0,1cm} 

\begin{prop} {\em Soit $\Gamma$ un sous-groupe de $\mathrm{Diff}_+^{1+\tau}([0,1])$ isomorphe 
\`a $\mathbb{Z}^d$, avec $\tau > 1/d$ et $d \geq 2$. Si la restriction \`a $]0,1[$ de 
l'action de $\Gamma$ est libre, alors elle est minimale et topologiquement 
conjugu\'ee \`a l'action d'un groupe de translations.}
\label{no-semi}
\end{prop}

\vspace{0.08cm}

Avant de passer \`a la d\'emonstration, rappelons qu'un th\'eor\`eme classique d\^u \`a H\"older 
stipule que tout groupe d'hom\'eomorphismes de l'intervalle $]0,1[$ qui agit librement est 
topologiquement semi-conjugu\'e \`a un groupe de translations (voir \cite{ghys} pour une 
preuve compl\`ete de ce r\'esultat). Lorsque l'action est minimale, cette semi-conjugaison 
est forc\'ement une conjugaison ; dans le cas contraire, on voit appara\^{\i}tre des 
intervalles errants pour la dynamique.

Supposons maintenant que $f_1,\ldots,f_d$ soient les g\'en\'erateurs d'un groupe 
$\Gamma \!\sim\! \mathbb{Z}^d$ qui agit librement sur $]0,1[$ mais qui n'est pas 
conjugu\'e \`a un groupe de translations. Si l'on identifie les points des orbites 
par $f_1$, alors $f_2,\ldots,f_d$ deviennent les g\'en\'erateurs d'un groupe agissant 
librement sur le cercle et qui n'est pas conjugu\'e \`a un groupe de rotations ; bien 
s\^ur, le th\'eor\`eme A implique que les $f_i$ ne peuvent pas \^etre tous de classe 
$C^{1\!/\!(d-1) + \varepsilon}$. Notons que cet argument n'utilise que la diff\'erentiabilit\'e 
des applications \`a l'int\'erieur ; dans ce contexte la r\'egularit\'e de l'obstruction 
pr\'ec\'edente est en fait optimale. N\'eanmoins, ce n'est pas le cas dans le cadre des 
groupes de diff\'eomorphismes de l'intervalle ferm\'e $[0,1]$ (et m\^eme de $[0,1[$). 
Pour chercher \`a d\'emontrer cela nos devrons donc tenir compte de la r\'egularit\'e 
des applications aux extr\'emit\'es.
  
\vspace{0.3cm}

\noindent{\bf D\'emonstration de la proposition \ref{no-semi}.} Encore une fois, nous ne 
donnons la preuve que pour le cas $d=2$. Soient $f_1$ et $f_2$ les g\'en\'erateurs d'un 
groupe $\Gamma \!\sim\! \mathbb{Z}^2$ de diff\'eomorphismes de classe $C^{1+\tau}$ de $[0,1]$ qui 
agit librement \`a l'int\'erieur. Quitte \`a les \'echanger par leurs inverses, nous pouvons 
supposer qu'ils contractent topologiquement vers l'origine. Supposons que l'action de $\Gamma$ 
sur $]0,1[$ ne soit pas conjugu\'ee \`a une action par des translations. Dans ce cas, un 
argument simple de contr\^ole de distorsion hyperbolique montre que tous les \'el\'ements 
de $\Gamma$ doivent \^etre tangents \`a l'identit\'e \`a l'origine. En effet, supposons par 
contradiction qu'il existe des intervalles errants ainsi qu'un \'el\'ement $f \in \Gamma$ 
tel que $f'(0) < 1$. Fixons $\lambda < 1$ et $c \!\in ]0,1[$ tels que $f'(x) \leq \lambda$ 
pour tout $x \!\in \![0,c[$, et fixons un intervalle errant ouvert et maximal $I=]a,b[$ 
contenu dans $]0,c[$. Si l'on d\'esigne par $K$ l'intervalle $[f(b),b]$ alors 
$$\sum_{n \geq 0} |f^n(K)|^{\tau} \leq |K|^{\tau} \sum_{n \geq 0} \lambda^{n \tau} 
= \frac{|K|^{\tau}}{1-\lambda^{\tau}} = \bar{M}.$$
Par cons\'equent, \thinspace $(g^n)'(x) / (g^n)'(y) \leq \exp(C \bar{M})$ 
\thinspace pour tout $x,y$ dans $K$, o\`u $C>0$ est une constante de 
$\tau$-h\"olderianit\'e pour $\log(f_1')$ et $\log(f_2')$. Or, cette derni\`ere 
estim\'ee permet d'appliquer les arguments de la preuve du lemme 3.2 
de \cite{Na-subexp}, en parvenant ainsi \`a une 
contradiction\footnote{Remarquons que dans cette partie de la preuve nous n'avons 
utilis\'e que le fait que $\tau > 0$. Une mani\`ere plus conceptuelle d'expliquer 
le ph\'enom\`ene sous-jacent consiste \`a rappeler que le th\'eor\`eme de lin\'earisation 
de Sternberg \cite{St} est encore valable pour des germes hyperboliques et de classe 
$C^{1+\tau}$ de l'intervalle : la preuve donn\'ee dans \cite{yoccoz} en classe $C^2$ 
s'\'etend ais\'ement dans ce contexte (voir aussi \cite{cha}). Or, le centralisateur 
d'un germe lin\'eaire non trivial est le groupe \`a un param\`etre des germes lin\'eaires ; 
en particulier, l'existence d'intervalles errants pour la dynamique de sous-groupes denses 
de ce centralisateur est interdite.}.

Identifions maintenant le semigroupe $\Gamma^+$ engendr\'e par $f_1$ et $f_2$ \`a 
$\mathbb{N}_0 \times \mathbb{N}_0$, et consid\'erons le processus 
markovien du \S \ref{denjoy-gen}. Lorsque $\tau > 1/2$, la preuve du lemme 
\ref{2d} montre la finitude de l'esp\'erance de la fonction
$$\omega \mapsto \ell_{\tau}(\omega) = \sum_{k \geq 0} |h_k(\omega)(I)|^{\tau}.$$
Prenons $M>0$ suffisamment grand de fa\c{c}on \`a ce que l'ensemble \thinspace 
$\Omega(M) = \{\omega \in \Omega: \ell_{\tau}(\omega) \leq M \}$ \thinspace ait une 
probabilit\'e strictement positive, et notons \thinspace $\bar{L} = |I| / \exp(2^{\tau}CM)$. 
\thinspace D'apr\`es la premi\`ere partie de la preuve du lemme \ref{?}, si $I''$ 
d\'esigne l'intervalle adjacent \`a gauche \`a $I$ et de longueur $\bar{L}$, alors pour tout 
$x$ et $y$ appartenant \`a $J = \bar{I}'' \cup \bar{I}$, tout $\omega \in \Omega(M)$ et tout 
$n \in \mathbb{N}$, 
\begin{equation}
\frac{h_n(\omega)'(x)}{h_n(\omega)'(y)} \leq \exp(2^{\tau} C M).
\label{controlado}
\end{equation}
L'intervalle $I$ n'\'etant pas contenu dans aucun autre intervalle errant et 
ouvert, il existe $h \in \Gamma$ tel que $h(I) \subset I''$ et \thinspace 
$|h(I)| < |I| \thinspace \exp(- 2^{\tau} C M)$. \thinspace 
Fixons un point arbitraire $y$ dans $I$. Puisque $h'(0)=1$ et que $h_n(\omega)(y)$ 
tend vers l'origine pour tout $\omega \in \Omega(M)$, \`a partir de l'\'egalit\'e
$$h'(y) = \frac{h_n(\omega)'(y)}{h_n(\omega)'(h(y))} \thinspace 
h' \big( h_n(\omega)(y) \big)$$ 
et de (\ref{controlado}) on conclut que \thinspace 
$h'(y) \geq \exp(-2^{\tau} C M).$ \thinspace 
Par suite, 
$$|h(I)| \geq |I| \thinspace \exp(-2^{\tau} C M),$$ 
ce qui est en contradiction avec le choix de $h$. \esp $\square$


\section{Le th\'eor\`eme de Sacksteder g\'en\'eralis\'e}
\label{sac}

\subsection{Pseudo-groupes de diff\'eomorphismes et feuilletages : le cas $C^{1+\tau}$}
\label{sac-ps}

Dans ce paragraphe nous donnons la preuve du th\'eor\`eme C en classe 
$C^{1+\tau}$. La cas $C^1$ n\'ecessite l'introduction d'autres techniques, et il 
sera trait\'e au paragraphe suivant.

Comme nous l'avons d\'ej\`a signal\'e, la preuve du th\'eor\`eme C se r\'eduit \`a d\'emontrer 
que, en pr\'esence d'une feuille ressort, il existe des feuilles ressort hyperboliques 
(voir l'appendice \ref{ap0}).  En consid\'erant une transversale au feuilletage qui 
\begin{tiny}$^{_\ll}$\end{tiny}capture\hspace{0.05cm}\begin{tiny}$^{_\gg}$\end{tiny} 
cette feuille ressort topologique, nous sommes amen\'es \`a consid\'erer la dynamique 
sur un intervalle $[a,a']$ de deux diff\'eomorphismes locaux $f$ et $h$ de classe 
$C^{1+\tau}$ qui v\'erifient~:

\vspace{0.1cm}

\noindent{-- \thinspace $f$ est defini sur tout l'intervalle 
$[a,a'[$ et il contracte topologiquement vers le point fixe $a$ ;}

\vspace{0.1cm}

\noindent{-- \thinspace $h$ est defini sur un 
voisinage de ce point fixe et $h(a) \!\in ]a,a'[$.}

\vspace{0.1cm}

Posons $c = h(a)$ et fixons $d' \!\in ]c,a'[$. Quitte \`a remplacer 
$f$ par $f^n$ pour $n \in \mathbb{N}$ assez grand, nous pouvons supposer 
que $f(d') < c$, que $f(d')$ appartient au domaine de d\'efinition de $h$, et 
que $hf (d')\!\in ]c,d'[$. Cette derni\`ere condition implique en particulier 
que $hf$ poss\`ede des points fixes dans $]c,d'[$. Soit $d$ le premier point 
fixe de $hf$ \`a droite de $c$, et soit $b = f(d)$. L'intervalle ouvert 
$I\!= ]b,c[$ correspond au premier {\em gap} d'un ensemble de Cantor 
$\Lambda$ qui est invariant par $f$ et $g = hf$ (voir la figure 3).

\vspace{0.3cm}


\beginpicture

\setcoordinatesystem units <1cm,1cm>

\putrule from -7 0 to -2 0 

\putrule from -7 -0.2 to -7 0.2  

\putrule from -2 -0.2 to -2 0.2 

\putrule from -5 -0.13 to -5 0.13  

\putrule from -4 -0.13 to -4 0.13 

\putrule from -5 -0.015 to -4 -0.015 

\putrule from -5 0.015 to -4 0.015

\putrule from -6.2 -0.07 to -6.2 0.07

\putrule from -5.8 -0.07 to -5.8 0.07

\putrule from -6.2 -0.015 to -5.8 -0.015 

\putrule from -6.2 0.015 to -5.8 0.015

\putrule from -3.2 -0.07 to -3.2 0.07

\putrule from -2.8 -0.07 to -2.8 0.07

\putrule from -3.2 -0.015 to -2.8 -0.015 

\putrule from -3.2 0.015 to -2.8 0.015

\putrule from -5.5 1.1 to -3.5 1.1

\plot -5.5 1.1 
-5.35 1.15 /

\plot -5.5 1.1 
-5.35 1.05 /

\putrule from -6 0.7 to -3 0.7 

\putrule from -6 0.3 to -6 0.7

\putrule from -3 0.3 to -3 0.7 

\plot -3 0.3 
-3.05 0.45 /

\plot -3 0.3 
-2.95 0.45 /


\putrule from 2 -2 to 7 -2
\putrule from 2 3 to 7 3
\putrule from 2 -2 to 2 3
\putrule from 7 -2 to 7 3

\putrule from 2.02 0 to 2.02 1 
\putrule from 1.98 0 to 1.98 1 

\putrule from 1.85 0 to 2.15 0
\putrule from 1.85 1 to 2.15 1

\putrule from 2.02 -1.2 to 2.02 -0.8
\putrule from 1.98 -1.2 to 1.98 -0.8

\putrule from 2.02 1.8 to 2.02 2.2
\putrule from 1.98 1.8 to 1.98 2.2

\putrule from 1.9 -1.2 to 2.1 -1.2   
\putrule from 1.9 -0.8 to 2.1 -0.8  
\putrule from 1.9 1.8 to 2.1 1.8 
\putrule from 1.9 2.2 to 2.1 2.2 

\plot 
2 1 
2.125 1.05296
2.25 1.109
2.375 1.16056
2.5 1.216 
2.625 1.26056 
2.75 1.3009
2.875 1.35296 
3 1.4 /

\plot 
3 1.4 
3.125 1.45296
3.25 1.509
3.375 1.56056
3.5 1.616 
3.625 1.66056 
3.75 1.7009
3.875 1.75296 
4 1.8 /

\plot 
4 1.8 
4.125 1.85296
4.25 1.909
4.375 1.96056
4.5 2.016 
4.625 2.06056 
4.75 2.1009
4.875 2.15296 
5 2.2 /

\plot 
5 2.2 
5.125 2.25296
5.25 2.309
5.375 2.36056
5.5 2.416 
5.625 2.46056 
5.75 2.5009
5.875 2.55296 
6 2.6 /

\plot 
6 2.6 
6.125 2.65296
6.25 2.709
6.375 2.76056
6.5 2.816 
6.625 2.86056 
6.75 2.9009
6.875 2.95296 
7 3 /


\plot 
7 0 
6.875 -0.05296
6.75 -0.109
6.625 -0.16056
6.5 -0.216 
6.375 -0.26056 
6.25 -0.3009
6.125 -0.35296 
6 -0.4 /

\plot 
6 -0.4 
5.875 -0.45296
5.75 -0.509
5.625 -0.56056
5.5 -0.616 
5.375 -0.66056 
5.25 -0.7009
5.125 -0.75296 
5 -0.8 /

\plot 
5 -0.8 
4.875 -0.85296
4.75 -0.909
4.625 -0.96056
4.5 -1.016 
4.375 -1.06056 
4.25 -1.1009
4.125 -1.15296 
4 -1.2 /

\plot 
4 -1.2 
3.875 -1.25296
3.75 -1.309
3.625 -1.36056
3.5 -1.416 
3.375 -1.46056 
3.25 -1.5009
3.125 -1.55296 
3 -1.6 /

\plot 
3 -1.6 
2.875 -1.65296
2.75 -1.709
2.625 -1.76056
2.5 -1.816 
2.375 -1.86056 
2.25 -1.9009
2.125 -1.95296 
2 -2 /


\setdots

\putrule from 4 -2 to 4 3 
\putrule from 5 -2 to 5 3
\putrule from 2 0 to 7 0
\putrule from 2 1 to 7 1

\putrule from 2 -1.2 to 4 -1.2
\putrule from 2 -0.8 to 5 -0.8
\putrule from 2 1.8 to 4 1.8
\putrule from 2 2.2 to 5 2.2

\plot 2 -2 
7 3 /

\put{Figure 3} at 0 -2.4 
\put{} at -8.1 0

\small

\put{$a$} at -7 -0.5
\put{$f(c)\!=\!b$} at -5.2 -0.5 
\put{$c\!=\!h(b)$} at -3.8 -0.5
\put{$d$} at -2 -0.5
\put{$f$} at -4.5 1.3 
\put{$h$} at -3.6 0.5 
\put{$I$} at -4.5 0.3


\put{$a$} at 1.7 -2.3
\put{$b$} at 4 -2.3 
\put{$c$} at 5 -2.3
\put{$d$} at 7 -2.3

\put{$b$} at 1.5 0
\put{$c$} at 1.5 1
\put{$d$} at 1.7 3
\put{$I$} at 1.75 0.5
\put{$f$} at 4.7 -1.2
\put{$g$} at 4.7 1.8

\endpicture


\vspace{0.4cm}

Puisque les orbites du pseudo-groupe engendr\'e par $f$ et $g$ sont denses sur $\Lambda$, 
la preuve du th\'eor\`eme se r\'eduit \`a d\'emontrer la proposition suivante.

\vspace{0.05cm}

\begin{prop} {\em Sous les hypoth\`eses pr\'ecedentes, le pseudo-groupe engendr\'e 
par $f$ et $g$ contient des \'el\'ements avec des points fixes hyperboliques 
appartenant \`a l'ensemble $\Lambda$.}
\label{resorte}
\end{prop}

\noindent{\bf D\'emonstration.} Consid\`erons l'espace 
$\Omega = \{ f,g \}^{\mathbb{N}}$ muni de la mesure de Bernoulli 
$\mathbb{P}$ de poids $(1/2,1/2)$. Pour chaque $\omega = (g_1,g_2,\ldots) \in \Omega$ et 
chaque $n \in \mathbb{N}$ notons $h_n(\omega) = g_n \cdots g_1$, et posons $h_0(\omega)=id$. 
Puisque les intervalles de l'ensemble $\big\{ g_n \cdots g_1 (I): n \geq 0, \thinspace 
(g_1,\ldots,g_n) \in \{f,g\}^n \big\}$ sont deux \`a deux disjoints,
$$\sum_{n \geq 0} \sum_{(g_1,\ldots,g_n) \in \{ f, g \}^n} |g_n \cdots g_1 (I)| 
\leq d\!-\!a < \infty.$$
Or, par le th\'eor\`eme de Fubini,
$$\sum_{n \geq 0} \sum_{(g_1,\ldots,g_n) \in \{ f, g \}^n} |g_n \cdots g_1 (I)| 
= \sum_{n \geq 0} 2^n \Big( \int_{\Omega} |h_n(\omega)(I)| 
\thinspace d\mathbb{P}(\omega) \Big) = \int_{\Omega} 
\Big( \sum_{n \geq 0} 2^n |h_n(\omega)(I)| \Big) d\mathbb{P}(\omega).$$
Donc, pour $\mathbb{P}$-presque toute suite $\omega \in \Omega$, on a la 
convergence de la s\'erie
$$\sum_{n \geq 0} 2^n |h_n (\omega)(I)|.$$
Fixons $\varepsilon \! \in ]0,1[$, et pour chaque $B>0$ 
consid\'erons l'ensemble $\Omega(B,\varepsilon)$ d\'efini par
$$\Omega(B,\varepsilon) = \big\{ \omega \in \Omega: |h_n(\omega)(I)| \leq 
B/(2-\varepsilon)^n \mbox{ pour tout } n \geq 0 \big\}.$$
D'apr\`es ce qui pr\'ec\`ede, $\mathbb{P} [\Omega(B,\varepsilon)]$ converge vers $1$ 
lorsque $B$ tend vers l'infini. Nous pouvons en particulier fixer $B$ suffisamment 
grand de fa\c{c}on \`a ce que $\mathbb{P} [\Omega(B,\varepsilon)] > 0$. Remarquons 
que si $\omega$ appartient \`a $\Omega(B,\varepsilon)$ alors
\begin{equation}
\sum_{n \geq 0} |h_n(\omega)(I)|^{\tau} \leq 
B^{\tau} \sum_{n \geq 0} \frac{1}{(2-\varepsilon)^{n \tau}} = M < \infty.
\label{cota}
\end{equation}

Consid\'erons l'intervalle $J' = [b-L,c+L]$ contenant l'intervalle errant $I$, o\`u \thinspace 
$L \!=\! L(\tau,M,|I|;\{f,g\})$ \thinspace est la constante qui appara\^{\i}t dans le lemme \ref{?}. 
Si $N \in \mathbb{N}$ est suffisamment grand alors \esp $f^N g$ \esp et \esp $g^N f$ \esp envoient 
tout l'intervalle $[0,1]$ sur $J' \setminus I$. La loi $0\!-\!1$ de Borel donne ainsi 
$$\mathbb{P}[h_n(\omega)(I) \subset J' \setminus I \mbox{ une infinit\'e de fois}] = 1.$$
Si $\omega \in \Omega(B,\varepsilon)$ et $n \in \mathbb{N}$ satisfont 
$h_n(\omega)(I) \subset J' \setminus I$, alors le  lemme \ref{?} montre que $h_n(\omega)$ 
poss\`ede un point fixe hyperbolique. Finalement, puisque l'ensemble $\Lambda$ 
est invariant par le pseudo-groupe et que le point fixe que l'on 
a trouv\'e attire une partie de cet ensemble, ce point fixe 
appartient n\'ecessairement \`a $\Lambda$. \esp $\square$ 

\begin{rem} Si l'on examine la preuve pr\'ec\'edente de plus pr\`es, on s'aper\c{c}oit 
que la m\^eme technique permet de d\'emontrer que pour presque tout $\omega \in \Omega$ 
il existe une infinit\'e d'entiers $n \in \mathbb{N}$ tels que $h_n(\omega)$ poss\`ede 
un point fixe hyperbolique $p_n \in \Lambda$ de telle sorte que, asymptotiquement, 
la contraction $h_n(\omega)' (p_n)$ est au moins de l'ordre de $1/(2-\varepsilon)^n$ 
pour tout $\varepsilon \!>\! 0$. 
\end{rem}

\vspace{0.05cm}

\begin{rem} En classe $C^{1+lip}$, la proposition \ref{resorte} d\'ecoule directement du th\'eor\`eme 
de Sacksteder classique. Soulignons cependant que, dans ce contexte, on dispose d'une version 
beaucoup plus fine. En effet, consid\'erons un intervalle $]\bar{a},\bar{d}[$ contenant 
strictement $[a,d]$, et fixons une application $T: [\bar{a},\bar{d}] \rightarrow [\bar{a},\bar{d}]$ 
de classe $C^{1+lip}$, n'ayant que deux points critiques, et dont la restriction \`a $[a,b]$ 
(resp. $[c,d]$) co\"{\i}ncide avec $f^{-1}$ (resp. $g^{-1}$). Il r\'esulte alors du th\'eor\`eme 
d'hyperbolicit\'e de Ma\~n\'e \cite{man} que tous les points p\'eriodiques de $T$ appartenant \`a 
l'ensemble de Cantor $\Lambda$ et dont la p\'eriode est suffisamment grande sont hyperboliques 
(dilatants). Or, pour tout $n \in \mathbb{N}$ la restriction de $T^n$ \`a $\Lambda$ co\"{\i}ncide 
localement avec des \'el\'ements du pseudo-groupe engendr\'e par $f$ et $g$.
\end{rem}

\vspace{0.05cm}

Nous finissons ce paragraphe en donnant le sch\'ema de la preuve d'une 
version raffin\'ee de l'un des r\'esultats contenus dans \cite{GT}.

\vspace{0.05cm}

\begin{prop} {\em Soient $\Gamma_1$ et $\Gamma_2$ deux groupes de diff\'eomorphismes 
du cercle dont la classe de diff\'erentiabi- lit\'e $C^{\alpha}$ est strictement 
sup\'erieure \`a $C^1$. Supposons qu'ils ne poss\`edent pas d'orbite finie et 
qu'ils ne soient pas semi-conjugu\'es \`a des groupes de rotations. Alors tout 
diff\'eomorphisme de classe $C^1$ qui conjugue $\Gamma_1$ avec $\Gamma_2$ 
est n\'ecessairement de classe $C^{\alpha}$.}
\label{gh-ts-tau}
\end{prop}

\noindent{\bf D\'emonstration.} Supposons que $\varphi$ soit un diffeomorphisme de 
classe $C^1$ qui conjugue $\Gamma_1$ avec $\Gamma_2$. Les hypoth\`eses faites sur 
ces groupes impliquent, par des arguments bien connus (voir par exemple les appendices 
\ref{ap0} et \ref{ap1}), qu'ils contiennent des \'el\'ements v\'erifiant les hypoth\`eses 
de la proposition \ref{resorte}. Par suite, ces groupes poss\`edent des \'el\'ements avec 
des points fixes hyperboliques. De plus, lorsqu'il existe un minimal exceptionnel, 
de tels points fixes peuvent \^etre pris appartenant \`a cet ensemble minimal. Le 
diff\'eomorphisme $\varphi$ \'etant de classe $C^1$, il envoie des points fixes 
hyperboliques sur des points fixes hyperboliques. En particulier, il conjugue les 
germes correspondants. D'apr\`es le th\'eor\`eme de lin\'earisation de Sternberg, 
l'application $\varphi$ est localement un diff\'eomorphisme de classe $C^{\alpha}$ 
autour de ces points fixes (voir la note au bas de la page 17). Puisque 
l'ensemble des points autour desquels $\varphi$ est un diff\'eomorphisme local 
de classe $C^{\alpha}$ est invariant par $\Gamma_1$, lorsque les orbites sont denses 
cela suffit pour d\'emontrer que $\varphi$ est un diff\'eomorphisme de classe 
$C^{\alpha}$ sur tout le cercle. C'est le cas aussi lorsqu'il y a un minimal 
exceptionnel, car toute orbite s'accumule sur ce minimal. 
Ceci termine la d\'emonstration. \esp $\square$


\subsection{Pseudo-groupes de diff\'eomorphismes et feuilletages : le cas $C^{1}$}
\label{sac-ps2}

M\^eme en pr\'esence d'une dynamique dilatante, il est en g\'en\'eral 
impossible de contr\^oler les distorsions en classe $C^1$ (une 
illustration classique de ce fait appara\^{\i}t dans \cite{bowen}). 
N\'eanmoins, lorsqu'on sait {\em a priori} qu'il y a de l'hyperbolicit\'e 
quelque part, des arguments \'el\'ementaires permettent d'en d\'eduire 
qu'il y a de l'hyperbolicit\'e partout : c'est l'id\'ee sous-jacente de la 
preuve que nous donnons dans la suite \footnote{Il est possible 
que par une m\'ethode similaire on puisse traiter le cas critique 
({\em i.e.} en classe $C^{1+1\!/\!d}$) de nos r\'esultats pour des 
groupes ab\'eliens de diff\'eomorphismes.}. Une formulation plus 
conceptuelle de cette id\'ee (en termes d'exposants de Lyapunov et de 
mesures stationnaires) sera esentielle au \S \ref{morse-smale}.

Tout en gardant les notations du \S \ref{sac-ps}, fixons une fois pour toutes 
une constante $\varepsilon \!\in ]0,1/3[$. Nous savons que pour $\mathbb{P}$-presque 
tout $\omega \in \Omega$ il existe $B = B(\omega) \geq 1$ tel que
\begin{equation}
|h_n(\omega)(I)| \leq \frac{B}{(2-\varepsilon)^n} 
\quad \mbox{ pour tout} \quad n \geq 0.
\label{int-exp}
\end{equation}


\begin{lem} {\em Il existe une constante $\bar{C}$ ne d\'ependant que de $f$ et $g$ 
telle que, si $\omega = (g_1,g_2,\ldots) \in \Omega$ satisfait} (\ref{int-exp}), 
{\em alors pour tout $x \in I$ et tout entier $n \geq 0$ on a}
\begin{equation}
h_n(\omega)'(x) \leq \frac{B \bar{C}}{(2-2\varepsilon)^n}.
\label{hyp-inf}
\end{equation}
\label{aaa}
\end{lem}
\vspace{-0.45cm}
\noindent{\bf D\'emonstration.} Fixons $\varepsilon_0 > 0$ suffisamment petit de 
fa\c{c}on \`a ce que pour tout $y,z$ dans $[a,d]$ \`a distance inf\'erieure ou 
\'egale \`a $\varepsilon_0$ on ait
\begin{equation}
\frac{f'(y)}{f'(z)} \leq \frac{2-\varepsilon}{2-2\varepsilon} \qquad 
\mbox{et} \qquad \frac{g'(y)}{g'(z)} \leq \frac{2-\varepsilon}{2-2\varepsilon}.
\label{continuity}
\end{equation}
Il est facile de voir qu'il existe $N \in \mathbb{N}$ tel que, pour tout 
$\omega \in \Omega$ et tout $i \geq 0$, la taille de l'intervalle 
$h_{N+i}(I)$ est inf\'erieure ou \'egale \`a $\varepsilon_0$. Nous affirmons 
alors que (\ref{hyp-inf}) a lieu pour $\bar{C} = \max\{ A,\bar{A} \}$, o\`u 

$$A = \sup_{x\in I, n\leq N, \omega \in \Omega} \frac{h_n(\omega)'(x) \esp (2-2\varepsilon)^n}{B}, 
\qquad \bar{A} = \sup_{x,y \in I, \omega \in \Omega} \frac{h_N(\omega)'(x)}
{h_N(\omega)'(y) \esp |I|} \esp \Big( \frac{2-2\varepsilon}{2-\varepsilon} \Big)^N.$$
En effet, si $n \leq N$ alors (\ref{hyp-inf}) a lieu \`a cause de l'in\'egalit\'e 
\esp $\bar{C} \geq A$. \esp Supposons donc $n > N$, et fixons $y = y(n) \in I$ 
tel que \esp $|h_n(\omega)(I)| = h_n(\omega)'(y) \esp |I|$. Si $x$ appartient 
\`a $I$ alors la distance entre les points $h_{N+i}(\omega)(x)$ et 
$h_{N+i}(\omega)(y)$ est inf\'erieure ou \'egale \`a $\varepsilon_0$ 
pour tout $i \geq 0$. Par (\ref{continuity}),
$$\frac{h_n(\omega)'(x)}{h_n(\omega)'(y)} = \frac{h_N(\omega)'(x)}{h_N(\omega)'(y)} \esp 
\frac{g_{N+1}'(h_N(\omega)(x))}{g_{N+1}'(h_N(\omega)(y))} \cdots 
\frac{g_n'(h_{n-1}(\omega)(x))}{g_n'(h_{n-1}(\omega)(y))} 
\leq \frac{h_N(\omega)'(x)}{h_N(\omega)'(y)} \esp 
\Big( \frac{2-\varepsilon}{2-2\varepsilon} \Big)^{n-N},$$
et donc
$$h_n(\omega)'(x) 
\leq \frac{h_N(\omega)'(x)}{h_N(\omega)'(y)} \esp \frac{|h_n(\omega)(I)|}{|I|} \esp 
\Big( \frac{2-\varepsilon}{2-2\varepsilon} \Big)^{n-N} 
\leq \frac{h_N(\omega)'(x)}{h_N(\omega)'(y)} \esp \frac{B}{|I|\esp(2-\varepsilon)^n} 
\esp \Big( \frac{2-\varepsilon}{2-2\varepsilon} \Big)^{n-N} \leq 
\frac{B\bar{C}}{(2-2\varepsilon)^n},$$ 
o\`u la derni\`ere in\'egalit\'e d\'ecoule de la condition 
\esp $\bar{C} \geq \bar{A}$. \esp $\square$

\vspace{0.2cm}

Pour obtenir de l'hyperbolicit\'e au del\`a de l'intervalle $I$ on doit 
utiliser un argument 
\hspace{0.05cm}\begin{tiny}$^{_\ll}$\end{tiny}dual\hspace{0.05cm}\begin{tiny}$^{_\gg}$\end{tiny} 
mais l\'eg\`erement plus \'elabor\'e 
que celui du lemme pr\'ec\'edent. Pour le formuler fixons une 
constante $\varepsilon_1 > 0$ suffisamment petite de sorte 
que pour tout $y,z$ dans $[a,d]$ \`a distance inf\'erieure 
ou \'egale \`a $\varepsilon_1$ on ait
\begin{equation}
\frac{f'(y)}{f'(z)} \leq \frac{2-2\varepsilon}{2-3\varepsilon} \qquad 
\mbox{et} \qquad \frac{g'(y)}{g'(z)} \leq \frac{2-2\varepsilon}{2-3\varepsilon}.
\label{continuity-2}
\end{equation}

\vspace{0.02cm}

\begin{lem} {\em Soient $C \geq 1$, $\omega = (g_1,g_2,\ldots) \in \Omega$ 
et $x \in [a,d]$ tels que}
\begin{equation}
h_n(\omega)'(x) \leq \frac{C}{(2-2\varepsilon)^n} \quad \mbox{ pour tout} 
\quad n \geq 0.
\label{hypot}
\end{equation}
{\em Alors pour tout point $y \in [a,d]$ \`a distance inf\'erieure ou \'egale 
\`a \esp $\varepsilon_1 / C$ \esp de $x$ et tout $n \geq 0$ on a}
\begin{equation}
h_n(\omega)'(y) \leq \frac{C}{(2 - 3\varepsilon)^n}.
\label{tesis}
\end{equation}
\label{bbb}
\end{lem}
\vspace{-0.5cm}
\noindent{\bf D\'emonstration.} La v\'erification de l'in\'egalit\'e (\ref{tesis}) 
se fait par r\'ecurrence. Pour $n = 0$ elle a lieu \`a cause de l'hypoth\`ese 
$C \geq 1$. Admettons qu'elle soit valable pour tout $j \in \{0,\ldots,n\}$, et 
notons $y_j = h_{j}(\omega)(y)$ et $x_j = h_j(\omega)(y)$. Supposons que 
$y \leq x$, l'autre cas \'etant analogue. Chaque point $y_j$ appartient 
alors \`a l'intervalle \esp $h_j(\omega)([x-\varepsilon_1/C,x])$. \esp Or, 
par hypoth\`ese de r\'ecurrence,
$$\big| h_j(\omega)([x-\varepsilon_1/C,x]) \big| \leq \frac{C}{(2-3\varepsilon)^j} \esp 
\big| [x-\varepsilon_1/C,x] \big| \leq C \esp \frac{\varepsilon_1}{C} = \varepsilon_1.$$
Par la d\'efinition de $\varepsilon_1$ on conclut que, pour tout $j \leq n$ on a \esp 
$g_{j+1}'(y_j) \leq g_{j+1}'(x_j) \esp \Big( \frac{2-2\varepsilon}{2-3\varepsilon} \Big)$. 
\esp Donc, d'apr\`es l'hypoth\`ese (\ref{hypot}),
\begin{eqnarray*}
h_{n+1}(\omega)'(y) 
&=& g_1'(y_0) \cdots g_{n+1}'(y_n) \esp\esp \leq \esp\esp g_1'(x_0) \cdots g_{n+1}'(x_n) \esp 
\Big( \frac{2-2\varepsilon}{2-3\varepsilon} \Big)^{n+1}\\
&\leq& \esp \esp \frac{C}{(2-2\varepsilon)^{n+1}} \esp 
\Big( \frac{2-2\varepsilon}{2-3\varepsilon} \Big)^{n+1} \esp \esp 
\leq \esp \esp \frac{C}{(2-3\varepsilon)^{n+1}}.
\end{eqnarray*}
Ceci ach\`eve la v\'erification par r\'ecurrence de (\ref{tesis}). \esp $\square$

\vspace{0.2cm}

Nous sommes maintenant en mesure de terminer la preuve du th\'eor\`eme C 
en classe $C^1$. Pour cela remarquons que, d'apr\`es le lemme \ref{aaa}, 
si $C$ est suffisamment grand alors la probabilit\'e de l'ensemble
$$\Omega(C) = \Big\{ \omega \in \Omega : h_n(\omega)'(x) \leq \frac{C}{(2-2\varepsilon)^n} 
\quad \mbox{pour tout } \esp n \geq 0 \esp \mbox{ et tout } \esp x \in I \Big\}$$
est strictement positive. Fixons un tel $C \geq 1$ et notons $L = \min \{\varepsilon_1/2 C, |I|/2 \}$. 
Si l'on d\'esigne par $J$ le $2L$-voisinage de $I$, alors le lemme \ref{bbb} entra\^{\i}ne 
que pour tout $\omega \in \Omega(C)$, tout $n \geq 0$ et tout $y \in J$,

\begin{equation}
h_n(\omega)'(y) \leq \frac{C}{(2-3\varepsilon)^n}.
\label{mato}
\end{equation}
Si l'on d\'esigne par $J'$ le $L$-voisinage de $I$, alors \esp\esp\esp 
$\mathbb{P} \big[ h_n(\omega)(I) \subset J' \setminus I 
\esp \mbox{ une infinit\'e de fois} \big] \!=\! 1.$ \esp\esp\esp  
Par cons\'equent, et \'etant donn\'e que $\varepsilon < 1/3$, il existe 
$\omega \in \Omega(C)$ et $m \in \mathbb{N}$ tels que $h_m(\omega)(I) 
\subset J' \setminus I$ et $(2-3\varepsilon)^m > C$. D'apr\`es (\ref{mato}), 
et puisque $L \leq |I|/2$, cela implique que $h_m(\omega)$ envoie $J$ sur 
l'une des deux composantes connexes de $J \setminus I$ de telle mani\`ere 
que $h_m(\omega)$ poss\`ede un point fixe dans cette composante connexe, 
lequel est n\'ecessairement hyperbolique (et appartient \`a $\Lambda$).


\subsection{Groupes de diff\'eomorphismes du cercle : m\'ethode d\'eterministe}
\label{sac-circ}

Avant de nous plonger dans la preuve du th\'eor\`eme D, 
nous voudrions montrer que par des m\'ethodes 
\begin{tiny}$^{_\ll}$\end{tiny}d\'eterministes\hspace{0.05cm}\begin{tiny}$^{_\gg}$\end{tiny} 
on peut d\'ej\`a obtenir une am\'elioration du th\'eor\`eme de Sacksteder concernant 
le nombre de points fixes hyperboliques.

\vspace{0.1cm}

\begin{prop} {\em Soit $\Gamma$ un sous-groupe de $\mathrm{Diff}_+^{1+lip}(\clo)$. Si 
$\Gamma$ ne pr\'eserve aucune mesure de probabilit\'e du cercle, alors il contient 
des \'el\'ements avec (au moins) deux points fixes hyperboliques, l'un 
contractant et l'autre dilatant.}
\label{sac-det}
\end{prop}

\vspace{0.05cm}

Remarquons que ce r\'esultat est tr\`es naturel, car le th\'eor\`eme C fournit d\'ej\`a 
un \'el\'ement avec un point fixe hyperbolique, et tout diff\'eomorphisme du cercle 
avec un tel point poss\`ede n\'ecessairement d'autres points fixes. La proposition 
ci-dessus stipule donc l'existence d'\'el\'ements pour lesquels au moins l'un 
de ces autres points fixes est hyperbolique. Pour aboutir \`a une 
d\'emonstration, nous chercherons \`a v\'erifier l'existence d'une 
suite $(h_n)$ d'\'el\'ements dans $\Gamma$ telle que $h_n$ et 
$h_n^{-1}$ satisfassent les hypoth\`eses du lemme \ref{!}.

\vspace{0.1cm}


\beginpicture

\setcoordinatesystem units <1cm,1cm>

\circulararc 360 degrees from 3 0 center at 0 0

\circulararc -28 degrees from 3.25 0.8 center at 0 0

\circulararc -40 degrees from 3.4 -1.8 center at 0 0

\circulararc 45 degrees from 3.4 1.8 center at 0 0

\put{$\bullet$} at 2.89 0.8  
\put{$\bullet$} at 2.89 -0.8 
\put{$a$} at 3.1 -1.05
\put{$b$} at 3.1 1.05 
\put{$I$} at 3.55 0
\put{$f^{-1}$} at 2.8 3.2 
\put{$f$} at 2.8 -3.1 
\put{$a'$} at -1.2 -2.2
\put{$b'$} at -1.2 2.2
\put{$g$} at 1.5 1
\put{$g$} at 1.5 -1
\put{$K$} at 2.05 2.9 

\plot -1.4 -2.8 
-1.2 -2.6 /

\plot -1.4 -2.6 
-1.2 -2.8 /

\plot -1.4 2.8 
-1.2 2.6 /

\plot -1.4 2.6 
-1.2 2.8 /

\plot 3.41 -1.78
3.34 -2.15 /

\plot 3.41 -1.78
3.1 -2.05 /

\plot 3.41 1.78
3.34 2.15 /

\plot 3.41 1.78
3.1 2.05 /

\plot 2.5 0.6 
2.15 0.65 /

\plot 2.5 0.6 
2.25 0.85 /

\plot 1.5 -1.98  
1.4 -1.67 /

\plot 1.5 -1.98 
1.7 -1.7 /

\plot 2 2 
2.25 2.25 /

\plot 1.2 2.55
1.385 2.85 /

\plot 2 -2 
2.25 -2.25 /

\plot 1.2 -2.55
1.385 -2.85 /

\circulararc -65 degrees from 2.5 0.6 center at 3 2

\circulararc 70 degrees from 2.5 -0.6 center at 3 -2

\circulararc -17 degrees from 1.53 2.95 center at 0 0 

\circulararc 17 degrees from 1.53 -2.95 center at 0 0 


\put{Figure 4} at 0 -3.6 
\put{} at -8.1 4.18

\endpicture

\vspace{0.1cm}

\begin{lem} {\em Soit $\Gamma$ un groupe de diff\'eomorphismes de classe $C^{1+lip}$ 
du cercle qui admet un ensemble invariant ferm\'e et non vide $\Lambda \neq \clo$. 
Supposons qu'il existe un \'el\'ement $f \in \Gamma$ et une composante connexe $I\!= ]a,b[$ 
de $\clo \setminus \Lambda$ tels que pour certains points $a'$ et $b'$ dans 
$\clo\setminus [a,b]$ on ait $\lim_{n \rightarrow \infty} f^n(x)=a$ pour tout $x\!\in ]a',a[$ 
et $\lim_{n \rightarrow \infty} f^{-n}(y) = b$ pour tout $y \!\in ]b,b'[$. Supposons 
aussi que $\Gamma$ contienne un \'el\'ement $g$ tel que $g(I) \!\subset ]a',a[$ 
et $g^{-1}(I) \!\subset ]b,b'[$. Alors pour $n$ assez grand, l'\'el\'ement 
$f^n g \in \Gamma$ poss\`ede un point fixe hyperbolique contractant et 
un autre point fixe hyperbolique dilatant.}
\label{cc}
\end{lem}

\noindent{\bf D\'emonstration.} Par hypoth\`ese, pour $\varepsilon\!>\!0$ assez petit 
les intervalles $]a,a\!+\!\varepsilon[=\!J,g(J),fg(J),\ldots,f^ng(J)$ sont deux \`a 
deux disjoints et convergent vers l'extr\'emit\'e gauche de $I$ (voir la figure 4). 
D'apr\`es le lemme \ref{!}, pour $n$ assez grand l'application $h_{n+1} = f^n g$ 
poss\`ede un point fixe hyperbolique contractant proche de $a$. Consid\'erons 
maintenant l'intervalle $K\!=\!g^{-1}(]b-\varepsilon,b[)$. Encore par hypoth\`ese, 
si $\varepsilon\!>\!0$ est assez petit alors les intervalles 
$K,f^{-1}(K),\ldots,f^{-n}(K)$ et $g^{-1}f^{-n}(K)$ sont deux \`a deux disjoints 
et convergent vers l'extr\'emit\'e droite de $K$. Comme pr\'ec\'edemment, 
$h_{n+1}^{-1} = g^{-1}f^{-n}$ poss\`ede un point fixe hyperbolique proche de 
(l'extr\'emit\'e droite de) $K$ d\`es que $n$ est assez grand. \esp $\square$

\vspace{0.35cm}

Si l'on est capable de s\'eparer les ensembles de points fixes de deux \'el\'ements d'un 
sous-groupe $\Gamma$ de $\mathrm{Diff}_+^{1+lip}(\clo)$, alors le lemme pr\'ec\'edent 
permet d'obtenir des \'el\'ements dans $\Gamma$ avec deux points fixes hyperboliques.

\vspace{0.1cm}

\begin{lem} {\em Soit $\Gamma$ un sous-groupe de $\mathrm{Diff}_+^{1+lip}(\clo)$. 
Supposons que $\Gamma$ contienne deux \'el\'ements $g_1$ et $g_2$ qui poss\`edent 
au moins deux points fixes et tels qu'il existe deux intervalles ouverts $U$ et $V$ 
dans $\clo$ contenant $Fix(g_1)$ et $Fix(g_2)$ respectivement et dont les fermetures 
sont disjointes. Alors $\Gamma$ contient un \'el\'ement avec un point fixe 
hyperbolique contractant et un autre point fixe hyperbolique dilatant.}
\label{pp}
\end{lem}

\noindent{\bf D\'emonstration.} Soient $p_1,q_1$ (resp. $p_2,q_2$) les points fixes de $g_1$ 
(resp. $g_2$) qui determinent un intervalle ferm\'e contenant $Fix(g_1)$ (resp. $Fix(g_2)$), 
et soient $U \!= ]u,u'[$ et $V\! = ]v,v'[$ les intervalles donn\'es par l'hypoth\`ese. 
Quitte \`a changer $g_1$ et $g_2$ par leurs inverses, et quitte \`a prendre des it\'er\'es 
suffisamment grands, nous pouvons supposer que $g_2([q_2,u']) \!\subset ]q_2,v']$ et 
$g_1([v,p_1]) \!\subset ]u,p_1]$ (voir la  figure 5). Soit $\Lambda$ l'ensemble ferm\'e 
non vide invariant et minimal pour le groupe engendr\'e par $g_1$ et $g_2$. Un argument 
facile de type ping-pong montre que $\Lambda$ n'est pas tout le cercle ; il est en fait 
hom\'eomorphe \`a l'ensemble de Cantor.

Consid\'erons maintenant les composantes connexes $I\!= ]a,b[$ et $K$ du compl\'ementaire 
de $\Lambda$ qui contiennent $]v',u[$ et $]u',v[$ respectivement. Notons que $g_1(K)=g_2(K)=I$. 
En particulier, l'application $f = g_2 g_1^{-1}$ fixe $a$ et $b$. Nous affirmons que 
ce sont les seuls points fixes de $f$ sur $[b,a]$. En effet, on v\'erifie ais\'ement 
que $x \!< \!f(x) \!< \!a$ pour tout $x \!\in ]b,a[$. Pour finir la preuve de la 
proposition, il suffit de remarquer que pour $g = g_2$ et pour $f$, $I$ et $K$ 
comme pr\'ec\'edemment, les hypoth\`eses du lemme \ref{cc} sont satisfaites. 
\esp $\square$

\vspace{0.4cm}


\vspace{0.65cm}


\beginpicture

\setcoordinatesystem units <1cm,1cm>

\circulararc 360 degrees from 4 0 center at 0 0
\circulararc 307 degrees from 2.7 1.4 center at 0 0 

\circulararc -177.8 degrees from 4.8 1.3 center at 0 0 
\circulararc 176.8 degrees from 4.57 -1.25 center at 0 0 

\put{$\bullet$} at 1.85 3.54
\put{$\bullet$} at 1.85 -3.54
\put{$\bullet$} at -1.85 -3.54
\put{$\bullet$} at -1.85 3.54

\put{$*$} at 3.9 0.9
\put{$*$} at 3.9 -0.9
\put{$*$} at -3.9 -0.9
\put{$*$} at -3.9 0.9

\plot 2.75 -1.3 
2.57 -1.45 /

\plot 2.75 -1.3 
2.75 -1.5 /

\plot 4.8 1.3 
4.78 1.15 /

\plot 4.8 1.3 
4.907 1.18 /

\plot 4.57 -1.25
4.54 -1.07 /

\plot 4.57 -1.25 
4.71 -1.1 /

\put{$I$} at 3.75 0 
\put{$K$} at -3.77 0
\put{$f$} at 0 -2.8
\put{$a$} at 3.4 -1.2
\put{$b$} at 3.4 1.25
\put{$u$} at 4.22 0.9
\put{$v$} at -4.22 -0.9
\put{$u'$} at -4.25 0.93
\put{$v'$} at 4.25 -0.83
\put{$p_1$} at 1.7 3.25 
\put{$q_1$} at -1.7 3.25 
\put{$p_2$} at -1.7 -3.25 
\put{$q_2$} at 1.7 -3.25 
\put{$g_1$} at 0 -4.75
\put{$g_2$} at 0 4.52

\circulararc 35 degrees from 3.4 -1.07 center at 0 0 
\circulararc -35 degrees from -3.4 -1.07 center at 0 0 

\circulararc 6.3 degrees from -3.78 -1.24 center at 0 0 
\circulararc 6.3 degrees from -3.82 -1.24 center at 0 0 

\circulararc 6.3 degrees from 3.78 1.24 center at 0 0 
\circulararc 6.3 degrees from 3.82 1.24 center at 0 0 

\circulararc -6.3 degrees from -3.78 1.24 center at 0 0 
\circulararc -6.3 degrees from -3.82 1.24 center at 0 0 

\circulararc -6.3 degrees from 3.78 -1.24 center at 0 0 
\circulararc -6.3 degrees from 3.82 -1.24 center at 0 0 

\circulararc -6 degrees from 1.85 3.523 center at 0 0 
\circulararc -6 degrees from 1.85 3.56 center at 0 0 

\circulararc 6 degrees from 1.85 -3.523 center at 0 0 
\circulararc 6 degrees from 1.85 -3.56 center at 0 0 

\circulararc 6 degrees from -1.85 3.523 center at 0 0 
\circulararc 6 degrees from -1.85 3.56 center at 0 0 

\circulararc -6 degrees from -1.85 -3.523 center at 0 0 
\circulararc -6 degrees from -1.85 -3.56 center at 0 0 

\circulararc 6.7 degrees from 3.21 2.42 center at 0 0 
\circulararc 6.7 degrees from 3.2 2.375 center at 0 0 

\circulararc -6.7 degrees from 3.2 -2.375 center at 0 0 
\circulararc -6.7 degrees from 3.21 -2.42 center at 0 0 

\circulararc -6.7 degrees from -3.20 2.375 center at 0 0 
\circulararc -6.7 degrees from -3.21 2.42 center at 0 0 

\circulararc 6.7 degrees from -3.20 -2.375 center at 0 0 
\circulararc 6.7 degrees from -3.21 -2.42 center at 0 0 

\plot 3.7 1.2 
3.885 1.26 /

\plot -3.7 1.2 
-3.885 1.26 /

\plot 3.7 -1.2 
3.885 -1.26 /

\plot -3.7 -1.2 
-3.885 -1.26 /
\plot 3.55 1.62
3.727 1.701 /

\plot 3.55 -1.62
3.727 -1.701 /

\plot -3.55 1.62
-3.727 1.701 /

\plot -3.55 -1.62
-3.727 -1.701 /
\plot 3.1 2.3 
3.32 2.45 /

\plot -3.1 2.3 
-3.32 2.45 /

\plot 3.1 -2.3 
3.32 -2.45 /

\plot -3.1 -2.3 
-3.32 -2.45 /
\plot 2.8 2.7 
2.98 2.86 /

\plot -2.8 -2.7 
-2.98 -2.86 /

\plot -2.8 2.7 
-2.98 2.86 /

\plot 2.8 -2.7 
2.98 -2.86 /
\plot 2.13 3.21 
2.3 3.42 /

\plot -2.13 -3.21 
-2.3 -3.42 /

\plot -2.13 3.21 
-2.3 3.42 /

\plot 2.13 -3.21 
2.3 -3.42 /
\plot 1.78 3.45 
1.93 3.67 /

\plot -1.78 3.45 
-1.93 3.67 /

\plot -1.78 -3.45 
-1.93 -3.67 /

\plot 1.78 -3.45 
1.93 -3.67 /


\put{Figure 5} at 0 -5.6 
\put{} at -8.1 4.3

\endpicture


\vspace{0.65cm}

Pour finir la preuve de la proposition \ref{sac-det} 
on pourrait essayer de trouver deux \'el\'ements v\'erifiant les hypoth\`eses 
du lemme pr\'ec\'edent. N\'eanmoins, m\^eme pour des groupes qui ne pr\'eservent 
pas de mesure de probabilit\'e du cercle, de tels \'el\'ements peuvent ne pas exister. 
C'est le cas par exemple des rev\^etements finis de sous-groupes non 
m\'etab\'eliens de $\mathrm{PSL}(2,\mathbb{R})$. L'id\'ee consiste alors \`a 
montrer que ces \begin{tiny}$^{_\ll}$\end{tiny}extensions 
finies~\begin{tiny}$^{_\gg}$\end{tiny} 
sont les seules obstructions \`a la s\'eparation d'ensembles de points 
fixes. Pour cela nous utilisons quelques arguments de la preuve de Ghys de 
l'alternative de Tits faible pour les groupes d'hom\'eomorphismes 
du cercle (voir l'appendice \ref{ap1}).

\vspace{0.45cm}

\noindent{\bf D\'emonstration de la proposition \ref{sac-det}.} Puisque $\Gamma$ 
ne pr\'eserve pas de mesure de probabilit\'e, il n'a pas d'orbite finie. Donc, 
soit toutes ses orbites sont denses, soit il pr\'eserve un ensemble 
ferm\'e et minimal hom\'eomorphe \`a l'ensemble de Cantor. Dans la 
suite nous ne consid\'ererons que le cas minimal ({\em i.e.} lorsque 
toutes les orbites sont denses), et nous laisserons au lecteur le soin 
d'adapter les arguments ci-dessous au cas o\`u il existe un ensemble minimal 
exceptionnel (pour cela on utilise la technique classique qui consiste \`a 
\begin{tiny}$^{_\ll}$\end{tiny}\'ecraser\hspace{0.05cm}\begin{tiny}$^{_\gg}$\end{tiny} 
les composantes connexes du compl\'ementaire de cet ensemble).

D'apr\`es le \S \ref{ap1}, il existe un 
rev\^etement (topologique) fini $\pi: \clo \rightarrow \clo/\!\sim$ tel que l'action 
de $\Gamma$ sur $\clo$ induit une action de $\Gamma$ par hom\'eomorphismes du cercle 
topologique $\clo/\!\sim$. De plus, si l'on d\'esigne par $\hat{\Gamma}$ le sous-groupe 
correspondant de $\mathrm{Homeo}_+(\clo/\!\sim)$, alors la 
{\em propri\'et\'e d'expansivit\'e forte} suivante est v\'erifi\'ee : pour 
chaque sous-intervalle ferm\'e $[a,b]$ de $\clo/\!\sim$ il existe une suite 
$(\hat{h}_n)$ dans $\hat{\Gamma}$ telle que $\hat{h}_n([a,b])$ converge vers un seul point. 
On conclut en particulier qu'il existe un \'el\'ement non trivial $\hat{g}_1 \in \hat{\Gamma}$ 
poss\'edant au moins deux points fixes. En effet, en prenant $x \in \clo/\!\sim$, un petit 
intervalle $]c,d[$ contenant $x$, et $y \notin [c,d]$, il existe une suite $(\hat{h}_n)$ dans 
$\hat{\Gamma}$ telle que les compl\'ementaires de $\hat{h}_n([c,d])$ convergent vers le point 
$y$, et ceci implique que pour $n$ assez grand l'\'el\'ement $\hat{h}_n$ a un point fixe dans 
$]c,d[$ et un autre point fixe proche de $y$.

Soit $\hat{U}$ un intervalle ouvert contenant l'ensemble des points fixes de 
$\hat{g}_1$ et dont la fermeture ne soit pas tout le cercle $\clo/\!\sim$. 
En vertu des propri\'et\'es de minimalit\'e et d'expansivit\'e forte 
de $\hat{\Gamma}$, il existe $\hat{h} \in \hat{\Gamma}$ tel que 
$\hat{h}(\hat{U}) \cap \hat{U} = \emptyset$. Soient $\hat{V} = \hat{h}(\hat{U})$ 
et $\hat{g}_2 = \hat{h} \hat{g}_1 \hat{h}^{-1}$. Tous les arguments topologiques 
de la preuve du lemme \ref{pp} peuvent \^etre appliqu\'es au sous-groupe de 
$\mathrm{Homeo}_+(\clo/\!\sim)$ engendr\'e par (les it\'er\'es correspondants de) 
$\hat{g}_1$ et $\hat{g}_2$. Ceci donne en particulier un intervalle ouvert $\hat{I}$ 
tel que les seuls points fixes de l'\'el\'ement $\hat{g}_2 \hat{g}_1^{-1} \in \hat{\Gamma}$ 
dans le compl\'ementaire de $\hat{I}$ ce sont ses extr\'emit\'es...

D\'esignons par $\kappa \in \mathbb{N}$ le degr\'e du rev\^etement 
$\pi: \clo \rightarrow \clo/\!\sim$, et notons $\phi$ le morphisme correspondant 
de $\Gamma$ vers $\hat{\Gamma}$. Consid\'erons deux \'el\'ements $g_1$ et $g_2$ 
de $\Gamma$ tels que $\phi(g_1) = \hat{g}_1$ et $\phi(g_2) = \hat{g}_2$, et soit 
$f = (g_2 g_1^{-1})^{\kappa}$. Cet \'el\'ement $f \in \Gamma$ fixe toutes les pr\'eimages 
de $\hat{I}$ par $\pi$. Notons $I\!=\!I_1,I_2,\ldots,I_{\kappa}$ ces pr\'eimages, tout 
en respectant leur ordre cyclique sur le cercle. Soit $\bar{g} \in \Gamma$ 
tel que $\phi(\bar{g}) = \hat{g}_2$, et soient $g = \bar{g}^{\kappa}$ et $K = g^{-1}(I)$. 
Remarquons que $K$ co\"{\i}ncide avec la pr\'eimage par $\phi$ de 
$\hat{g}_2^{-\kappa}(\hat{I})$ qui est situ\'ee entre $I\!=\!I_1$ 
et $I_2$ (voir la figure 6 pour une illustration du cas 
$\kappa \! = \! 3$). On v\'erifie rapidement que les \'el\'ements $f$ et $g$ de $\Gamma$, 
ainsi que les intervalles $I$ et $K$, satisfont les hypoth\`eses du lemme \ref{cc}, 
ce qui permet d'achever la d\'emonstration. \esp $\square$ 


\vspace{-0.6cm}


\beginpicture

\setcoordinatesystem units <1cm,1cm>

\circulararc 360 degrees from 2.5 0 center at 0 0 
\circulararc 24 degrees from 2.43 -0.5 center at 0 0 
\circulararc -24 degrees from -1.7 1.818 center at 0 0 
\circulararc 24 degrees from 2.47 -0.5 center at 0 0 
\circulararc -24 degrees from -1.7 1.858 center at 0 0 
\circulararc 24 degrees from -1.7 -1.818 center at 0 0 
\circulararc 24 degrees from -1.7 -1.857 center at 0 0 

\put{$\bullet$} at 2.45 0.5  
\put{$\bullet$} at 2.45 -0.5 
\put{$\bullet$} at -1.7 1.85
\put{$\bullet$} at -1.7 -1.858
\put{$\bullet$} at -0.78 2.375 
\put{$\bullet$} at -0.78 -2.4

\put{$I=I_1$} at 3.2 0
\put{$I_2$} at -1.5 2.5
\put{$I_{\kappa}$} at -1.5 -2.4 
\put{$f^{-1}$} at 2.1 2.62 
\put{$f$} at 2 -2.5 
\put{$g$} at 1 0.75
\put{$g$} at 1 -0.75
\put{$K$} at 1.25 2.435

\plot 1.45 1.9 
1.6 2.1 /

\plot 0.65 2.3 
0.75 2.53 /

\plot 1.45 -1.9 
1.6 -2.1 /

\plot 0.65 -2.3 
0.75 -2.53 /

\plot -2.32 0.48 
-2.6 0.53 /

\plot -2.32 -0.48 
-2.6 -0.53 /

\circulararc -82 degrees from 2.8 -1 center at 0 0  
\circulararc 82 degrees from 2.8 1 center at 0 0 
\circulararc -75 degrees from 2 0.1 center at 2.5 1.5
\circulararc 80 degrees from 2 -0.1 center at 2.5 -1.5

\plot 2 0.1 
1.87 0.22 /

\plot 2 0.1 
1.83 0.1 /

\plot 1.02 -1.76 
0.93 -1.6 /

\plot 1.02 -1.76 
1.1 -1.6 /

\plot 2.8 -1 
2.82 -1.2 /

\plot 2.8 -1 
2.63 -1.18 /

\plot 2.8 1 
2.82 1.2 /

\plot 2.8 1 
2.63 1.18 /

\put{Figure 6} at 0 -3.5 
\put{} at -8.1 4.3

\endpicture
 

\vspace{0.65cm}

Pour finir ce paragraphe, remarquons que la preuve pr\'ec\'edente montre en fait 
que dans $\Gamma$ il existe des \'el\'ements qui poss\`edent (au moins) $\kappa (\Gamma)$ 
points fixes hyperboliques contractants (resp. dilatants), o\`u $\kappa (\Gamma)$ d\'esigne 
le {\em degr\'e de} $\Gamma$ (voir le \S \ref{ap1}). Dans le paragraphe suivant, et 
en utilisant des m\'ethodes probabilistes, nous montrerons que, m\^eme en classe 
$C^{1}$, il existe des \'el\'ements avec exactement $2 \kappa (\Gamma)$ points 
fixes hyperboliques (toujours sous l'hypoth\`ese de non existence 
de mesure de probabilit\'e invariante).


\subsection{Groupes de diff\'eomorphismes du cercle : m\'ethode probabiliste}
\label{morse-smale}

Les techniques des paragraphes \ref{sac-ps} et \ref{sac-circ} (resp. des 
\S \ref{sac-ps2} et \S \ref{sac-circ}) permettent d\'ej\`a de d\'emontrer 
que si $\Gamma$ est un sous-groupe de $\mathrm{Diff}_+^{1+\tau}(\clo)$ (resp. 
de $\mathrm{Diff}_+^1(\clo)$) qui ne pr\'eserve aucune mesure de probabilit\'e 
du cercle, alors il contient des \'el\'ements avec (au moins) $2 \kappa(\Gamma)$ 
points fixes hyperboliques. Nous en donnons un sch\'ema de preuve, car ceci 
permettra d'illustrer la difficult\'e de la d\'emonstration du th\'eor\`eme D.

Au cours de la preuve de la proposition \ref{sac-det}, on s'est ramen\'e 
\`a \'etudier une situation qui (\`a indice fini pr\`es) est illustr\'ee par 
la figure 4. Or, si l'on regarde l'action de $f$ et $g$ sur l'intervalle $[a',a]$, 
on s'aper\c{c}oit que leur dynamique est analogue \`a celle qui est illustr\'ee 
par la figure 3. Les m\'ethodes du \S \ref{sac-ps} (resp. du \S \ref{sac-ps2}) 
permettent ainsi de trouver des 
\'el\'ements avec des points fixes hyperboliques contractants. Pour 
trouver des points fixes hyperboliques dilatants, l'id\'ee consiste 
\`a examiner l'action de $f^{-1}$ et $g^{-1}$ sur l'intervalle $[b',b]$, 
laquelle aussi est analogue \`a celle illustr\'ee par la figure 3. Nous 
sommes donc confront\'es \`a la difficult\'e d'\'etudier simultan\'ement 
deux processus al\'eatoires, dont l'un va dans le \begin{tiny}$^{_\ll}$\end{tiny}sens 
inverse\hspace{0.05cm}\begin{tiny}$^{_\gg}$\end{tiny} de l'autre. Celui-ci 
n'est qu'un probl\`eme de nature probabiliste et relativement simple du 
point de vue technique.

\vspace{0.2cm}

Le probl\`eme qui consiste \`a trouver des \'el\'ements avec exactement $2 \kappa(\Gamma)$ 
points fixes hyperboliques est bien plus d\'elicat. La strat\'egie que nous proposons 
pour le r\'esoudre diff\`ere radicalement des m\'ethodes employ\'ees dans les 
paragraphes pr\'ec\'edents. Pour sa mise en \oe uvre nous aurons fortement besoin 
des r\'esultats contenus dans l'appendice de ce travail, dont la lecture 
pr\'eliminaire est fondamentale pour la compr\'ehension de ce qui suit.

\vspace{0.1cm}

Si $\Gamma$ est un sous-groupe d\'enombrable de $\mathrm{Hom\acute{e}o}_+(\clo)$ qui ne 
pr\'eserve aucune mesure de probabilit\'e du cercle, alors un argument simple montre 
l'existence d'un sous-groupe de type fini de $\Gamma$ qui ne fixe aucune probabilit\'e 
sur $\clo$. Donc, dans l'\'enonc\'e du th\'eor\`eme D, nous pouvons supposer (sans 
perdre en g\'en\'eralit\'e) que $\Gamma$ est de type fini. Consid\'erons une mesure 
de probabilit\'e $\mu$ sur $\Gamma$ qui soit non d\'eg\'en\'er\'ee, sym\'etrique et 
de support fini, et d\'esignons par $\Omega$ l'espace $\Gamma^{\mathbb{N}}$ 
(muni de la mesure $\mathbb{P} = \mu^{\mathbb{N}}$). Notons 
$\sigma: \Omega \rightarrow \Omega$ le d\'ecalage, et consid\'erons 
la transformation $T: \Omega \times \clo \rightarrow \Omega \times \clo$ 
d\'efinie par $T(\omega,x) = (\sigma(\omega),h_1(\omega)(x))$. Soit $\nu$ 
une probabilit\'e sur $\clo$ {\em stationnaire} par rapport \`a $\mu$, {\em i.e.} telle 
que la transformation $T$ pr\'eserve la mesure de probabilit\'e $\mathbb{P}\!\times\!\nu$. 
Le th\'eor\`eme ergodique de Birkhoff donne l'existence (presque partout) de la limite
$$\lambda_{(\omega,x)}(\nu) = \lim_{n \rightarrow \infty} \frac{\log \big( h_n(\omega)'(x) \big)}{n}.$$
Ce nombre est par d\'efinition l'exposant de Lyapunov du point $(\omega,x) \in \Omega \times \clo$. 

Lorsque $\Gamma$ ne pr\'eserve pas de mesure de probabilit\'e du cercle, la mesure 
stationnaire (par rapport \`a $\mu$) est unique et ne poss\`ede pas d'atome ; par 
ailleurs, la transformation $T$ est ergodique. L'exposant de Lyapunov est donc 
$\mathbb{P}\!\times\!\nu$-presque partout constant ; de plus, la valeur 
$\lambda = \lambda(\nu)$ de cette constante est strictement n\'egative 
(voir les appendices \ref{ap1}, \ref{ap2} et \ref{ap3}).

\vspace{0.2cm}

Pour d\'emontrer le th\'eor\`eme D nous examinons d'abord le cas le plus 
simple, \`a savoir celui d'une action minimale qui v\'erifie la propri\'et\'e 
d'expansivit\'e forte. Dans ce cas, nous savons d'apr\`es le \S \ref{ap1} que 
pour tout $\omega$ appartenant \`a un sous-ensemble $\Omega^*$ de probabilit\'e 
totale de $\Omega$, le {\em coefficient de contraction} $c(h_n(\omega))$ converge vers 
z\'ero lorsque $n$ tend vers l'infini. Cela se traduit par le fait que pour tout 
$\omega \in \Omega^*$ il existe des intervalles ferm\'es $I_n(\omega)$ et 
$J_n(\omega)$ dont la taille converge vers z\'ero et tels que 
$h_n(\omega) \big( \overline{\clo \setminus I_n(\omega)} \big) = J_n(\omega)$. 
De plus, l'intervalle $I_n(\omega)$ converge vers un point $\varsigma_+(\omega)$ 
(et l'application $\varsigma_+: \Omega^* \rightarrow \clo$ ainsi d\'efinie est mesurable). 
Si $I_n(\omega)$ et $J_n(\omega)$ sont disjoints alors les points fixes de $h_n(\omega)$ 
sont inclus \`a l'int\'erieur de 
ces deux intervalles. Pour d\'emontrer l'unicit\'e et l'hyperbolicit\'e 
d'au moins le point fixe contractant ({\em i.e.} celui situ\'e dans $J_n(\omega)$), nous 
aimerions nous servir du fait que l'exposant de Lyapunov de l'action est 
strictement n\'egatif (ce 
qui correspond \`a une propri\'et\'e de contraction diff\'erentielle au niveau local). 
N\'eanmoins, on voit appara\^{\i}tre imm\'ediatement deux difficult\'es 
techniques : les intervalles $I_n(\omega)$ et $J_n(\omega)$ ne sont pas 
toujours disjoints, et la rapidit\'e de contraction diff\'erentielle locale 
d\'epend du point initial $(\omega,x)$. Pour r\'esoudre le premier probl\`eme 
il suffit de remarquer que les instants $n \in \mathbb{N}$ pour lesquels 
$I_n(\omega) \cap J_n(\omega) \neq \emptyset$ sont assez rares (en effet, 
la densit\'e de cet ensemble d'entiers est g\'en\'eriquement \'egale \`a 
z\'ero). La deuxi\`eme difficult\'ee est surmont\'ee en utilisant le fait 
que la transformation $T$ est ergodique, ce qui implique que presque tout 
point initial $(\omega,x)$ tombera au bout d'un certain moment sur un point 
o\`u la rapidit\'e de contraction est bien contr\^ol\'ee. Soulignons finalement 
que, pour obtenir l'unicit\'e et l'hyperbolicit\'e du point fixe dilatant, il 
est n\'ecessaire d'utiliser un argument analogue \`a celui qui pr\'ec\`ede 
mais pour les compositions des inverses et dans l'ordre oppos\'e. Pour cela 
nous passons aux distributions en temps finis ; en effet, puisque la 
mesure de d\'epart $\mu$ est sym\'etrique, en temps fini les 
distributions pour ces deux processus co\"{\i}ncident.

\vspace{0.2cm}

Pour mettre en \oe uvre les id\'ees ci-dessus, commen\c{c}ons par fixer une 
constante $C \geq 1$ suffisamment grande de fa\c{c}on \`a ce que l'ensemble 
$$E (C) = \Big\{ (\omega,x) \in \Omega \times \clo: \esp h_n(\omega)'(x) \leq C 
\exp \Big( \frac{2 \lambda n}{3} \Big) \esp \mbox{ pour tout } \esp n \geq 0 \Big\}$$
ait une mesure strictement positive \footnote{\`A l'aide du lemme classique de Pliss \cite{pliss}, 
on peut d\'emontrer que la mesure de l'ensemble $E(1)$ est d\'ej\`a strictement positive. 
Cependant, nous n'avons pas besoin de cela pour poursuivre notre argumentation.}. Consid\'erons 
$\bar{\varepsilon} > 0$ tel que pour tout $g$ dans le support de $\mu$ et tout 
\esp $x,y$ \esp \`a distance au plus $\bar{\varepsilon}$ on ait
$$\frac{g'(x)}{g'(y)} \leq \exp \Big( \frac{\lambda}{6} \Big).$$
La preuve du lemme suivant est analogue \`a celle du lemme 
\ref{bbb}, et nous la laissons au lecteur.

\vspace{0.1cm}

\begin{lem} {\em Si $(\omega,y)$ appartient \`a $E(C)$ et $x \in \clo$ est un 
point \`a distance inf\'erieure ou \'egale \`a \esp $\varepsilon = 
\bar{\varepsilon}/C$ \esp de $y$, alors}
$$h_n(\omega)'(x) \leq C \esp \exp \Big( \frac{\lambda \esp n}{2} \Big) \quad 
\mbox{ pour tout} \quad n \geq 0.$$
\label{encore-continuity}
\end{lem}

\vspace{-0.3cm}

Pour $\delta \! \in ]0,1[$ posons 
$$\alpha_{\nu}(\delta) = \inf \big\{ \nu(I): I \esp \mbox{ intervalle 
de taille sup\'erieure ou \'egale \`a }  \esp 1\!-\!\delta \big\}.$$
Puisque $\nu$ est sans atome, $\alpha_{\nu}(\delta)$ tend vers $1$ lorsque $\delta$ 
tend vers z\'ero. Pour tout $\omega \in \Omega$ et chaque entier $k > 10$ 
notons\footnote{La condition $k > 10$ est impos\'ee afin que les 
arguments qui suivent ne tombent pas en d\'efaut par des raisons 
\begin{tiny}$^{_\ll}$\end{tiny}stupides\hspace{0.05cm}\begin{tiny}$^{_\gg}$\end{tiny} 
(impossibilit\'e d'existence d'intervalles avec les propri\'et\'es d\'emand\'ees, etc).} 
$$\mathbb{N}(\omega,k) = 
\big\{ n \in \mathbb{N}: \esp 
dist \big( \varsigma_+(\omega),J_n(\omega) \big) \geq 4/k \big\},$$
et posons 
$$\Omega^*_{1/k} = \big\{ \omega \in \Omega^*: \esp 
dens \big( \mathbb{N} (\omega,k) \big) \geq \alpha_{\nu}(10/k) \big\},$$
o\`u $dens$ d\'esigne la {\em densit\'e} de l'ensemble d'entiers strictement positifs 
correspondant\footnote{Nous entendons par densit\'e d'un sous-ensemble $X$ de 
$\mathbb{N}$ la valeur 
$$dens(X) = \liminf_{N \rightarrow \infty} \frac{card ( X \cap \{1,\ldots,N\} )}{N}.$$}.

\vspace{0.1cm}

\begin{lem} {\em Pour chaque entier $k > 10$ la probabilit\'e 
de l'ensemble $\Omega_{1/k}^*$ est totale.}
\end{lem}

\noindent{\bf D\'emonstration.} Lorsque l'action est minimale, le support de 
toute probabilit\'e stationnaire est total. Par cons\'equent, d'apr\`es le 
th\'eor\`eme ergodique de Birkhoff, pour chaque $k > 10$ on peut fixer un point 
$y_k$ dans $\clo \setminus [\varsigma_+(\omega)-5/k,\varsigma_+(\omega)+5/k]$ de sorte 
que pour un sous-ensemble g\'en\'erique $\bar{\Omega}_{1/k}^*$ de $\Omega^*$ on ait
\begin{equation}
\lim_{N \rightarrow \infty} \frac{card \big\{n \leq N: \esp h_n(p_k) \notin 
 [\varsigma_+(\omega)-5/k,\varsigma_+(\omega)+5/k] \big\} }{N} = 
\nu \big( \clo \setminus [\varsigma_+(\omega)-5/k,\varsigma_+(\omega)+5/k] \big).
\label{birk}
\end{equation}
Pour chaque $\omega \in \bar{\Omega}^*_{1/k}$ 
fixons $n(\omega) \in \mathbb{N}$ de sorte que 
\esp $c \big( h_n(\omega) \big) \leq 1/k$ \esp et \esp 
$I_{n}(\omega) \subset [\varsigma_+(\omega)-1/k,\varsigma_+(\omega)+1/k]$ \esp 
pour tout $n \geq n(\omega)$. Si $n \geq n(\omega)$ alors $y_k$ est \'evidemment  
dans $\clo \setminus I_n(\omega)$, et donc $h_n(\omega)(y_k)$ appartient \`a 
$J_{n}(\omega)$. Par suite, si $n \geq n(\omega)$ est tel que $h_n(\omega)(y_k)$ 
n'est pas dans $[\varsigma_+(\omega)-5/k,\varsigma_+(\omega)+5/k]$, alors 
$dist \big( \varsigma_+(\omega),J_n(\omega) \big) \geq 4/k$. D'apr\`es (\ref{birk}), 
ceci implique que l'ensemble $\bar{\Omega}^*_{1/k}$ est contenu dans 
$\Omega_{1/k}^*$, ce qui permet de conclure la preuve du lemme. \esp $\square$

\vspace{0.35cm}

Puisque l'application $T$ est ergodique, si l'on d\'esigne par $\Omega(C)$ l'ensemble 
des $\omega \!\in\! \Omega$ tels que pour $\nu$-presque tout \esp $y\in\clo$ \esp on a 
\esp $\big( \sigma^n(\omega),h_n(\omega)(y) \big) \!\in \!E(C)$ pour une infinit\'e 
d'entiers $n \in \mathbb{N}$, alors la probabilit\'e de $\Omega(C)$ est totale. Par 
suite, pour tout entier $k > 10$ l'ensemble $\Omega_{1/k}^{**} = \Omega_{1/k}^* \cap 
\Omega(C)$ poss\`ede lui aussi une probabilit\'e \'egale \`a $1$. Pour chaque 
$\eta > 0$ d\'esignons par $D_{\eta}(\clo)$ l'ensemble des diff\'eomorphismes $g$ 
de classe $C^1$ du cercle qui v\'erifient la propri\'et\'e suivante : il existe 
deux intervalles ferm\'es $I'$ et $J'$ de taille au plus $\eta$ et \`a distance 
au moins $2 \eta$ entre eux tels que $g$ envoie $\clo \setminus I'$ dans $J'$ et la 
contraction $g: \clo \setminus I' \rightarrow J'$ est uniform\'ement hyperbolique 
({\em i.e.} la d\'eriv\'ee de $g$ sur $\overline{\clo \setminus I'}$ est 
strictement inf\'erieure \`a $1$).

\vspace{0.1cm}

\begin{lem} {\em Si $\omega$ appartient \`a $\Omega_{1/k}^{**}$ alors l'ensemble 
$\mathbb{N}^*(\omega,k)$ des entiers $n \in \mathbb{N}$ tels que $h_n(\omega)$ 
appartient \`a $D_{1/k}(\clo)$ a une densit\'e au moins \'egale \`a \esp $\alpha_{\nu}(10/k)$.}
\label{entrar}
\end{lem}

\noindent{\bf D\'emonstration.} Posons $\varepsilon_k = \min \{1/2k,\varepsilon\}$, 
et pour chaque $\omega \in \Omega_{1/k}^{**}$ prenons :

\vspace{0.15cm}

\noindent{-- \esp \thinspace un entier positif $n_0 = n_0(\omega)$ 
tel que \esp $c \big(h_n(\omega) \big) \leq \varepsilon_k$ \esp et \esp 
$I_n(\omega) \subset [\varsigma_+(\omega)-1/k,\varsigma_+(\omega)+1/k]$ 
pour tout $n \geq n_0$~;}

\vspace{0.15cm}

\noindent{-- \esp un point $y$ n'appartenant pas 
\`a $[\varsigma_+(\omega) - 1/k,\varsigma_+(\omega) + 1/k]$ 
et un entier positif $n_1 = n_1(\omega) \geq n_0$ tels que 
$\big( \sigma^{n_1}(\omega),h_{n_1}(\omega)(y) \big) \in E(C)$ ;}

\vspace{0.15cm}

\noindent{-- \esp un entier positif $n_2 = n_2(\omega)$ tel que si l'on d\'esigne 
par $M$ le supremum des valeurs de $g'(z)$ avec $z \in \clo$ et $g$ dans le 
support de $\mu$, alors}
$$C \exp \Big( \frac{\lambda \esp n_2}{2} \Big) \esp M^{n_1} < 1.$$

\vspace{0.1cm}

Supposons que $n \in \mathbb{N}$ soit sup\'erieur ou \'egal 
\`a $n_1+n_2$ et qu'il appartienne \`a $\mathbb{N}(\omega,k)$. 
Le point $h_{n_1}(\omega)(y)$ appartient \`a l'intervalle 
$J_{n_1}(\omega)$, dont la taille n'est pas plus grande que $\varepsilon$. 
Par le lemme \ref{encore-continuity}, pour tout $x \in J_{n_1}(\omega)$ 
et tout $m \geq 0$,
$$h_m \big( \sigma^{n_1}(\omega) \big)'(x) \leq C \exp 
\Big(\frac{\lambda \esp m}{2} \Big).$$
En particulier, si $m \geq n_2$ alors pour tout $z \notin I_{n_1}(\omega)$ on a
$$h_{n_1+m}(\omega)'(z) = h_{n_1}(\omega)'(z) \cdot 
h_m \big( \sigma^{n_1}(\omega) \big)' \big( h_{n_1}(\omega)(z) \big) 
\leq M^{n_1} \esp C \esp \exp \Big( \frac{\lambda \esp m}{2} \Big) < 1.$$
Donc, si $m \geq n_2$ est tel que $n = n_1+m$ appartient \`a 
$\mathbb{N}(\omega,k)$, alors :

\vspace{0.15cm}

\noindent{-- \esp la d\'eriv\'ee de \esp $h_n(\omega)$ \esp sur \esp 
$\overline{\clo \setminus I_{n_1}(\omega)}$ \esp est born\'ee par une 
constante strictement inf\'erieure \`a $1$ ;} 

\vspace{0.15cm}

\noindent{-- \esp la taille de l'intervalle \esp $h_n(\omega)(I_{n_1}(\omega)) 
= h_{m}(\sigma^{n_1}(\omega))(J_{n_1}(\omega))$ \esp est inf\'erieure 
ou \'egale \`a \esp $|J_{n_1}(\omega)| \leq \varepsilon_k$ ;}

\vspace{0.15cm}

\noindent{-- \esp on a \esp 
$h_n(\omega) \big( \overline{\clo \setminus I_{n}(\omega)} \big) = J_n(\omega)$.}

\vspace{0.1cm}

Puisque les quatre intervalles concern\'es ci-dessus ont une taille au 
plus \'egale \`a \esp $\varepsilon_k \leq 1/2k$, \esp des arguments \'el\'ementaires 
montrent que \esp $I_{n_1}(\omega)$ \esp et \esp $I_{n}(\omega)$ \esp s'intersectent, et de 
m\^eme pour \esp $J_n(\omega)$ \esp et \esp $h_n(\omega)(I_{n_1}(\omega))$ ; de 
plus si l'on d\'esigne par \esp $I_n'(\omega)$ \esp (resp. \esp $J_n'(\omega)$) 
\esp l'intervalle \esp $I_{n_1}(\omega) \cup I_{n}(\omega)$ \esp 
(resp. \esp $J_n(\omega) \cup h_n(\omega)(I_{n_1}(\omega))$), \esp 
alors \esp $I_n'(\omega)$ \esp et \esp $J_n'(\omega)$ \esp ont une taille 
au plus \'egale \`a $1/k$, ils sont \`a une distance 
sup\'erieure ou \'egale \`a $2/k$ et \esp 
$h_n(\omega) \big( \clo \setminus I_{n}'(\omega) \big) \subset J_n'(\omega)$. 

Nous avons donc montr\'e que tout entier \esp $n \geq n_1(\omega) + n_2(\omega)$ 
\esp dans $\mathbb{N}(\omega,k)$ appartient \`a $\mathbb{N}^*(\omega,k)$. Puisque 
la densit\'e de $\mathbb{N}(\omega,k)$ est sup\'erieure ou \'egale \`a \esp 
$\alpha_{\nu}(10/k)$, \esp ceci conclut la preuve du lemme. \esp $\square$

\vspace{0.3cm}

Fixons maintenant une constante $\rho > 1/2$ et un entier $k_0 > 10$ tels que 
$\alpha_{\nu} (10/k_0) > \rho$. Prenons une autre constante $\varrho \!\in ]0,1[$ 
telle que $\rho \varrho > 1/2$, et pour chaque $N \in \mathbb{N}$ d\'esignons 
par $\Omega(k_0,N)$ l'ensemble des $\omega \in \Omega_{1/k_0}^{**}$ tels que 
$$\frac{card \big\{ n \in \{1,\ldots,N\}: h_n(\omega) \in D_{1/k_0}(\clo) \big\}}{N} 
\geq \rho.$$
Bien \'evidemment, pour $N$ suffisamment grand nous avons 
\esp $\mathbb{P} \big[ \Omega(k_0,N) \big] \geq \varrho$. \esp Pour  
$n \in \mathbb{N}$ notons $\mathbb{P}_n$ la probabilit\'e $\mu^n$ 
d\'efinie sur $\Gamma^n$.

\vspace{0.15cm}

\begin{lem} {\em Si $\mathbb{P} \big[ \Omega (k_0,N) \big] \geq \varrho$ alors il 
existe $n \in \{1,\ldots,N\}$ tel que}
\begin{equation}
\mathbb{P}_n \big[ (g_1,\ldots,g_n): g_n \cdots g_1 \in D_{1/k_0}(\clo) \big] 
\geq \rho \varrho.
\label{ro}
\end{equation}    
\label{densidades}
\end{lem}
\vspace{-0.5cm}
\noindent{\bf D\'emonstration.} Pour chaque $n \in \mathbb{N}$ notons $\xi_n$ la 
variable al\'eatoire d\'efinie par
$$\xi_n(\omega) = \mathcal{X}_{D_{1/k_0}(\clo)} \big( h_n(\omega) \big),$$
o\`u $\mathcal{X}$ d\'esigne la fonction caract\'eristique. Un \'el\'ement 
$\omega \in \Omega_{1/k_0}^{**}$ appartient \`a $\Omega(k_0,N)$ si et seulement si
$$\frac{1}{N} \sum_{n=1}^N \xi_n(\omega) \geq \rho.$$
Donc, si $\mathbb{P} \big[ \Omega(k_0,N) \big] \geq \varrho$ alors 
$$\mathbb{E} \Big( \frac{1}{N} \sum_{n=1}^N \xi_n \Big) \geq \rho \varrho.$$
Par suite, il existe $n \in \{1,\ldots,N\}$ tel que \esp 
$\mathbb{E}(\xi_n) \geq \rho \varrho$. \esp Puisque 
$$\mathbb{E}(\xi_n) = 
\mathbb{P}_n \big[ (g_1,\ldots,g_n): g_n \cdots g_1 \in D_{1/k_0}(\clo) \big],$$
ceci ach\`eve la preuve du lemme. \esp $\square$

\vspace{0.25cm}

Rappelons que la mesure $\mu$ sur $\Gamma$ est suppos\'ee sym\'etrique. Par suite, la distribution 
dans $\Gamma$ des \'el\'ements de la forme $g_n \cdots g_1$ obtenus \`a partir d'\'el\'ements 
$(g_1,\ldots,g_n) \in \Gamma^n$ est exactement la m\^eme que celle des $g_1^{-1} \cdots g_n^{-1}$. 
Donc, si $n \in \mathbb{N}$ v\'erifie (\ref{ro}) alors
$$\mathbb{P}_n \big[ (g_1,\ldots,g_n): \esp g_n \cdots g_1 \esp \mbox{ et } \esp 
(g_n \cdots g_1)^{-1 } \esp \mbox{ appartiennent \`a } \esp D_{1/k_0}(\clo) \big] 
\geq 2 \rho \varrho - 1 > 0,$$
La preuve du th\'eor\`eme D (dans le cas minimal et fortement expansif) est 
alors conclue par le lemme suivant.

\vspace{0.12cm}

\begin{lem} {\em Soit $g$ un diff\'eomorphisme de classe $C^1$ du cercle. Si $g$ 
et $g^{-1}$ appartiennent \`a $D_{1/k_0}(\clo)$ alors $g$ ne poss\`ede que deux points 
fixes, l'un hyperboliquement contractant et l'autre hyperboliquement dilatant.}
\label{dos-sentidos}
\end{lem}

\noindent{\bf D\'emonstration.} Il est facile de voir que l'intersection des 
intervalles $I'(g)$ et $J'(g^{-1})$ est d'int\'erieur non vide, et de m\^eme 
pour $J'(g)$ et $I'(g^{-1})$. Puisque tous ces intervalles ont une taille au 
plus \'egale \`a $1/k_0$ et \esp $dist \big( I'(g),J'(g) \big) \geq 2/k_0$, \esp 
les intervalles $K(g) = I'(g) \cup J'(g^{-1})$ et $K'(g)= J'(g) \cup I'(g^{-1})$ 
sont disjoints. De plus, $g$ envoie $\clo \setminus K(g)$ dans $K'(g)$, et 
$g^{-1}$ envoie $\clo \setminus K'(g)$ dans $K(g)$. Puisque la d\'eriv\'ee de $g$ 
(resp. de $g^{-1}$) dans $\clo \setminus K(g)$ (resp. dans $\clo \setminus K'(g)$) 
est strictement inf\'erieure \`a $1$, ceci d\'emontre le lemme. \esp $\square$ 

\vspace{0.07cm}

\begin{rem} Une bonne lecture de la preuve pr\'ec\'edente permet de conclure que, 
lorsqu'on fait des compositions al\'eatoires, tr\`es probablement \`a partir d'un certain 
moment on tombera assez souvent sur un diff\'eomorphisme avec deux points fixes hyperboliques. 
\end{rem}

\vspace{0.07cm}

L'extension de la d\'emonstration pr\'ec\'edente au cas minimal et expansif (mais non 
fortement expansif) ne pr\'esente pas de difficult\'e majeure. En effet, nous savons 
d'apr\`es le \S \ref{ap1} qu'il existe un hom\'eomorphisme $\theta: \clo \rightarrow \clo$ 
d'ordre fini $\kappa \!=\! \kappa(\Gamma)$ et qui commute avec tous les \'el\'ements 
de $\Gamma$. De plus, si l'on d\'esigne par $\clo\!/\!\sim$ le cercle topologique obtenu 
en identifiant les points des orbites de $\theta$ et par $\hat{\Gamma}$ le sous-groupe 
de $\mathrm{Hom\acute{e}o}_+(\clo\!/\!\sim)$ induit, alors l'action de $\hat{\Gamma}$ est 
minimale et fortement expansive. Remarquons que $\theta$ n'est pas n\'ecessairement 
diff\'erentiable. Cependant, tous les arguments topologiques de la preuve donn\'ee pour 
le cas minimal et fortement expansif s'appliquent \`a l'action de $\hat{\Gamma}$ sur 
 $\mathrm{Hom\acute{e}o}_+(\clo\!/\!\sim)$. D'autre part, les arguments de nature 
diff\'erentielle doivent \^etre l\^us dans l'action originale de $\Gamma$ sur $\clo$. 
On d\'emontre ainsi l'existence d'un \'el\'ement 
$g = g_n \cdots g_1$ dans $\Gamma$ tel qu'il existe des intervalles 
disjoints $K(g)$ et $K'(g)$ dans $\clo\!/\!\sim$ qui se r\'el\`event respectivement 
en des intervalles $K_1(g),\ldots,K_{\kappa}(g)$ et $K_1'(g),\ldots,K_{\kappa}'(g)$ 
qui sont deux \`a deux disjoints et qui v\'erifient :

\vspace{0.15cm}

\noindent{-- \esp $g$ envoie $\clo \setminus \cup_{i=1}^{\kappa} K_i(g)$ sur 
$\cup_{i=1}^{\kappa} K_i'(g)$ avec d\'eriv\'ee strictement inf\'erieure \`a $1$ ;}

\vspace{0.15cm}

\noindent{-- \esp $g^{-1}$ envoie $\clo \setminus \cup_{i=1}^{\kappa} K_i'(g)$ sur 
$\cup_{i=1}^{\kappa} K_i(g)$ avec d\'eriv\'ee strictement inf\'erieure \`a $1$.}

\vspace{0.15cm}

Il est alors facile de s'apercevoir que l'\'el\'ement $g^{\kappa} \in \Gamma$ ne contient 
que $2\kappa$ points fixes, dont la moiti\'e sont hyperboliquement contractants et l'autre 
moiti\'e sont hyperboliquement dilatants.

\vspace{0.1cm}

\begin{rem} Le dernier argument de la preuve pr\'ec\'edente ({\em i.e.} celui qui consiste \`a 
remplacer $g$ par $g^{\kappa}$) n'est pas de nature probabiliste, dans le sens que si bien $g$ 
correspond \`a un produit al\'eatoire g\'en\'erique d'\'el\'ements de $\Gamma$, l'\'el\'ement 
$g^{\kappa} \in \Gamma$ ne jouit pas de cette propri\'et\'e de g\'en\'ericit\'e. Il est donc 
naturel d'essayer de d\'eterminer la fr\'equence avec laquelle on voit appara\^{\i}tre des 
diff\'eomorphismes n'ayant que des points fixes hyperboliques dans une suite al\'eatoire 
typique. Ce probl\`eme a une r\'eponse \'evidente -cette fr\'equence est (g\'en\'eriquement) 
\'egale \`a $1 / \kappa$- mais la d\'emonstration de ce fait n'est pas compl\'etement 
\'el\'ementaire (elle utilise une loi de type $0\!-\!2$).
\end{rem}

\vspace{0.1cm}

Le cas o\`u il existe un minimal exceptionnel (et il n'y a pas de probabilit\'e invariante) 
est plus compliqu\'e. En effet, dans cette situation le support de l'unique mesure 
stationnaire co\"{\i}ncide avec cet ensemble minimal $\Lambda$ ; le fait que l'exposant 
de Lyapunov soit strictement n\'egatif ne donne donc aucune information d'hyperbolicit\'e 
loin de l'ensemble $\Lambda$.

En \'ecrasant les composantes connexes de $\clo \setminus \Lambda$ on obtient un cercle 
topologique $\clo_{\Lambda}$ muni d'une action (non n\'ecessairement fid\`ele) de $\Gamma$ 
par hom\'eomorphismes (pour chaque $g \in \Gamma$ nous noterons encore $g$ l'hom\'eomorphisme 
induit sur $\clo_{\Lambda}$). Si l'on d\'esigne par $\Gamma_{\Lambda}$ le sous-groupe 
correspondant de $\mathrm{Homeo}_+(\clo_{\Lambda})$, alors l'action de $\Gamma_{\Lambda}$ 
est minimale et expansive. Pour simplifier, nous supposerons dans la suite que cette action 
est fortement expansive, en laissant le traitement du cas g\'en\'eral au lecteur. L'id\'ee 
de la preuve consiste alors \`a appliquer les arguments topologiques du d\'ebut de ce 
paragraphe \`a l'action de $\Gamma_{\Lambda}$, et puis de lire les propri\'et\'es 
de nature diff\'erentielle dans l'action originale sur $\clo$. Bien \'evidemment, 
pour aboutir \`a l'hyperbolicit\'e nous devrons travailler 
\begin{tiny}$^{_\ll}$\end{tiny}suffisamment 
loin\hspace{0.05cm}\begin{tiny}$^{_\gg}$\end{tiny} 
des composantes connexes de $\clo\setminus \Lambda$ qui sont 
\begin{tiny}$^{_\ll}$\end{tiny}trop 
grandes\hspace{0.05cm}\begin{tiny}$^{_\gg}$\end{tiny}.

Remarquons que le cercle topologique $\clo_{\Lambda}$ peut \^etre muni d'une structure m\'etrique 
naturelle : on peut le param\'etrer en utilisant la mesure $\nu$. Nous entendons alors par 
\begin{tiny}$^{_\ll}$\end{tiny}taille\hspace{0.05cm}\begin{tiny}$^{_\gg}$\end{tiny} 
d'un intervalle $V$ de $\clo_{\Lambda}$ la valeur de $\nu(V)$, et par 
\begin{tiny}$^{_\ll}$\end{tiny}distance\hspace{0.05cm}\begin{tiny}$^{_\gg}$\end{tiny} 
entre deux points $p$ et $q$ de $\clo_{\Lambda}$ la valeur \esp 
$dist_{\Lambda}(p,q) = \min \big\{ \nu \big( [p,q] \big),\nu \big( [q,p] \big) \big\}$. 

Notons encore une fois 
$$E (C) = \Big\{ (\omega,x) \in \Omega \times \clo: \esp h_n(\omega)'(x) \leq C 
\exp \Big( \frac{2 \lambda n}{3} \Big) \esp \mbox{ pour tout } \esp n \geq 0 \Big\}$$
et fixons une constante $C \geq 1$ telle que \esp 
$\mathbb{P}\!\times\!\nu \big[ E(C) \big] > 0$. Consid\'erons 
$\bar{\varepsilon} > 0$ tel que pour tout $g$ dans le support de $\mu$ et tout 
\esp $x,y$ \esp \`a distance au plus $\bar{\varepsilon}$ on ait \esp 
$g'(x)/g'(y) \leq \exp (\lambda / 6)$, \esp et posons 
$\varepsilon = \bar{\varepsilon} / C$. 
Pour chaque $\eta > 0$ les composantes connexes de $\clo \setminus \Lambda$ de taille 
sup\'erieure ou \'egale \`a \esp $\eta / 3$ \esp se projettent en des points distincts 
$p_1,\ldots,p_{m(\eta)}$ de $\clo_{\Lambda}$. Il n'est pas difficile de voir que 
si $\eta \leq \varepsilon$ est suffisamment petit (et strictement positif), alors 
il existe des constantes strictement positives $\rho$, $\bar{\rho}$, $\varrho$, 
$\bar{\varrho}$, $\delta$ et $\bar{\delta}$ telles que :

\vspace{0.15cm}

\noindent{$(i)$ \esp $1 > \bar{\varrho} > \varrho > \bar{\rho} > \rho > 1/2$, \esp \esp 
$\rho \varrho > 1/2$~;}

\vspace{0.15cm}

\noindent{$(ii)$ \esp si $V$ est un intervalle quelconque de $\clo_{\Lambda}$ 
de taille au plus $\delta$ et dont le centre est \`a distance au moins \esp 
$(1-\bar{\varrho})/2m(\eta)$ \esp de tous les $p_i$, alors la taille 
de la pr\'eimage de $V$ dans $\clo$ est inf\'erieure ou \'egale \`a $\eta/2$ ;}

\vspace{0.15cm}

\noindent{$(iii)$ \esp si $V$ est un intervalle quelconque de $\clo_{\Lambda}$ de taille 
au moins \'egale \`a $\bar{\delta}$, alors la taille de la pr\'eimage de $V$ dans $\clo$ 
est sup\'erieure ou \'egale \`a $2 \eta$~;}

\vspace{0.15cm}

\noindent{$(iv)$ \esp on a l'in\'egalit\'e \esp 
$\bar{\varrho} - 2(\delta + \bar{\delta}) \geq \bar{\rho}$.}

\vspace{0.15cm}

Fixons de telles constantes et pour chaque $i \in \{1,\ldots,m(\eta)\}$ d\'esignons par $U_i$ 
(resp. $U_i^+$) l'intervalle ferm\'e dans $\clo_{\Lambda}$ de centre $p_i$ et de taille \esp 
$(1-\bar{\varrho})/m(\eta)$ \esp (resp. de taille \esp $(1-\bar{\varrho} +\delta)/m(\eta)$). 
\esp Bien \'evidement,
\begin{equation}
\nu \big( \clo_{\Lambda} \setminus \cup_{i=1}^{m(\eta)} U_i \big) = \bar{\varrho}, \qquad 
\nu \big( \clo_{\Lambda} \setminus \cup_{i=1}^{m(\eta)} U_i^+ \big) = \bar{\varrho} - \delta.
\label{mesure}
\end{equation}

Pour chaque $\omega \in \Omega^*$ notons $I_n(\omega)$ et $J_n(\omega)$ les intervalles 
de $\clo_{\Lambda}$ donn\'es par la propri\'et\'e d'expansivit\'e forte de l'action de 
$\Gamma_{\Lambda}$ sur $\clo_{\Lambda}$, et d\'esignons par 
$\tilde{I}_n(\omega)$ et $\tilde{J}_n(\omega)$ ses pr\'eimages dans $\clo$. Posons 
$$\mathbb{N}_*(\omega,\eta)  = \big\{ n \in \mathbb{N} : \esp \varsigma_{+}(\omega) \in 
\clo_{\Lambda} \setminus \cup_{i=1}^{m(\eta)} U_i, \esp |\tilde{I}_n(\omega)| \leq \eta/2, 
\esp |\tilde{J}_n(\omega)| \leq \eta/2, \esp 
dist \big( \tilde{I}_n(\omega),\tilde{J}_n(\omega) \big) \geq 2 \eta \big\},$$
$$\Omega_*(\eta) = \{\omega \in \Omega^*: \esp dens \big( \mathbb{N}_*(\omega,\eta) \big) 
\geq \bar{\varrho}\}.$$

\vspace{0.1cm}

\begin{lem} {\em La probabilit\'e de l'ensemble $\Omega_*(\eta)$ est sup\'erieure ou \'egale 
\`a \esp $\bar{\varrho}$.}
\label{sep+}
\end{lem}

\noindent{\bf D\'emonstration.} La distribution des points $\varsigma_+(\omega)$ co\"{\i}ncide 
avec $\nu$ (voir la preuve de la proposition \ref{unicidad} plus loin). Par suite, 
$$\mathbb{P} \big[ \omega \in \Omega^*: \esp 
\varsigma_+(\omega) \in \clo_{\Lambda} \setminus \cup_{i=1}^{m(\eta)} U_i \big] = 
\nu \big( \clo_{\Lambda} \setminus \cup_{i=1}^{m(\eta)} U_i \big) = \bar{\varrho}.$$ 
Nous affirmons que presque tout $\omega$ dans l'ensemble
$$\bar{\Omega}_*(\eta) = \big\{\omega \in \Omega^*: \esp \varsigma_+(\omega) \in 
\clo_{\Lambda} \setminus \cup_{i=1}^{m(\eta)} U_i \big\}$$
appartient \`a $\Omega_*(\eta)$, ce qui implique \'evidemment l'affirmation du 
lemme. Pour v\'erifier cela, notons d'abord que si $\varsigma_+(\omega)$ appartient 
\`a \esp $\clo_{\Lambda} \setminus \cup_{i=1}^{m(\eta)} U_i$ \esp alors $I_n(\omega)$ est 
contenu dans \esp $\clo_{\Lambda} \setminus \cup_{i=1}^{m(\eta)} U_i$ \esp pour tout $n$ 
assez grand. Donc, d'apr\`es $(ii)$, \esp on a $|\tilde{I}_n(\omega)| \leq \eta/2$ 
\esp pour tout $n$ suffisamment grand. Pour chaque $\omega \in \Omega^*$ d\'esignons 
par $V(\omega)$ l'intervalle de centre $\varsigma_+(\omega)$ et de taille \esp 
$2\bar{\delta} + \delta$. \esp D'apr\`es le th\'eor\`eme ergodique de Birkhoff et 
(\ref{mesure}), pour presque tout $\omega \in \bar{\Omega}_*(\eta)$ il existe \esp 
$p \in \clo_{\Lambda} \setminus [\varsigma_+(\omega)-1/3,\varsigma_+(\omega)+1/3]$ tel que 
\begin{equation}
\lim_{N \rightarrow \infty} 
\frac{card \big\{ n \in \{1,\ldots,N\}: \esp h_n(\omega) \notin \cup_{i=1}^{m(\eta)} 
U_i^+ \cup V(\omega) \big\}}{N} 
\geq \bar{\varrho} - 2(\delta + \bar{\delta}).
\label{minor}
\end{equation}
Pour $n$ assez grand on a $h_n(\omega)(p) \in J_n(p)$ et $\nu \big( J_n(\omega) \big) \leq \delta$, 
et donc la condition \esp $h_n(\omega)(p) \notin \cup_{i=1}^{m(\eta)} U_i^+$ \esp (resp. \esp 
$h_n(\omega)(p) \notin V(\omega)$) \esp implique d'apr\`es $(ii)$ que \esp 
$|\tilde{J}_n(\omega)| \leq \eta/2$ \esp (resp. implique d'apr\`es $(iii)$ que la 
valeur de \esp $dist \big( \tilde{I}_n(\omega), \tilde{J}_n(\omega) \big)$ est 
sup\'erieure ou \'egale \`a $2 \eta$). \esp Le lemme s'en suit alors de (\ref{minor}) 
et de la condition \esp $\bar{\varrho} - 2(\delta + \bar{\delta}) \geq \bar{\rho}$. 
\esp $\square$

\vspace{0.35cm}

La d\'emonstration continue de mani\`ere analogue \`a celle donn\'ee 
pour le cas minimal (et fortement expansif). On note $\Omega(C)$ l'ensemble 
des $\omega \in \Omega$ tels que pour $\nu$-presque tout $y \in \clo$ on ait 
$\big( \sigma^n(\omega),h_n(\omega)(y) \big) \!\in\! E(C)$ pour une infinit\'e 
d'entiers $n \in \mathbb{N}$, et on consid\`ere 
l'ensemble \esp $\Omega_{**}(\eta) = \Omega_*(\eta) \cap \Omega(C)$ \esp (dont 
la probabilit\'e est encore sup\'erieure ou \'egale \`a  $\bar{\varrho}$). Si l'on 
tient compte du fait que $\eta \leq \varepsilon$ (et que la conclusion du 
lemme \ref{encore-continuity} est donc valable), alors par des arguments 
analogues \`a ceux du lemme \ref{entrar} on constate ais\'ement que pour tout 
$\omega \in \Omega_{**}(\eta)$ la densit\'e des entiers $n \in \mathbb{N}$ tels 
que $h_n(\omega)$ appartient \`a $D_{\eta}(\clo)$ est sup\'erieure ou \'egale \`a 
$\bar{\rho}$. Par suite, pour $N$ suffisamment grand on a 
$$\mathbb{P} \Big[ \omega \in \Omega_{**}(\eta): \esp \frac{card \big\{n \in \{1,\ldots,N\}: 
\esp h_n(\omega) \in D_{\eta}(\clo) \big\}}{N}  \geq \rho \Big] \geq \varrho.$$
La preuve du lemme \ref{densidades} montre alors qu'il existe 
$n \in \{1,\ldots,N\}$ tel que
$$\mathbb{P}_n \big[ (g_1,\ldots,g_n): \esp g_n \cdots g_1 \in D_{\eta}(\clo) \big] 
\geq \rho \varrho > 1/2.$$
En passant aux inverses on conclut que
$$\mathbb{P}_n \big[ (g_1,\ldots,g_n) : \esp \esp g_n \cdots g_1 \esp \mbox{ et } \esp 
(g_n \cdots g_1)^{-1 } \esp \mbox{ appartiennent \`a } \esp D_{\eta}(\clo) \big] 
\geq 2 \rho \varrho - 1 > 0,$$
ce qui permet de conclure la d\'emonstration par un argument analogue \`a celui 
du lemme \ref{dos-sentidos}.


\section{Appendice}


\subsection{Une alternative topologique pour des pseudo-groupes d'hom\'eomorphismes}
\label{ap0}

Soient $X$ une vari\'et\'e unidimensionnelle compacte, $\Gamma$ un pseudo-groupe d'hom\'eomorphismes 
de $X$ et $\Lambda$ un sous-ensemble de $X$ qui est compact, invariant par $\Gamma$, non vide 
et minimal. L'orbite $O(x)$ d'un point $x \in X$ est dite {\em de type ressort} si d'une part 
il existe un \'el\'ement $g \in \Gamma$ dont le domaine de d\'efinition $dom(g)$ contient 
un demi-intervalle $[x,p[$ ou $]p,x]$ pour lequel $x$ est un point fixe attractif 
topologique ({\em i.e.} $g^n(p)$ converge vers $x$ quand $n$ tend vers l'infini), et d'autre 
part $O(x)$ intersecte ce demi-intervalle. La proposition suivante s'appuie fortement sur la 
partie moins connue (et pourtant tr\`es int\'eressante) de l'article de Sacksteder~\cite{Sa}.

\vspace{0.1cm}

\begin{prop} {\em Avec les notations pr\'ec\'edentes, soit il existe une mesure de probabilit\'e 
support\'ee sur $\Lambda$ et invariante par tous les \'el\'ements de $\Gamma$, soit il existe 
dans $\Lambda$ une orbite de type ressort.}
\label{sacc}
\end{prop}

\vspace{0.1cm}

Pour d\'emontrer cette proposition, commen\c{c}ons par remarquer qu'il y a trois possibilit\'es 
pour l'ensemble minimal $\Lambda$~: soit il est fini, soit c'est un ensemble de Cantor, soit 
$\Lambda=X$. Ceci se montre en remarquant que la fronti\`ere $Fr(\Lambda)$ et l'ensemble 
$\Lambda^{ac}$ des points d'accumulation de $\Lambda$ sont des ferm\'es invariants contenus 
dans $\Lambda$. Donc soit ils sont vides, soit il sont \'egaux \`a $\Lambda$. En \'etudiant 
toutes les possibilit\'es on d\'emontre le fait escompt\'e.

Dans le cas o\`u $\Lambda$ est fini, la moyenne des masses de Dirac aux points de $\Lambda$ est 
une mesure invariante par $\Gamma$. Dans le cas o\`u $\Lambda = X$, la proposition \ref{sacc} 
a \'et\'e d\'emontr\'ee par Sacksteder (voir le lemme 9.1 de \cite{Sa}). Notre d\'emonstration 
pour le cas o\`u $\Lambda$ est exceptionnel ({\em i.e.} c'est un ensemble de Cantor) suit sa 
strat\'egie~: nous d\'efinissons une notion d'\'equicontinuit\'e pour la restriction de 
$\Gamma$ \`a $\Lambda$ et puis nous prouvons que si cette restriction n'est pas 
\'equicontinue, alors il existe une orbite de type ressort dans $\Lambda$. De 
plus, dans le cas \'equicontinu nous d\'emontrons l'existence d'une mesure de 
probabilit\'e invariante par $\Gamma$ et support\'ee sur $\Lambda$. La 
subtilit\'e de la preuve vient de la d\'efinition suivante.

\begin{defn} Fixons une mesure de probabilit\'e $\nu$ sur $\Lambda$ qui soit sans atome et 
dont le support co\"{\i}ncide avec tout l'ensemble $\Lambda$. Nous dirons que la restriction de 
$\Gamma$ \`a $\Lambda$ est \'equicontinue si pour tout $\varepsilon >0$ il existe $\delta >0$ 
tel que pour tout intervalle $I$ de $X$ de mesure $\nu(I)$ inf\'erieure \`a $\delta$ et 
pour tout \'el\'ement $g$ de $\Gamma$ dont le domaine de d\'efinition contient $I$, 
on a $\nu \big( g(I) \big) < \varepsilon$. 
\end{defn}

Il est facile de v\'erifier que cette notion d'\'equicontinuit\'e ne d\'epend pas 
du choix de la mesure $\nu$. 

\vspace{0.1cm}

\begin{lem} Si la restriction de $\Gamma$ \`a $\Lambda$ est non \'equicontinue 
alors $\Lambda$ contient une orbite de type ressort. 
\end{lem} 

\noindent {\bf D\'emonstration.} Supposons que la restriction de $\Gamma$ \`a $\Lambda$ ne soit pas 
\'equicontinue. Il existe donc $\varepsilon_0 >0$, une suite d'intervalle $I_n$ dont les mesures 
tendent vers $0$, et une suite d'\'el\'ements $g_n \!\in\! \Gamma$ dont les domaines de d\'efinition 
contiennent $I_n$, tels que $\nu \big( g_n(I_n) \big) \geq \varepsilon_0$. Supposons que les longueurs 
des intervalles 
$I_n$ soient minor\'ees par une constante strictement positive ind\'ependante de $n$. Dans ce cas, nous 
pouvons extraire de la suite des $I_n$ une sous-suite d'intervalles qui converge vers un intervalle $I$. 
Puisque $\nu$ n'a pas d'atome, la mesure de $I$ est nulle, ce qui signifie que $I$ est contenu dans une 
composante de $X \setminus \Lambda$. Pour tout $n \in \mathbb{N}$ l'ext\'erieur de $I$ dans $I_n$ est 
form\'e d'au plus deux intervalles $J_n^1$ et $J_n^2$ dont les longueurs tendent vers $0$. De plus, 
nous avons 
$$\nu \big( g_n(J_n^1 \cup J_n^2) \big) \geq \varepsilon_0.$$
Ainsi, l'un d'entre eux, par exemple $J_n^1$, v\'erifie $\nu \big( g_n(J_n^1 ) \big) \geq \varepsilon_0/2$, 
et sa longueur tend vers $0$ lorsque $n$ tend vers l'infini. Si les longueurs des $I_n$ ne sont pas 
minor\'ees par une constante strictement positive, on peut en extraire une sous-suite d'intervalles 
dont les longueurs tendent vers $0$ et tels que leurs images respectives sont de mesure 
sup\'erieure ou \'egale \`a $\varepsilon _0$. Dans tous les cas, nous avons trouv\'e une 
suite d'intervalles $J_m$ et des \'el\'ements $g_m$ de $\Gamma$ d\'efinis sur $J_m$ tels que 
$$\lim_{m \rightarrow \infty} |J_m| = 0 \qquad  \mbox{ et } 
\qquad \nu \big( g_m(J_m) \big) \geq \varepsilon_0/2 \quad \mbox{pour tout } m \in \mathbb{N}.$$  
De plus, ces intervalles $J_m$ peuvent \^etre pris ferm\'es et avec leurs extr\'emit\'es 
appartenant \`a $\Lambda$. Quitte \`a extraire une sous-suite de la suite d'intervalles 
$g_m(J_m)$, nous pouvons supposer qu'ils contiennent tous dans leur int\'erieur un intervalle 
$J$ de mesure sup\'erieure ou \'egale \`a $\varepsilon _0/2$. Puisque toutes les orbites 
sont denses dans $\Lambda$ et que $\Lambda$ est compact, il existe des intervalles 
$U_1,\ldots, U_k$ de $X$ recouvrant $\Lambda$, ainsi qu'une famille $\{ h_1,\ldots,h_k \}$ 
d'\'el\'ements de $\Gamma$, tels que chaque $h_i$ est d\'efini sur $U_i$ et $h_i(U_i) \!\subset\! J$ 
pour tout $i \!\in\! \{1,\ldots,k\}$. Lorsque $m$ est assez grand, les intervalles $J_m$ sont 
de longueurs assez petites pour qu'ils soient enti\`erement contenus dans l'un des intervalles 
$U_i$, disons $U_{i_m}$. La transformation $g = h_{i_m} \circ g_m ^{-1} $ est alors d\'efinie 
sur l'intervalle ferm\'e $g_m (J_m)$ et son image est contenue dans l'int\'erieur de 
$g_m (J_m)$. \'Ecrivons $g_m(J_m)= [p,p']$, avec $p$ et $p'$ dans $\Lambda$. Nous avons 
donc $p<g(p)<g(p')<p'$. L'orbite du point \esp $x= \lim g^n (p)$ \esp est alors de type 
ressort, puisque d'une part les it\'er\'es de tout point $q$ de $[p,x[$ par $g$ 
convergent vers $x$, et d'autre part l'ensemble $]p,x[ \cap \Lambda$ est non 
vide (\`a cause de la minimalit\'e de $\Lambda$). \esp $\square$

\vspace{0.1cm}

\begin{lem} Si la restriction de $\Gamma$ \`a $\Lambda$ est \'equicontinue, alors il existe une mesure 
de probabilit\'e support\'ee sur $\Lambda$ qui est invariante par tous les \'el\'ements de $\Gamma$. 
\end{lem} 

\noindent{\bf D\'emonstration.} La preuve repose sur les id\'ees de Weil pour construire la mesure 
de Haar d'un groupe localement compact~\cite{We}. Notons $C_{\Lambda}(X)$ l'ensemble des fonctions 
continues sur $X$ qui sont constantes sur les composantes connexes de $X \setminus \Lambda$, et 
d\'esignons par $C^+_{\Lambda}(X)$ le sous-ensemble des fonctions positives et non identiquement 
nulles\footnote{Rappelons 
que, en fran\c{c}ais, un nombre $a$ est dit {\em positif} (resp. {\em n\'egatif}) lorsque 
$a \geq 0$ \esp (resp. $a \leq 0$)}. Si $g$ est un 
\'el\'ement de $\Gamma$ et $\psi$ un \'el\'ement de $C_{\Lambda}^+(X)$, notons $\psi \circ g^{-1}$ la fonction 
qui vaut $\psi \circ g^{-1} $ sur l'image de $g$ et $0$ en dehors. Si le domaine de d\'efinition de $g$ 
contient le support de $\psi$, alors la fonction $\psi \circ g^{-1}$ est \`a nouveau continue. Lorsque la 
restriction de $\Gamma$ \`a $\Lambda$ est \'equicontinue, nous allons construire une fonctionnelle 
$L: C^+_{\Lambda}(X) \rightarrow ]0,\infty[$ qui est lin\'eaire et homog\`ene, et qui v\'erifie 
$$L(  \zeta \circ g^{-1} ) = L (\zeta)$$
pour toute fonction $\zeta \!\in\! C^+_{\Lambda}(X)$ et tout $g \!\in\! \Gamma$ dont le domaine de 
d\'efinition contient le support de $\zeta$. Le lecteur pourra s'assurer que cela d\'efinit une 
mesure de probabilit\'e support\'ee sur $\Lambda$ qui est invariante par $\Gamma$. 

Soient $\psi,\zeta$ des \'el\'ements de $C^+_{\Lambda}(X)$. Puisque $\psi$ est non nulle sur $\Lambda$, 
il est ais\'e de trouver un nombre fini de r\'eels positifs $c_i$, ainsi que des \'el\'ements 
$g_i$ de $\Gamma$ dont le domaine de d\'efinition soit connexe, de telle sorte que
$$\zeta \leq  \sum_i c_i \esp \psi \circ g_i^{-1}.$$
Nous d\'efinissons $(\zeta : \psi)$ comme \'etant l'infimum des valeurs des $\sum_i c_i$ parmi 
tous choix possibles des $c_i$ pour lesquels une telle in\'egalit\'e a lieu. Nous laissons au 
lecteur la v\'erification des propri\'et\'es suivantes~:

\vspace{0.08cm}

\noindent -- \esp $(\zeta : \psi) \geq \frac{\| \zeta \|_{\infty}}{ \| \psi \|_{\infty}}$~;

\vspace{0.08cm}

\noindent -- \esp $(c \zeta : d \psi)= \frac{c}{d} (\zeta : \psi)$ \esp pour toute paire 
de r\'eels strictement positifs $c,d$~;

\vspace{0.08cm}

\noindent -- \esp $(\zeta_1 + \zeta_2 : \psi ) \leq (\zeta_1 : \psi) + (\zeta_2 : \psi)$~;

\vspace{0.08cm}

\noindent -- \esp $(\zeta : \psi) \leq (\zeta : \xi) (\xi : \psi)$~;

\vspace{0.08cm}

\noindent -- \esp $(\zeta \circ g^{-1} : \psi) \leq (\zeta : \psi)$, \esp avec 
\'egalit\'e si le domaine de $g$ contient le support de $\psi$.  

\vspace{0.08cm}

Pour chaque $\psi \in C^+_{\Lambda}(X)$ consid\'erons la fonctionnelle 
$L_{\psi}: C^+_{\Lambda}(X) \rightarrow \mathbb{R}$ donn\'ee par
$$L_{\psi} (\zeta) = \frac{(\zeta : \psi)}{(1 : \psi)}.$$
Ces fonctionnelles sont positives, homog\`enes et sous-additives. 
Nous allons v\'erifier que lorsque le diam\`etre du support de $\psi$ tend 
vers $0$, et que la restriction de $\Gamma$ \`a $\Lambda$ est \'equicontinue, 
ces fonctionnelles deviennent de plus en plus lin\'eaires. Choisissons deux fonctions 
$\xi_1$ et $\xi_1$ dans $C^+_{\Lambda}(X)$ telles que $\xi_1 + \xi_2 =1$. Nous allons d\'emontrer que 
$$\big| L_{\psi} (\xi_1 \zeta) + L_{\psi} (\xi_2  \zeta) - L_{\psi} (\zeta) \big|$$
tend vers $0$ lorsque le diam\`etre du support de $\psi$ tend vers $0$. Pour cela consid\'erons un 
r\'eel $\eta \!>\! 0$. Puisque $\xi_1$ et $\xi_2$ sont continues et que $X$ est compact, il existe 
$\varepsilon>0$ tel que si $dist (x,x') \leq \varepsilon$ alors $|\xi_1 (x) - \xi_1(x')| < \eta$ 
et $|\xi_2 (x) - \xi_2 (x')| < \eta$. De plus, puisque la restriction de $\Gamma$ \`a $\Lambda$ 
est \'equicontinue, il existe $\delta > 0$ tel que si $I$ est un intervalle v\'erifiant 
$\nu(I) \leq \delta$, alors pour tout $g \in \Gamma$ dont le domaine de d\'efinition 
est connexe on a $\nu \big( g(I \cap dom(g)) \big) < \varepsilon$. Supposons que le 
support de $\psi$ soit contenu dans un intervalle de longueur plus petite que $\delta$. 
Consid\'erons une famille de r\'eels positifs $c_i$ et des $g_i \!\in\! \Gamma$ dont le 
domaine de d\'efinition soit connexe et tels que 
$$\zeta \leq \sum_i c_i \esp \psi \circ g_i ^{-1}.$$
Si le domaine de d\'efinition de $g_i$ n'intersecte pas 
le support de $\psi$, alors la fonction $\psi \circ g_i^{-1}$ 
est identiquement nulle, et on peut l'enlever de la somme en conservant l'in\'egalit\'e. Nous 
supposerons donc que le domaine de $g_i$ intersecte le support de $\psi$. Choisissons un point 
$y_i$ dans $g_i \big( supp (\psi) \cap dom(g_i) \big)$. Pour tout $x$ dans $X$ nous avons 
$$\zeta(x) \xi_1(x) \leq \sum_i c_i \esp \xi_1 (x) \esp \psi \circ g_i^{-1} (x) 
= \sum_i {}' \esp c_i \esp \xi_1 (x) \esp \psi \big( g_i^{-1} (x) \big),$$
o\`u la deuxi\`eme somme ne porte que sur les termes pour lesquels $g_i^{-1} (x)$ appartient au 
support de $\psi$, c'est-\`a-dire tels que \esp $x \!\in\! g_i \big( supp(\psi) \cap dom(g_i) \big)$. 
\esp D'apr\`es l'hypoth\`ese d'\'equicontinuit\'e nous avons 
$dist \big( g_i^{-1} (x), y_i \big) < \varepsilon$, et donc 
$$\big|  \xi_1 (y_i) - \xi_1 ( g_i^{-1} (x)) \big| \leq \eta.$$
Par suite,
$$\xi_1 \zeta \leq \sum_i  c_i \esp \big( \xi(y_i) + \eta \big) \esp \psi \circ g_i ^{-1},$$
ce qui implique l'in\'egalit\'e 
$$(\xi_1 \zeta : \psi)  \leq  \sum_i c_i \esp \big( \xi_1 (y_i) + \eta \big).$$
Puisque nous avons les m\^emes in\'egalit\'es pour $\xi_2$, nous obtenons 
$$(\xi_1 \zeta:\psi) + (\xi_2 \zeta:\psi) \leq  \sum_i c_i \esp 
\big( \xi_1 (y_i) + \xi_2 (y_i) + 2 \eta \big) = (1 + 2\eta) \sum_i c_i.$$
En passant \`a l'infimum sur les $\sum_i c_i$, et puis en divisant par $(1:\psi)$, 
nous en d\'eduisons que 
$$L_{\psi} (\xi_1 \zeta) + L_{\psi} (\xi_2  \zeta) \leq (1 + 2\eta) \esp L_{\psi} (\zeta)$$
lorsque le diam\`etre du support de $\psi$ est contenu dans un intervalle de taille inf\'erieure 
\`a $\delta$. Comme nous savons d\'ej\`a que $L_{\psi}$ est sous-additive, ceci montre que 
$$\big| L_{\psi} (\xi_1 \zeta) + L_{\psi} (\xi_2  \zeta) - L_{\psi} (\zeta) \big|$$
tend vers $0$ lorsque le diam\`etre du support de $\psi$ tend vers $0$. 

Pour conclure la d\'emonstration, il nous suffit de trouver une suite de fonctions $\psi_n$ dans 
$C^+_{\Lambda}(X)$ telle que les fonctionnelles $L_{\psi_n}$ convergent simplement vers une 
fonctionnelle $L$. Pour cela, remarquons que pour tout $\psi,\zeta$ dans $C^+_{\Lambda}(X)$ on a 
$$\frac{1}{(1:\zeta)} \leq L_{\psi}(\zeta) \leq (\zeta:1).$$
Ainsi, les fonctionnelles $L_{\psi}$ d\'efinissent un point 
$(L_{\psi}(\zeta))_{\zeta \in C^+_{\Lambda}(X)}$ dans l'espace produit
$$\Pi = \prod _{\zeta \in C^+_{\Lambda}(X) } \left[ \frac{1}{(1:\zeta)},(\zeta:1) \right].$$
Cet espace $\Pi$ muni de la topologie produit est compact. Donc, si $\psi_n$ est une suite 
d'\'el\'ements de $C^+_{\Lambda}(X)$ dont les supports tendent vers un point, de la suite d'\'el\'ements 
$(L_{\psi_n}(\zeta))_{\zeta \in C^+_{\Lambda}(X)}$ de $\Pi$ on peut extraire une sous-suite qui converge 
vers un \'el\'ement $\big( L(\zeta) \big)_{\zeta \in C^+_{\Lambda}(X)}$. D'apr\`es ce qui pr\'ec\`ede, 
$L$ d\'efinit une fonctionnelle sur $C^+_{\Lambda}(X)$ qui v\'erife les conditions suivantes~:

\vspace{0.08cm} 

\noindent -- \esp $L$ est strictement positive~;

\vspace{0.08cm}

\noindent -- \esp pour tout $\zeta \in C^+_{\Lambda}(X)$ et tout $g \in \Gamma$ tel que $dom(g)$ 
contient le support de $\zeta$ on a $L(\zeta \circ g^{-1} ) = L(\zeta)$~; 

\vspace{0.08cm}

\noindent -- \esp $L$ est homog\`ene~;

\vspace{0.08cm}

\noindent -- pour toutes fonctions $\zeta,\xi_1,\xi_2$ de $C^+_{\Lambda}(X)$ qui v\'erifient 
$\xi_1 + \xi_2= 1$ on a $L(\xi_1 \zeta)+ L(\xi_2 \zeta) = L(\zeta)$.

\vspace{0.08cm}

Il n'est pas difficile de s'assurer que les deux derni\`eres conditions impliquent que $L$ est 
en fait lin\'eaire. Ceci ach\`eve la preuve du lemme, et donc celle de la proposition 
\ref{sacc}. \esp $\square$

\vspace{0.3cm}

La terminologie 
\begin{tiny}$^{_\ll}$\end{tiny}ressort\hspace{0.05cm}\begin{tiny}$^{_\gg}$\end{tiny} vient de la th\'eorie 
des feuilletages (voir par exemple~\cite{GLW} pour les d\'etails ainsi qu'un joli dessin d'un ressort). 
Dans ce contexte ce qui pr\'ec\`ede se traduit par le r\'esultat suivant, \`a comparer avec \cite{IT}.

\vspace{0.1cm}

\begin{prop} {\em Soit $\mathcal{F}$ un feuilletage de codimension 1 d'une vari\'et\'e 
compacte. Si $\mathcal{F}$ n'admet pas de mesure (de probabilit\'e) transverse 
invariante, alors il contient une feuille ressort (topologique).} 
\end{prop}

\noindent{\bf D\'emonstration.} Quitte \`a passer au rev\^etement double orientable, nous pouvons supposer 
que $\mathcal{F}$ est transversalement orientable. Il suffit alors d'appliquer la conclusion de la 
proposition \ref{sacc} au pseudo-groupe d'holonomie de $\mathcal{F}$. \esp $\square$


\subsection{Contraction presque s\^ure pour les groupes d'hom\'eomorphismes du cercle}
\label{ap1}

Une \'elaboration plus simple et conceptuelle des id\'ees du paragraphe pr\'ec\'edent peut \^etre faite 
pour les groupes d'hom\'eomorphismes du cercle. En suivant Ghys \cite{ghys} (pages 360-362 ; voir aussi 
\cite{mar-CRAS}), consid\'erons un tel groupe $\Gamma$ et supposons d'abord que ses orbites soient 
denses. Si son action sur le cercle est {\em \'equicontinue} (dans le sens que pour tout 
$\varepsilon > 0$ il existe $\delta > 0$ tel que si $dist (x,y) \leq \delta$ 
alors $dist \big( g(x),g(y) \big) \leq \varepsilon$ pour tout $g \in \Gamma$), 
alors $\Gamma$ est topologiquement conjugu\'e \`a un groupe 
de rotations. Sinon, c'est que son action est {\em expansive}, dans le sens que 
pour tout $x \in \clo$ il existe un intervalle ouvert $I$ contenant $x$ et une suite 
$(h_n)$ d'\'el\'ements de $\Gamma$ tels que la taille de $h_n(I)$ tend vers z\'ero. 
En particulier, il existe des points $y \in \clo$ tels que l'intervalle $[x,y[$ est 
\begin{tiny}$^{_\ll}$\end{tiny}contractable\hspace{0.05cm}\begin{tiny}$^{_\gg}$\end{tiny}. 
Si l'on denote par $\theta(x)$ le 
\begin{tiny}$^{_\ll}$\end{tiny}supremum\hspace{0.05cm}\begin{tiny}$^{_\gg}$\end{tiny} 
de ces points, alors l'application $x \mapsto \theta(x)$ s'av\`ere \^etre un 
hom\'eomorphisme d'ordre fini $\kappa (\Gamma)$ qui commute avec tous les \'el\'ements 
de $\Gamma$ (le nombre $\kappa (\Gamma)$ 
est appel\'e le {\em degr\'e de} $\Gamma$). Par suite, si l'on 
identifie les points des orbites par $\theta$, alors on obtient un cercle topologique 
$\clo /\!\sim$ sur lequel le groupe $\Gamma$ agit par hom\'eomorphismes. Remarquons 
que cette derni\`ere action n'est pas n\'ecessairement fid\`ele. Cependant, elle 
v\'erifie une {\em propri\'et\'e d'expansivit\'e forte}, \`a savoir tout intervalle 
dont le compl\'ementaire n'est ni vide ni r\'eduit \`a un seul point est contractable 
par une suite d'\'el\'ements du groupe.

Si les orbites de notre groupe original ne sont pas denses, alors il peut se pr\'esenter 
deux cas : soit il poss\`ede une orbite finie, soit il existe un (unique) ensemble de 
Cantor invariant minimal (et sur lequel les orbites de tous les points du cercle 
s'accumulent). En \'ecrasant les composantes connexes du compl\'ementaire de cet 
ensemble, ce dernier cas peut \^etre ramen\'e \`a celui d'orbites denses.
 
Cherchons maintenant \`a donner des versions probabilistes de ce qui pr\'ec\`ede. 
Pour simplifier, supposons que $\Gamma$ soit d\'enombrable, et fixons une mesure 
de probabilit\'e $\mu$ sur $\Gamma$ qui soit {\em non d\'eg\'en\'er\'ee}, c'est-\`a-dire 
telle que son support engendre $\Gamma$ en tant que semigroupe. Une mesure de 
probabilit\'e $\nu$ sur le cercle est dite {\em stationnaire} (par rapport \`a 
$\mu$) si $\mu * \nu = \nu$, c'est-\`a-dire si pour toute fonction continue $\psi$ 
d\'efinie sur le cercle on a 
$$\int_{\clo} \psi(x) \thinspace d \nu(x) = 
\int_{\Gamma} \int_{\clo} \psi \big( g(x) \big) \thinspace 
d \nu(x) \thinspace d \mu(g).$$
Pour mieux comprendre cette d\'efinition, consid\'erons l'op\'erateur 
de diffusion agissant sur les fonctions continues du cercle par la formule 
\begin{equation}
D \psi (x) = \int_{\Gamma} \psi \big( g (x) \big) \thinspace d\mu(g).
\label{difundir}
\end{equation}
Cet op\'erateur agit de mani\`ere duale sur les mesures de probabilit\'e en pr\'eservant 
le compact convexe des mesures de probabilit\'e : cette action duale correspond pr\'ecis\'ement 
\`a celle donn\'ee par \thinspace $\upsilon \mapsto \mu * \upsilon$. \thinspace L'existence 
de (au moins) une mesure stationnaire $\nu$ sur le cercle est ainsi assur\'ee par le 
th\'eor\`eme du point fixe de Kakutani (bien s\^ur, on peut aussi d\'emontrer cela 
par un argument de moyennes de Birkhoff). 

\vspace{0.1cm} 

\begin{lem} {\em Si les orbites de $\Gamma$ sont denses, alors $\nu$ est \`a support 
total et sans atome. Si $\Gamma$ admet un ensemble de Cantor invariant et minimal, alors 
cet ensemble co\"{\i}ncide avec le support de $\nu$, et $\nu$ est encore sans atome.}
\label{soporte}
\end{lem}

\noindent{\bf D\'emonstration.} Montrons d'abord que si $\nu$ poss\`ede des atomes, 
alors $\Gamma$ a des orbites finies (signalons en passant que, contrairemment 
\`a ce que l'on peut croire, ce dernier cas peut \^etre tr\`es compliqu\'e ; 
voir par exemple la proposition \ref{vg} et la remarque \ref{miok} plus loin). 
En effet, si $p$ est un point de mesure (positive et) 
maximale, alors \`a partir de l'\'egalit\'e   
$$\nu(p) = \int_{\Gamma} \nu \big( g^{-1}(p) \big) \thinspace d\mu(g),$$
on conclut que $\nu \big( g^{-1}(p) \big) = \nu(p)$ pour tout $g$ 
dans le support de $\nu$. Ceci reste vrai pour tout \'el\'ement de 
$\Gamma$, car $\mu$ est une mesure non d\'eg\'en\'er\'ee. La masse totale 
de $\nu$ \'etant finie, l'orbite de $p$ par $\Gamma$ doit 
n\'ecessairement \^etre finie.

Si l'action de $\Gamma$ est minimale, alors le support de $\nu$ est 
tout le cercle puisque c'est un ensemble ferm\'e invariant. Si $\Gamma$ 
admet un minimal exceptionnel $\Lambda$ alors cet ensemble \'etant unique, il 
doit \^etre contenu dans le support de $\nu$. Donc, pour montrer que $\Lambda$ 
et $supp (\nu)$ co\"{\i}ncident, nous devons montrer que $\nu(I)\!=\!0$ pour 
toute composante connexe de $\clo \setminus \Lambda$. Or, si ce n'est 
pas le cas, en prenant une telle composante de masse maximale, on conclut 
(par un argument analogue \`a celui donn\'e dans le cas d'existence 
d'atomes) que l'orbite de $I$ est finie. N\'eanmoins, ceci est absurde, car 
les orbites des extr\'emit\'es de $I$ sont denses dans $\Lambda$. \esp $\square$

\vspace{0.25cm}

D\'esignons par $\Omega$ l'espace des suites $(g_1, g_2, \ldots) \in \Gamma^{\mathbb{N}}$ 
(muni de la mesure $\mathbb{P} = \mu^{\mathbb{N}}$). Si l'on d\'esigne par $\sigma$ 
le d\'ecalage sur $\Omega$, alors on v\'erifie ais\'ement qu'une mesure de probabilit\'e 
$\nu$ sur le cercle est stationnaire par rapport \`a $\mu$ si et seulement si la mesure 
$\mathbb{P} \times \nu$ est invariante par le {\em produit crois\'e} \esp \esp 
$T: \Omega \times \clo \rightarrow \Omega \times \clo$ \esp d\'efini par
$$T(\omega,x) = \big( \sigma(\omega),h_1(\omega)(x) \big) = 
\big( \sigma(\omega), g_1 (x) \big), \qquad \omega = (g_1,g_2,\ldots).$$
En suivant Furstenberg \cite{furst}, pour \'etudier l'\'evolution des 
compositions al\'eatoires on consid\`ere le {\em processus inverse} 
\thinspace $\bar{h}_n(\omega) = g_1 \cdots g_n$. 
\thinspace Si $\psi$ est une fonction continue d\'efinie sur le cercle et $\nu$ 
est une mesure stationnaire, alors la suite de variables al\'eatoires 
$$\xi_n(\omega) = \int_{\clo} \psi \thinspace d \big(g_1 \cdots g_n(\nu)\big)$$ 
est une martingale. Le th\'eor\`eme de convergence de martingales (plus un 
argument de densit\'e) montre ainsi que, pour un sous-ensemble de probabilit\'e 
total $\Omega_0$ de $\Omega$, la limite suivante (par rapport 
\`a la topologie faible) existe :
$$\lim_{n \rightarrow \infty} g_1 g_2 \cdots g_n (\nu) = 
\lim_{n \rightarrow \infty} \bar{h}_n(\omega)(\nu) = \omega(\nu).$$
De plus, l'application (bien d\'efinie presque partout) $\omega \mapsto \omega(\nu)$ 
est mesurable (voir \cite{margulis-libro}, page 199). Une proposition \'equivalente \`a 
celle pr\'esent\'ee \`a continuation a \'et\'e originalement d\'emontr\'ee par Antonov 
dans \cite{anton} (voir aussi \cite{victor}) ; elle est \`a rapprocher avec 
des r\'esultats contenus dans \cite{furst,ka-anals,woess} et le \S 2 du 
chapitre VI de \cite{margulis-libro}.

\vspace{0.1cm}

\begin{prop} {\em Supposons que $\Gamma$ soit d\'enombrable, que ses orbites soient denses 
et que la propri\'et\'e d'expansivit\'e forte soit v\'erifi\'ee. Alors pour presque toute 
suite $\omega$ dans $\Omega_0$, la mesure $\omega(\nu)$ est une mesure de Dirac.}
\label{cont-top}
\end{prop}

\noindent{\bf D\'emonstration.} Nous montrerons que pour tout $\varepsilon\!\in ]0,1]$ 
il existe un sous-ensemble de probabilit\'e total $\Omega_{\varepsilon}$ de $\Omega_0$ 
tel que pour tout $\omega \in \Omega_{\varepsilon}$ il existe un intervalle $I$ de 
longueur $\leq \varepsilon$ tel que $\omega(\nu)(I) \geq 1 - \varepsilon$. 
Ceci permet de conclure de mani\`ere \'evidente que, pour tout $\omega$ dans 
l'ensemble \thinspace $\Omega^* = \cap_{n \in \mathbb{N}} \Omega_{1/n}$ \thinspace 
(dont la probabilit\'e est totale), la mesure $\omega(\nu)$ est une mesure de Dirac.

Fixons donc $\varepsilon > 0$. Pour chaque $n \in \mathbb{N}$ d\'esignons par 
$\Omega^{n,\varepsilon}$ l'ensemble des suites $\omega \in \Omega_0$ telles que, 
pour tout $m \geq 0$ et tout intervalle $I$ de $\clo$ de longueur $\leq \varepsilon$, 
on a $\bar{h}_{n+m}(\omega)(\nu)(I) < 1-\varepsilon$. Nous devons d\'emontrer que 
la mesure de $\Omega^{n,\varepsilon}$ est nulle. Pour cela, nous allons commencer par 
exhiber un sous-ensemble fini $\mathcal{G}_{\varepsilon}$ de $supp (\mu)$ ainsi 
qu'un entier $l \in \mathbb{N}$ tels que, pour tout $r \in \mathbb{N}$ et 
tout $(g_1,\cdots,g_r) \in \Gamma^r$, il existe un intervalle $I$ de longueur 
$\leq \varepsilon$, un entier $\ell \leq 2l$, et des \'el\'ements 
$f_1,\ldots,f_{\ell}$ dans $\mathcal{G}_{\varepsilon}$ v\'erifiant 
\begin{equation}
g_1 \cdots g_r f_1 \cdots f_{\ell} (\nu) (I) \geq 1-\varepsilon.
\label{aprobar}
\end{equation}

Fixons deux points distincts $a$ et $b$ du cercle, ainsi qu'un entier 
$q > 1 / \varepsilon$, et prenons $q$ points distincts $a_1,\ldots,a_q$ de 
l'orbite de $a$ par $\Gamma$. Pour chaque $i \! \in \! \{1,\ldots,q\}$ fixons 
un \'el\'ement $h_i \in \Gamma$ et un intervalle ouvert $U_i$ contenant $a_i$ de 
fa\c{c}on \`a ce que les $U_i$ soient deux \`a deux disjoints, $h_i(a) = a_i$ 
et $h_i(U) = U_i$ pour un certain voisinage $U$ de $a$ qui ne contient pas $b$. 
Prenons maintenant un voisinage $V$ de $b$ disjoint de $U$ et tel que 
$\nu(\clo \setminus V) \geq 1 - \varepsilon$. Par la minimalit\'e et la 
propri\'et\'e d'expansivit\'e forte, il existe $h \! \in \!\Gamma$ tel 
que $h(\clo \setminus V) \subset U$. 
Chaque \'el\'ement dans $\{h_1,\ldots,h_q,h\}$ peut \^etre 
\'ecrit comme un produit d'\'el\'ements du support de $\mu$. 
Cela peut \^etre fait de plusieurs mani\`eres diff\'erentes, 
mais si l'on fixe une fois pour toutes une \'ecriture pour chaque \'el\'ement, alors 
l'ensemble $\mathcal{G}_{\varepsilon}$ des \'el\'ements de $supp (\mu)$ qui sont utilis\'es 
est fini. Soit $l$ le nombre maximal d'\'el\'ements qui apparaissent 
dans l'une des \'ecritures pr\'ec\'edentes. Pour v\'erifier (\ref{aprobar}), notons que 
pour $g = g_1 \cdots g_r$ les intervalles $g(U_i)$ sont deux \`a deux disjoints, et donc 
la longueur d'au moins l'un d'entre eux doit \^etre major\'ee par $\varepsilon$. Si l'on fixe un 
tel intervalle $I = g(U_i)$ alors on obtient
$$g_1 \cdots g_r h_i h (\nu) (I) = \nu \big( h^{-1} (U) \big) 
\geq \nu(\clo \setminus V) \geq 1 - \varepsilon,$$
ce qui ach\`eve la v\'erification de (\ref{aprobar}). 

Notons maintenant $\rho = \min \{\mu(f): f \in \mathcal{G}_{\varepsilon} \}$ et posons 
\vspace{-0.2cm}
$$\Omega^{\varepsilon}_{n + m} = \{ \omega \in \Omega_0: \mbox{ pour tout intervalle } I 
\mbox{ de longueur } |I| \leq \varepsilon \mbox{ et tout } k \leq m 
\mbox{ on a } \bar{h}_{n+k}(\omega)(\nu)(I) < 1 - \varepsilon\}.$$

\vspace{-0.2cm}

\noindent D'apr\`es (\ref{aprobar}) nous avons $\mathbb{P}(\Omega_{n+2 l t}^{\varepsilon}) 
\leq (1-\rho^{2 l})^t$. Donc, en passant \`a la limite lorsque $t$ tend vers l'infini, 
on conclut que $\mathbb{P}(\Omega^{n,\varepsilon}) = 0$, ce qui permet de finir la 
d\'emonstration. \esp $\square$

\vspace{0.1cm}

\begin{rem}
Il est tr\`es int\'eressant de remarquer que les seules propri\'et\'es de 
$\nu$ que l'on a utilis\'e dans la preuve pr\'ec\'edente sont le fait que 
son support est total et que, pour presque tout $\omega \in \Omega$, la suite 
de mesures de probabilit\'e $g_1 \cdots g_n (\nu)$ converge faiblement.
\end{rem}

On d\'efinit le {\em coefficient de contraction} $c(h)$ d'un hom\'eomorphisme $h$ du cercle comme 
\'etant l'infimum parmi tous les $\varepsilon > 0$ tels qu'il existe deux intervalles ferm\'es 
$I$ et $J$ du cercle de taille au plus $\varepsilon$ et tels que $h(\overline{\clo \setminus I}) 
= J$. Cette d\'efinition permet de donner une \begin{tiny}$^{_\ll}$\end{tiny}version 
topologique\hspace{0.05cm}\begin{tiny}$^{_\gg}$\end{tiny} 
de la proposition pr\'ec\'edente pour les compositions dans 
\begin{tiny}$^{_\ll}$\end{tiny}l'ordre naturel\hspace{0.05cm}\begin{tiny}$^{_\gg}$\end{tiny}.

\vspace{0.08cm}

\begin{prop} {\em Sous les hypoth\`eses de la proposition} \ref{cont-top}, {\em 
pour tout $\omega = (g_1,g_2,\ldots) \in \Omega^*$ le coefficient de contraction de 
$h_n(\omega)\!=\! g_n \cdots g_1$ converge vers z\'ero lorsque $n$ tend vers l'infini.}
\label{ultima}
\end{prop}

\noindent{\bf D\'emonstration.} Comme $\nu$ est sans atome et son support est total, il existe 
un hom\'eomorphisme $\varphi$ du cercle qui envoit $\nu$ sur la mesure de Lebesgue. Par suite, 
et puisque l'affirmation \`a d\'emontrer est de nature purement topologique, nous pouvons 
supposer que $\nu$ co\"{\i}ncide avec la mesure de Lebesgue. 

D\'esignons par $\bar{\mu}$ la mesure sur $\Gamma$ d\'efinie par $\bar{\mu}(g) = \mu(g^{-1})$, 
et notons $\bar{\Omega}$ l'espace de probabilit\'e $\Gamma^{\mathbb{N}}$ muni de la mesure 
$\bar{\mu}^{\mathbb{N}}$. Sur cet espace consid\'erons le processus 
\esp $\bar{h}_n(\bar{\omega})=g_1 \cdots g_n$, \esp o\`u $\bar{\omega}= (g_1,g_2,\ldots)$. 
\esp D'apr\`es la proposition 
\ref{cont-top}, pour tout $\varepsilon > 0$ il existe $n(\varepsilon) \in \mathbb{N}$ tel 
que si $n \geq n(\varepsilon)$ alors il existe un intervalle (ferm\'e) $I$ tel que 
$\nu(I) \leq \varepsilon$ et $\bar{h}_n(\bar{\omega})(\nu)(I) \geq 1 - \varepsilon$. 
Si l'on note \esp $J$ la fermeture de 
$\clo \setminus g_n^{-1} \cdots g_1^{-1}(I)$, \esp 
alors on voit que $|I| = \nu(I) \leq \varepsilon$, 
$$|J| = 1 - |g_n^{-1} \cdots g_1^{-1} (I)| 
= 1 - \nu \big( \bar{h}_n (\bar{\omega})^{-1}(I) \big) 
= 1 - \bar{h}_n(\bar{\omega})(\nu)(I) \leq \varepsilon$$
et \esp $g_n^{-1} \cdots g_1^{-1} (\overline{\clo \setminus I}) \!=\! J$. Par suite, 
\esp $c \big( g_n^{-1} \cdots g_1^{-1} \big) \!\leq\! \varepsilon$ \esp pour tout 
$n \!\geq\! n(\varepsilon)$. La preuve est conclue en remar- quant que la transformation 
\esp $(g_1,g_2,\ldots) \mapsto (g_1^{-1},g_2^{-1},\ldots)$ \esp est un isomorphisme entre 
$(\bar{\Omega},\bar{\mu}^{\mathbb{N}})$ et $(\Omega,\mu^{\mathbb{N}})$. \esp $\square$

\vspace{0.3cm}

Soulignons que le coefficient de contraction est toujours r\'ealis\'e, dans le sens que 
pour tout hom\'eomor- phisme $h$ du cercle il existe des intervalles ferm\'es $I$ et $J$ tels 
que \esp $\max \{ |I|,|J| \} = c(h)$ \esp et \esp $h(\overline{\clo \setminus I}) = J$ 
\esp (ces intervalles ne sont pas n\'ecessairement uniques). En cons\'equence, pour 
tout $\omega \in \Omega^*$ on peut choisir deux suites d'intervalles ferm\'es 
$I_n(\omega)$ et $J_n(\omega)$ dont la taille tend vers z\'ero et tels que 
$h_n(\omega) \big( \overline{\clo \setminus I_n(\omega)} \big) = J_n(\omega)$ 
pour tout $n \in \mathbb{N}$. Pour n'importe quel choix, les 
intervalles $I_n(\omega)$ convergent vers le point $\varsigma_+(\omega)$, tandis 
que (g\'en\'eriquement) les intervalles $J_n(\omega)$ se prom\`enent un peu 
partout sur le cercle.

\begin{rem} Les propositions \ref{cont-top} et \ref{ultima} 
sugg\`erent qu'une version faible du th\'eor\`eme 
d'Oseledets devrait \^etre valable pour la plupart des transformations de fibr\'es en 
cercles. Pour \^etre plus pr\'ecis, consid\'erons un espace m\'etrique quelconque 
$\bar{\Omega}$ muni d'une mesure de probabilit\'e $\bar{\mathbb{P}}$, et soit 
$T : \bar{\Omega} \rightarrow \bar{\Omega}$ un hom\'eomorphisme qui pr\'eserve 
$\mu$ de mani\`ere ergodique. \'Etant donn\'ee une application bor\'elienne 
$L: \bar{\Omega}  \rightarrow \mathrm{Hom\acute{e}o}_+(\clo)$, consid\'erons la 
transformation $\Phi: \bar{\Omega} \times \clo \rightarrow \bar{\Omega} \times \clo$ 
d\'efinie par 
$$\Phi(\omega,x) = ( T(\omega), L(\omega) (x) ),$$
et notons 
$$L_n(\omega)(x) = L(T^{n}(\omega)) \circ \cdots \circ L(T(\omega)) 
\circ L(\omega) (x), \qquad n \geq 0,$$
$$L_n(\omega)(x) = L(T^{n}(\omega))^{-1} \circ \cdots \circ 
L(T^{-2}(\omega))^{-1} \circ L(T^{-1}(\omega))^{-1} (x), \qquad n < 0,$$
Supposons que la transformation $\Phi$ ne pr\'eserve aucune mesure 
de probabilit\'e, et admettons pour simplifier qu'elle soit continue 
et que ses orbites soient denses. Sous ces hypoth\`eses, il devrait exister un 
hom\'eomorphisme d'ordre fini $\theta: \clo \rightarrow \clo$ qui commute avec 
presque tous les $L(\omega)$ et tel que l'action de $\Phi$ 
sur l'espace \thinspace $\bar{\Omega} \times \clo \!/ \!\sim$ \thinspace 
des orbites par l'application $(\omega,x) \mapsto (\omega,\theta(x))$ admette 
deux sections mesurables $\varsigma_+$ et $\varsigma_{-}$ telles que pour 
$\bar{\mathbb{P}}$-presque tout point $\omega \in \bar{\Omega}$ et tout $x \in \clo$ on ait
$$\lim_{n \rightarrow +\infty} dist \big( L_n(\omega)(x),\varsigma_{-}(T^n(\omega) \big) 
= 0 \qquad \lim_{n \rightarrow -\infty} dist \big(L_n(\omega)(x),\varsigma_+(T^n(\omega) \big) 
= 0.$$
\end{rem}


\subsection{Quelques remarques autour de la mesure stationnaire}
\label{ap2}

La proposition \ref{cont-top} permet (par un argument bien connu) de d\'emontrer l'unicit\'e de la 
mesure stationnaire. Signalons que ce r\'esultat est valable dans le cadre beaucoup plus g\'en\'eral 
des feuilletages de codimension 1 sans mesure transverse invariante \cite{DK} (la notion de mesure 
stationnaire dans ce contexte est celle de Garnett ; voir \cite{candel,garnett}).

\vspace{0.1cm}

\begin{prop} {\em Soit $\Gamma$ un groupe d\'enombrable d'hom\'eomorphismes 
du cercle muni d'une mesure de probabilit\'e non d\'eg\'en\'er\'ee $\mu$. 
Si $\Gamma$ ne pr\'eserve aucune probabilit\'e du cercle, alors la 
mesure stationnaire (par rapport \`a $\mu$) est unique.}
\label{unicidad}
\end{prop}

\noindent{\bf D\'emonstration.} Supposons d'abord que l'action de $\Gamma$ 
est minimale et satisfait la propri\'et\'e d'expansivit\'e forte, et fixons 
une mesure $\nu$ sur $\clo$ qui soit stationnaire par rapport \`a $\mu$. 
Pour chaque $\omega \in \Omega$ telle que la limite 
$\lim_{n\rightarrow \infty} \bar{h}_n(\omega) (\nu)$ existe et 
soit une mesure de Dirac, d\'esignons par $\varsigma_{\nu} (\omega)$ 
l'atome de la mesure $\omega(\nu)$, {\em i.e.} le point de $\clo$ 
tel que $\omega(\nu) = \delta_{\varsigma_{\nu}(\omega)}$. L'application 
$\varsigma_{\nu}: \Omega \rightarrow \clo$ est presque partout bien d\'efinie 
et mesurable. Nous affirmons que les mesures $\nu$ et $\varsigma_{\nu}(\mathbb{P})$ 
co\"{\i}ncident. En effet, d'apr\`es la stationnarit\'e de $\nu$,
$$\nu = \mu^{*n} * \nu = \sum_{g \in \Gamma} \mu^{*n}(g) \thinspace g(\nu) = 
\int_{\Omega} \bar{h}_n(\omega) (\nu) \thinspace d \mathbb{P}(\omega).$$
Donc, en passant \`a la limite lorsque $n$ tend vers l'infini (ce qui peut 
facilement \^etre justifi\'e par le th\'eor\`eme de convergence domin\'ee),
$$\nu = \int_{\Omega} \lim_{n \rightarrow \infty} h_n(\omega) (\nu) \thinspace d \mathbb{P}(\omega) 
= \int_{\Omega} \delta_{\varsigma_{\nu}(\omega)} \thinspace d \mathbb{P}(\omega),$$
c'est-\`a-dire \thinspace $\nu = \varsigma_{\nu}(\mathbb{P})$.

Consid\'erons maintenant deux mesures stationnaires $\nu_1$ et $\nu_2$. La mesure 
$\nu = (\nu_1 + \nu_2)/2$ est elle aussi stationnaire, et la fonction $\varsigma_{\nu}$ 
v\'erifie, pour $\mathbb{P}$-presque tout $\omega \in \Omega$, 
$$\frac{\delta_{\varsigma_{\nu_1}(\omega)} + \delta_{\varsigma_{\nu_2}(\omega)}}{2} = 
\delta_{\varsigma_{\nu}(\omega)}.$$
Bien \'evidemment, ceci n'est possible que si $\varsigma_{\nu_1}$ et $\varsigma_{\nu_2}$ co\"{\i}ncident 
presque partout. Par cons\'equent, d'apr\`es l'affirmation de la premi\`ere partie de la preuve, 
$$\nu_1 = \varsigma_{\nu_1}(\mathbb{P}) = \varsigma_{\nu_2}(\mathbb{P}) = \nu_2.$$

Supposons maintenant que $\Gamma$ agisse de fa\c{c}on minimale et expansive (mais non fortement 
expansive), et fixons une probabilit\'e $\nu$ qui soit stationnaire (par rapport \`a $\mu$). 
D'apr\`es le \S \ref{ap1}, il existe un hom\'eomorphisme $\theta: \clo \rightarrow \clo$ 
d'ordre fini, commutant avec (tous les \'el\'ements de) $\Gamma$, et tel que l'action induite 
sur le cercle topologique $\clo\!/\!\!\sim$ obtenu en tant qu'espace d'orbites de $\theta$ 
est (minimale et) fortement expansive. Pour chaque $x \in \clo$ notons \esp $\psi(x) = 
\nu \big( [x,\theta(x)[ \big)$. \esp Comme $\nu$ est sans atome, la fonction $\psi$ est 
continue, et puisque $\theta$ commute avec tous les \'el\'ements de $\Gamma$, elle est 
harmonique. Par suite, l'ensemble des points o\`u $\psi$ prend sa valeur maximale est 
invariant par $\Gamma$. Les orbites de $\Gamma$ \'etant denses, $\psi$ est une fonction 
constante ; en d'autres termes, $\nu$ est invariante par $\theta$. Par ailleurs, $\nu$ 
se projette sur une mesure de probabilit\'e stationnaire pour l'action de $\Gamma$ sur 
$\clo \! / \!\! \sim$. D'apr\`es la premi\`ere partie de la preuve, cette derni\`ere mesure 
est unique, ce qui montre alors l'unicit\'e de $\nu$.

Si $\Gamma$ admet un ensemble de Cantor invariant et minimal, alors cet ensemble 
co\"{\i}ncide avec le support de $\nu$. Si l'on \'ecrase les composantes connexes 
du compl\'ementaire de cet ensemble on obtient une action avec toutes ses orbites 
denses ; lorsque cette action est expansive, on peut appliquer les arguments 
pr\'ec\'edents afin de conclure l'unicit\'e de la mesure stationnaire. Pour 
compl\'eter la d\'emonstration on constate ais\'ement que, dans tous les cas 
qui n'ont pas encore \'et\'e consid\'er\'es, le groupe $\Gamma$ laisse 
invariante une mesure de probabilit\'e du cercle. \esp $\square$

\vspace{0.05cm}

\begin{rem} Sous les hypoth\`eses de la proposition pr\'ec\'edente, il est int\'eressant de 
chercher des conditions sous lesquelles le cercle $\clo$ (l'espace $\clo \!/ \!\sim $ le 
cas \'ech\'eant), muni de la mesure stationnaire, est un bord stochastique maximal pour 
$(\Gamma,\mu)$ ({\em i.e.} il co\"{\i}ncide avec son bord de Poisson-Furstenberg). 
\end{rem}

Comme nous l'avons d\'ej\`a signal\'e, lorsqu'il existe une mesure de probabilit\'e 
invariante on est ramen\'e \`a l'\'etude d'un groupe de rotations ou \`a celle d'un 
groupe avec des orbites finies. Dans le premier cas la mesure stationnaire $\nu$ est 
encore unique et co\"{\i}ncide avec l'unique mesure invariante. En effet, pour chaque 
bor\'elien $X \subset \clo$ la fonction $g \mapsto \nu \big( g^{-1}(X) \big)$ est 
harmonique (par rapport \`a $\mu$), et donc constante d'apr\`es \cite{CD}. 
Le cas o\`u il existe des orbites finies est plus int\'eressant. \`A 
indice fini pr\`es, il se ram\`ene \`a celui d'un groupe agissant sur l'intervalle. 
Or, dans ce contexte la mesure stationnaire ne donne aucune information dynamique 
(tout au moins lorsque la mesure de d\'epart $\mu$ est sym\'etrique). 

\vspace{0.1cm}

\begin{prop} {\em Soit $\Gamma$ un sous-groupe de $\mathrm{Hom\acute{e}o}_+([0,1])$ sans points 
fixe global \`a l'int\'erieur. Si $\mu$ est une mesure de probabilit\'e non d\'eg\'en\'er\'ee 
et sym\'etrique sur $\Gamma$, alors toute mesure de probabilit\'e sur $[0,1]$ qui est 
stationnaire par rapport \`a $\mu$ est support\'ee sur les extr\'emit\'es de $[0,1]$.}
\label{vg}
\end{prop}

\noindent{\bf Premi\`ere d\'emonstration.} Supposons par contradiction que $\overline{\nu}$ 
soit une mesure stationnaire telle que $\overline{\nu}(\{0,1\})<1$. Notons $\nu$ la partie 
de $\overline{\nu}$ support\'ee dans $]0,1[$ convenablement normalis\'ee~: 
c'est encore une mesure stationnaire, 
et elle v\'erifie $\nu(]0,1[)=1$. Comme dans la 
preuve du lemme \ref{soporte}, on montre ais\'ement que $\nu$ est 
sans atome. Donc, quitte \`a \'ecraser les composantes connexes du 
compl\'ementaire du support de $\nu$, nous pouvons supposer que le 
support de $\nu$ co\"{\i}ncide avec $]0,1[$. De plus, en reparam\'etrant 
l'intervalle, nous pouvons supposer que $\nu$ est la mesure de Lebesgue.

Puisque $\nu$ est une mesure stationnaire par rapport \`a la 
probabilit\'e sym\'etrique $\mu$, pour tout $s\!\in ]0,1[$ on a
$$s = \nu \big( [0,s] \big) = 
\int_{\Gamma} \nu \big( g^{-1}([0,s]) \big) \thinspace d\mu(g) 
= \int_{\Gamma} \frac{ \nu \big( g([0,s]) \big) +  
\nu \big( g^{-1}([0,s]) \big)}{2} \thinspace d\mu(g),$$
et donc 
$$s = \int_{\Gamma} \frac{g(s) + g^{-1}(s)}{2} \thinspace d\mu(g).$$
En int\'egrant entre 0 et un point arbitraire $t\!\in ]0,1[$ on obtient 
\begin{equation}
t^2 = \int_{\Gamma} \int_{0}^{t} \big( g(s) + g^{-1}(s) \big) 
\thinspace ds \thinspace d \mu(g).
\label{genial}
\end{equation}
Or, pour tout hom\'eomorphisme $f$ de l'intervalle et tout $t\!\in \![0,1]$ on a
$$\int_0^t \big( f(s) + f^{-1}(s) \big) \thinspace ds \geq t^2,$$
avec l'\'egalit\'e si et seulement si $f(t) = t$ (voir la figure 7). D'apr\`es 
(\ref{genial}), ceci implique que $g(t) = t$ pour tout $g$ appartenant au  
support de $\mu$. Puisque $\mu$ est une mesure non d\'eg\'en\'er\'ee, on conclut 
que $t \!\in ]0,1[$ est un point fixe global pour l'action de $\Gamma$, 
ce qui contredit notre hypoth\`ese. \esp $\square$

\vspace{0.45cm}


\beginpicture

\setcoordinatesystem units <1cm,1cm>

\putrule from 0 0 to 0 5
\putrule from 5 0 to 5 5
\putrule from 0 5 to 5 5
\putrule from 0 0 to 5 0

\plot 
0 0
0.7071 0.1 
1 0.2 
1.5811 0.5 
2 0.8 
2.236 1 
2.5495 1.3 
2.8482  1.6    
3.1622 2
3.3911 2.3 
3.6055 2.6 
3.8729 3 
4.1833 3.5 
4.4721 4 
4.7434 4.5 
5 5 /

\plot 
0 0
0.1 0.7071 
0.2 1  
0.5 1.5811 
0.8 2  
1 2.236  
1.3 2.5495 
1.6 2.8482    
2 3.1622 
2.3 3.3911  
2.6 3.6055  
3 3.8729  
3.5 4.1833  
4 4.4721  
4.5 4.7434  
5 5 /

\plot 0 0 0.1 0.1 /
\plot 0.2 0.2 0.3 0.3 /
\plot 0.4 0.4 0.5 0.5 /
\plot 0.6 0.6  0.7 0.7 /
\plot 0.8 0.8  0.9 0.9 /
\plot 1 1  1.1 1.1 /
\plot 1.2 1.2  1.3 1.3 /
\plot 1.4 1.4  1.5 1.5 /
\plot 1.6 1.6  1.7 1.7 /
\plot 1.8 1.8  1.9 1.9 /
\plot 2 2  2.1 2.1 /
\plot 2.2 2.2  2.3 2.3 /
\plot 2.4 2.4  2.5 2.5 /
\plot 2.6 2.6  2.7 2.7 /
\plot 2.8 2.8  2.9 2.9 /
\plot 3 3  3.1 3.1 /
\plot 3.2 3.2  3.3 3.3 /
\plot 3.4 3.4  3.5 3.5 /
\plot 3.6 3.6  3.7 3.7 /
\plot 3.8 3.8  3.9 3.9 /
\plot 4 4  4.1 4.1 /
\plot 4.2 4.2  4.3 4.3 /
\plot 4.4 4.4  4.5 4.5 /
\plot 4.6 4.6  4.7 4.7 /
\plot 4.8 4.8  4.9 4.9 /

\putrule from 3.5 0 to 3.5 0.2 
\putrule from 3.5 0.4 to 3.5 0.6 
\putrule from 3.5 0.8 to 3.5 1 
\putrule from 3.5 1.2 to 3.5 1.4 
\putrule from 3.5 1.6 to 3.5 1.8 
\putrule from 3.5 2 to 3.5 2.2 
\putrule from 3.5 2.4 to 3.5 2.6 
\putrule from 3.5 2.8 to 3.5 3 
\putrule from 3.5 3.2 to 3.5 3.4 
\putrule from 3.5 3.6 to 3.5 3.8 
\putrule from 3.5 4 to 3.5 4.2 

\putrule from 0 3.5 to 0.2 3.5 
\putrule from 0.4 3.5 to 0.6 3.5 
\putrule from 0.8 3.5 to 1 3.5 
\putrule from 1.2 3.5 to 1.4 3.5 
\putrule from 1.6 3.5 to 1.8 3.5
\putrule from 2 3.5 to 2.2 3.5 
\putrule from 2.4 3.5 to 2.6 3.5 
\putrule from 2.8 3.5 to 3 3.5 
\putrule from 3.2 3.5 to 3.4 3.5 

\putrule from 0 4.18 to 0.2 4.18 
\putrule from 0.4 4.18 to 0.6 4.18 
\putrule from 0.8 4.18 to 1 4.18 
\putrule from 1.2 4.18 to 1.4 4.18 
\putrule from 1.6 4.18 to 1.8 4.18
\putrule from 2 4.18 to 2.2 4.18 
\putrule from 2.4 4.18 to 2.6 4.18 
\putrule from 2.8 4.18 to 3 4.18 
\putrule from 3.2 4.18 to 3.4 4.18 

\put{$f$} at 2.05 1.2
\put{$f^{-1}$} at 1.5 2.3
\put{$A$} at 2.5 0.5
\put{$A$} at 0.45 2.6
\put{$B$} at 2.9 2.3
\put{$\Delta$} at 3.15 3.725 

\put{$0$} at 0 -0.3
\put{$t$} at 3.5 -0.3
\put{$1$} at 5 -0.3
\put{$t$} at -0.6 3.5
\put{$f^{-1}(t)$} at -0.6 4.18 

\put{$A = \int_{0}^{t} f(s) \thinspace ds$} at 8 3.5
\put{$B = \int_0^t f^{-1}(s) \thinspace ds$ } at 8 2.8
\put{$A + B = t^2 + \Delta$} at 8 1.2 


\setdots
 
\putrule from 0.7071 0.1 to 3.5 0.1 
\putrule from 1 0.2 to 3.5 0.2  
\putrule from 1.2 0.3 to 3.5 0.3  
\putrule from 1.4 0.4 to 3.5 0.4  
\putrule from 1.5811 0.5 to 3.5 0.5 
\putrule from 1.8 0.6 to 3.5 0.6 
\putrule from 1.9 0.7 to 3.5 0.7 
\putrule from 2 0.8 to 3.5 0.8 
\putrule from 2.15 0.9 to 3.5 0.9 
\putrule from 2.236 1 to 3.5 1 
\putrule from 2.35 1.1 to 3.5 1.1 
\putrule from 2.45 1.2 to 3.5 1.2 
\putrule from 2.5495 1.3 to 3.5 1.3 
\putrule from 2.67 1.4 to 3.5 1.4 
\putrule from 2.75 1.5 to 3.5 1.5 
\putrule from 2.8482 1.6 to 3.5 1.6 
\putrule from 2.9 1.7 to 3.5 1.7 
\putrule from 3.04 1.8 to 3.5 1.8 
\putrule from 3.1 1.9 to 3.5 1.9 
\putrule from 3.1622 2 to 3.5 2 
\putrule from 3.24 2.1 to 3.5 2.1 
\putrule from 3.32 2.2 to 3.5 2.2 
\putrule from 3.3911 2.3 to 3.5 2.3

\putrule from 0.1 0.7071 to 0.1 3.5  
\putrule from 0.2 1 to 0.2 3.5   
\putrule from 0.3 1.2 to 0.3 3.5   
\putrule from 0.4 1.4 to 0.4 3.5   
\putrule from 0.5 1.5811 to 0.5 3.5  
\putrule from 0.6 1.8 to 0.6 3.5  
\putrule from 0.7 1.9 to 0.7 3.5 
\putrule from 0.8 2 to 0.8 3.5 
\putrule from 0.9 2.15 to 0.9 3.5  
\putrule from 1 2.236 to 1 3.5  
\putrule from 1.1 2.35 to 1.1 3.5  
\putrule from 1.2 2.45 to 1.2 3.5  
\putrule from 1.3 2.5495 to 1.3 3.5  
\putrule from 1.4 2.67 to 1.4 3.5  
\putrule from 1.5 2.75 to 1.5 3.5  
\putrule from 1.6 2.8482 to 1.6 3.5  
\putrule from 1.7 2.9 to 1.7 3.5  
\putrule from 1.8 3.04 to 1.8 3.5  
\putrule from 1.9 3.1 to 1.9 3.5  
\putrule from 2 3.1622 to 2 3.5  
\putrule from 2.1 3.24 to 2.1 3.5  
\putrule from 2.2 3.32 to 2.2 3.5  
\putrule from 2.3 3.3911 to 2.3 3.5 

\putrule from 2.3 3.5 to 2.3 3.3911 
\putrule from  2.5 3.5 to 2.5 3.6055 
\putrule from  2.7 3.5 to 2.7 3.72 
\putrule from  2.9 3.5 to 2.9 3.8 
\putrule from 3.1 3.5 to 3.1 3.9  
\putrule from 3.3 3.5 to 3.3 4.1 

\begin{tiny}
\put{$\cdot$} at 3.37 2.27
\put{$\cdot$} at 2.16 3.41 
\end{tiny}

\put{Figure 7} at 4.8 -0.9

\begin{Large}
\put{$\Downarrow$} at 8 2 
\end{Large}

\put{} at -3.4 0

\endpicture


\vspace{0.35cm}

\noindent{\bf Deuxi\`eme d\'emonstration.} Pour une mesure stationnaire $\nu$ 
consid\'erons la fonction $x \mapsto \nu \big( [0,x] \big)$. Si $\mu$ est 
sym\'etrique, alors cette fonction est harmonique par rapport \`a $\mu$ 
le long de toute orbite. D'apr\`es 
le th\'eor\`eme ergodique de Garnett (Theorem 1. (b) dans \cite{garnett}), cette 
fonction est constante pour $\nu$ presque toute orbite. Autrement dit, pour $\nu$ 
presque tout point $x$ l'\'egalit\'e $\nu \big( [0,g(x)] \big) = \nu \big( [0,x] \big)$ 
est v\'erifi\'ee pour tout $g \in \Gamma$. L'hypoth\`ese de non existence 
de point fixe global pour l'action de $\Gamma$ sur $]0,1[$ implique 
alors que le support de $\nu$ est contenu dans $\{0,1\}$. \esp $\square$

\vspace{0.05cm}

\begin{rem} Comme nous l'avons d\'ej\`a sugg\'er\'e, l'hypoth\`ese de sym\'etrie pour 
la mesure $\mu$ est essentielle pour la validit\'e de la proposition pr\'ec\'edente, 
comme le montre par exemple \cite{ka-intervalo}.
\label{miok}
\end{rem}

\vspace{0.05cm}

Voici une proposition g\'en\'erale de r\'egularit\'e dont la preuve s'appuie sur la 
mesure stationnaire. Signalons en passant que, d'apr\`es \cite{harrison}, un tel 
r\'esultat ne peut \^etre valable qu'en dimension 1. 

\vspace{0.1cm}

\begin{prop} {\em Tout sous-groupe d\'enombrable de $\mathrm{Hom\acute{e}o}_+(\clo)$ 
(resp. de $\mathrm{Hom\acute{e}o}_+([0,1])$) est topologiquement conjugu\'e \`a un 
groupe d'hom\'eomorphismes lipschitziens.}
\end{prop} 

\noindent{\bf D\'emonstration.} Consid\'erons d'abord le cas d'un sous-groupe 
d\'enombrable $\Gamma$ de $\mathrm{Hom\acute{e}o}_+(\clo)$ dont les orbites 
sont denses. Munissons ce groupe d'une mesure de probabilit\'e non 
d\'eg\'en\'er\'ee et sym\'etrique $\mu$, et consid\'erons une mesure de probabilit\'e 
$\nu$ sur le cercle qui soit stationnaire par rapport \`a $\mu$ . Pour chaque 
intervalle $I \subset \clo$ et chaque \'el\'ement $g \in supp (\mu)$ on a
$$\nu(I) = \sum_{g \in supp (\mu)} \nu \big( g^{-1}(I) \big) \thinspace \mu(g)   
\geq \nu \big( g (I) \big) \thinspace \mu(g^{-1}),$$
et donc 
\begin{equation}
\nu \big( g(I) \big) \leq \frac{1}{\mu(g)} \thinspace \nu(I).
\label{lip}
\end{equation}
Prenons maintenant un hom\'eomorphisme $\varphi$ du cercle dans lui-m\^eme 
qui envoie $\nu$ sur la mesure de Lebesgue. Si $J$ est un intervalle quelconque 
du cercle alors, d'apr\`es (\ref{lip}), pour tout $g \in supp (\mu)$ on a
$$|\varphi \!\circ\! g \!\circ\! \varphi^{-1} (J)| = 
\nu \big( g\! \circ\! \varphi^{-1} (J) \big) \leq 
\frac{1}{\mu(g)} \thinspace \nu \big( \varphi^{-1}(J) \big) 
= \frac{1}{\mu(g)} \thinspace |J|.$$ 
Ainsi, pour tout $g$ dans $supp (\mu)$, la transformation 
$\varphi \!\circ\! g\! \circ\! \varphi^{-1}$ est lipschitzienne de 
rapport \thinspace $1/\mu(g)$. \thinspace Puisque la mesure $\mu$ est non 
d\'eg\'en\'er\'ee, cela d\'emontre la proposition dans le cas minimal. 

Si $\Gamma$ est un sous-groupe d\'enombrable quelconque de $\mathrm{Hom\acute{e}o}_+(\clo)$, 
alors en rajoutant une rotation d'angle irrationnel et en consid\'erant le groupe engendr\'e, 
on se ram\`ene au cas minimal. Les arguments plus haut montrent que ce nouveau groupe (et donc 
le groupe originel $\Gamma$) est topologiquement conjugu\'e \`a un groupe d'hom\'eomorphismes 
lipschitziens. Finalement, en identifiant les extr\'emit\'es de l'intervalle $[0,1]$, chaque 
sous-groupe $\Gamma$ de $\mathrm{Hom\acute{e}o}_+([0,1])$ induit un groupe 
d'hom\'eomorphismes du cercle (avec un point fixe global marqu\'e). D'apr\`es ce 
qui pr\`ec\`ede, si $\Gamma$ est d\'enombrable alors ce nouveau groupe est 
conjugu\'e par un \'el\'ement $\varphi$ de $\mathrm{Hom\acute{e}o}_+(\clo)$ 
\`a un groupe d'hom\'eomorphismes lipschitziens du cercle. Pour obtenir 
une vraie conjugaison dans $\mathrm{Hom\acute{e}o}_+([0,1])$, 
il suffit de composer $\varphi$ avec une rotation de fa\c{c}on \`a 
ramener le point marqu\'e du cercle sur lui-m\^eme. \esp $\square$

\vspace{0.25cm}

La proposition pr\'ec\'edente indique que l'on ne doit pas s'attendre 
\`a des r\'esultats de rigidit\'e pour des groupes d'hom\'eomorphismes du 
cercle qui soient sp\'ecifiques \`a une r\'egularit\'e strictement 
comprise entre $C^0$ et $C^{lip}$. Par contre, la compr\'ehension des passages 
$C^{lip}\!\rightarrow C^1$ et $C^1\!\rightarrow C^{1+\tau}$ semble \^etre 
un probl\`eme tr\`es int\'eressant. Une obstruction pour le premier est 
donn\'ee par le th\'eor\`eme de stabilit\'e de Thurston \cite{th}. Par ailleurs, 
la question de savoir si la proposition \ref{gh-ts-tau} est encore valable 
pour des conjugaisons lipschitziennes reste ouverte. Quant au deuxi\`eme 
passage, les r\'esultats de cet article sugg\`erent que le groupe de 
Grigorchuk-Maki consid\'er\'e dans \cite{ceuno} ne devrait pas \^etre 
isomorphe \`a un sous-groupe de $\mathrm{Diff}_+^{1+\tau}([0,1])$ 
pour aucun $\tau\!>\!0$ (m\^eme s'il agit fid\`element par 
diff\'eomorphismes de classe $C^1$ de l'intervalle !).


\subsection{Exposants de Lyapunov et mesures invariantes}
\label{ap3}

Consid\'erons un sous-groupe d\'enombrable $\Gamma$ de $\mathrm{Diff}_+^1(\clo)$. Pour simplifier, 
supposons que $\Gamma$ soit de type fini, et munissons-le d'une mesure de probabilit\'e $\mu$ 
qui soit non d\'eg\'en\'er\'ee et \`a support fini. Soit $\nu$ une probabilit\'e stationnaire 
(par rapport \`a $\mu$)
sur le cercle, et consid\'erons la transformation $T$ de $\Omega \times \clo$ donn\'ee par 
$T(\omega,x)=\big( \sigma(\omega),h_1(\omega)(x) \big)$. Puisque $T$ pr\'eserve 
$\mathbb{P}\!\times\!\nu$, en appliquant le th\'eor\`eme ergodique de Birkhoff \`a la 
fonction \thinspace $(\omega,x) \mapsto \log \big(h_1(\omega)'(x) \big)$ \thinspace 
on conclut que pour $\nu$-presque tout point $x$ du cercle et pour 
$\mathbb{P}$-presque tout chemin al\'eatoire $\omega$ dans $\Omega$, la limite
$$\lambda_{(\omega, x)}(\nu) = 
\lim_{n \rightarrow \infty} \frac{\log \big( h_n(\omega)'(x) \big)}{n}$$
existe~: c'est l'exposant de Lyapunov correspondant au point $(\omega,x)$. Si la mesure 
de probabilit\'e stationnaire $\nu$ est ergodique, c'est-\`a-dire si elle ne peut pas 
\^etre exprim\'ee comme une combinaison convexe de deux probabilit\'es stationnaires 
distinctes, alors la transformation $T$ est ergodique (au sens classique) par rapport 
\`a $\mathbb{P}\!\times\!\nu$ (c'est le th\'eor\`eme ergodique al\'eatoire de Kakutani : 
voir \cite{furman,kifer}). Dans ce cas, l'exposant de Lyapunov est 
$\mathbb{P}\!\times\!\nu$-presque partout constant et \'egal \`a 
$$ \lambda (\nu) = \int_{\Omega} \int _{\clo} \log \big( h_1(\omega)'(x) \big) 
\thinspace d\nu(x) \thinspace d \mathbb{P} (\omega) = 
\int_{\Gamma} \int _{\clo} \log \big( g'(x) \big) 
\thinspace d\nu(x) \thinspace d\mu(g).$$
Lorsque $\mu$ est sym\'etrique et la mesure $\nu$ est invariante par l'action de $\Gamma$, 
l'exposant de Lyapunov $\lambda_{(\omega,x)}(\nu)$ est presque partout nul. En effet, si 
$\nu$ est invariante alors la transformation $S$ de $\Omega \times \clo$ sur lui m\^eme 
d\'efinie par $S(\omega,x) = \big(\omega,h_1(\omega)^{-1}(x)\big)$ pr\'eserve 
$\mathbb{P}\!\times\!\nu$. Ceci implique que 
\begin{eqnarray*}
\lambda(\nu) &=& \int_{\Omega} \int _{\clo} \log \big( h_1(\omega)'(x) \big) 
\thinspace d\nu(x) \thinspace d \mathbb{P} (\omega) 
= \int_{\Omega} \int _{\clo} \log \big( h_1(\omega)'(h_1(\omega)^{-1}(x)) \big) 
\thinspace d\nu(x) \thinspace d \mathbb{P} (\omega) \\
&=& - \int_{\Omega} \int _{\clo} \log \big( (h_1(\omega)^{-1})'(x) \big) 
\thinspace d\nu(x) \thinspace d \mathbb{P} (\omega) = -\int_{\Gamma} \int_{\clo} 
\log \big( (g^{-1})'(x) \big) \esp d\nu(x) \esp d\mu(g) = -\lambda(\nu),
\end{eqnarray*}
o\`u la derni\`ere \'egalit\'e est une cons\'equence de la  sym\'etrie de $\mu$. Par 
suite, si $\nu$ est invariante et ergodique, alors son exposant de Lyapunov est nul. 
Le cas g\'en\'eral s'en d\'eduit par un argument de d\'ecomposition ergodique. 

L'objectif de cet appendice est de donner une d\'emonstration courte et autocontenue 
d'une affirmation r\'eciproque. Le r\'esultat suivant a \'et\'e originalement d\'emontr\'e 
par Baxendale dans un contexte bien plus g\'en\'eral \cite{bax} ; ici nous donnons une preuve 
qui est inspir\'ee de \cite{DK}.

\vspace{0.1cm}

\begin{prop} {\em Si la mesure $\mu$ est (non d\'eg\'en\'er\'ee, \`a support fini et) 
sym\'etrique et $\Gamma$ ne pr\'eserve aucune mesure de probabilit\'e du cercle, alors 
l'exposant de Lyapunov de l'unique mesure stationnaire est strictement n\'egatif.}
\label{neg}
\end{prop}

\vspace{0.1cm}

Pour la d\'emonstration notons
$$\psi (x) = \int_{\Gamma} \log \big( g'(x) \big) \thinspace d\mu(g),$$
et supposons par contradiction que \thinspace $\lambda(\nu) \geq 0$, 
\thinspace c'est-\`a-dire
\begin{equation}
\int_{\clo} \psi(x) \thinspace d \nu(x) \geq 0.
\label{suppp}
\end{equation}
Nous allons d\'emontrer dans ce cas que $\Gamma$ pr\'eserve une mesure de probabilit\'e, 
contredisant ainsi notre hypoth\`ese. Pour cela nous nous appuyons sur un lemme inspir\'e 
des travaux de Sullivan sur les cycles feuillet\'es \cite{Sullivan} (voir \'egalement 
\cite{Gh}). Rappelons que le {\em laplacien} $\Delta \zeta$ d'une fonction r\'eelle et 
continue $\zeta$ est d\'efini par $\Delta \zeta = D \zeta - \zeta$, o\`u $D$ d\'esigne 
l'op\'erateur de diffusion (\ref{difundir}).

\vspace{0.1cm}

\begin{lem} {\em Sous l'hypoth\`ese} (\ref{suppp}), {\em il existe une suite de fonctions 
continues $\zeta_n$ d\'efinies sur le cercle telles que, pour tout entier 
$n\in \mathbb{N}$ et tout point $x$ du cercle,}
\begin{equation}
\psi(x) + \Delta \zeta _n (x) \geq -\frac{1}{n}.
\label{casi-cociclo}
\end{equation}
\end{lem}

\noindent{\bf D\'emonstration.} D\'esignons par $C(\clo)$ l'espace des fonctions continues 
sur le cercle. Notons $E$ le sous-espace constitu\'e des fonctions qui sont des laplaciens de 
fonctions dans $C(\clo)$, et soit $C_+$ le c\^one convexe des fonctions partout positives. Nous 
devons d\'emontrer que si $\psi$ satisfait (\ref{suppp}), alors son image par la projection 
$\pi: C(\clo) \rightarrow C(\clo) / \bar{E}$ est contenue dans $\pi(C_+)$. Supposons que ce 
ne soit pas le cas. Le th\'eor\`eme de s\'eparation de Hahn-Banach donne alors l'existence 
d'une fonctionnelle continue $\bar{L}: C (\clo) / \bar{E} \rightarrow \mathbb{R}$ telle que 
$\bar{L}(\pi(\psi)) < 0 \leq \bar{L} (\pi(\Phi))$ pour tout $\Phi \in C_+$. Bien s\^ur, $\bar{L}$ 
induit une fonctionnelle continue $L: C (\clo) \rightarrow \mathbb{R}$ qui est identiquement 
nulle sur $E$ et telle que $L(\psi) < 0 \leq L(\Phi)$ pour tout $\Phi \in C_+$. Nous affirmons 
que \esp $L = c \esp \nu$ \esp pour certain $c \in \mathbb{R}$ (nous identifions les mesures 
de probabilit\'e aux fonctionnelles lin\'eaires qu'elles d\'efinissent sur l'espace des 
fonctions continues). Pour montrer cela, commen\c{c}ons par remarquer que, 
puisque $L$ est nul sur $E$, pour tout $\zeta \in C (\clo)$ on a 
$$\langle D L, \zeta \rangle = \langle L, D \zeta \rangle = 
\langle L ,\Delta \zeta + \zeta \rangle = \langle L , \zeta \rangle,$$
c'est-\`a-dire que $L$ est invariant par la diffusion. Supposons que la 
d\'ecomposition de Hahn de $L$ s'exprime sous la forme
$$L = \alpha \nu_1 - \beta \nu_2,$$
o\`u $\nu_1$ et $\nu_2$ sont des mesures de probabilit\'e de supports 
disjoints du cercle, $\alpha > 0$ et $\beta > 0$. Dans ce cas, l'\'egalit\'e 
$DL = L$ et l'unicit\'e de la d\'ecomposition de Hahn pour $DL$ montrent que 
$\nu_1$ et $\nu_2$ sont elles aussi invariantes par la diffusion. Par 
suite, \thinspace $\nu_1 = \nu_2 = \nu$, \thinspace ce qui contredit 
le fait que les supports de $\nu_1$ et $\nu_2$ sont disjoints. La 
fonctionnelle $L$ s'exprime donc sous la forme $L = c \esp \upsilon$ pour 
certaine mesure de probabilit\'e $\upsilon$ du cercle ; l'\'egalit\'e 
$L = D L$ donne \'evidemment $\upsilon =\nu$. 

Notons maintenant que, puisque  
$$0 > L(\psi) = c \esp \nu(\psi) = c \esp \int_{\clo} \psi(x) \esp d\nu(x),$$
l'hypoth\`ese (\ref{suppp}) entra\^{\i}ne que $\nu(\psi) > 0$ et $c < 0$. Or, 
comme la fonction constante \'egale \`a $1$ appartient \`a $C_+$, nous avons \esp 
$c = L(1) \geq 0$. \esp Cette contradiction conclut la d\'emonstration. \esp $\square$

\vspace{0.3cm}

Revenons \`a la preuve de la proposition \ref{neg}. Quitte \`a rajouter une constante \`a chaque 
$\zeta_n$, nous pouvons supposer que l'int\'egrale de $\exp(\zeta_n)$ est \'egale \`a $1$ pour 
tout $n \in \mathbb{N}$. Consid\'erons les mesures de probabilit\'e $\nu_n$ sur le cercle 
d\'efinies par
$$\frac{d \nu_n (s)}{ds}  =  \exp \big( \zeta_n (s) \big).$$
Prenons une sous-suite $\nu_{n_i}$ qui converge vers une mesure de probabilit\'e $\bar{\nu}$ 
sur le cercle. Nous allons montrer que $\bar{\nu}$ est invariante par $\Gamma$.

Commen\c{c}ons par d\'emontrer que $\bar{\nu}$ est une mesure harmonique. Pour cela, remarquons 
d'abord que si l'on d\'esigne par $Jac_n(g)$ le jacobien de $g \in \Gamma$ par rapport \`a 
$\nu_n$, alors la relation (\ref{casi-cociclo}) donne, pour tout point $x$ du cercle,
\begin{eqnarray*}
\int_{\Gamma} \log \big( Jac_{n} (g) (x) \big) \thinspace d\mu(g) 
&=& \int_{\Gamma} \log \big( g'(x) \big) \thinspace d\mu(g) + 
\int_{\Gamma} \big[ \zeta_n(g(x)) - \zeta_n(x) \big] \thinspace d \mu(g) \\
&=& \psi(x) + \Delta \zeta_n (x) \geq -\frac{1}{n}.
\end{eqnarray*}
Remarquons maintenant que puisque la diffusion agit contin\^ument sur l'espace des mesures 
de probabilit\'e sur le cercle, la suite de mesures $D \nu_{n_i}$ converge faiblement vers 
la mesure $D \bar{\nu}$. Or, la diffusion de $\nu_n$ est une mesure absolument continue 
par rapport \`a $\nu_n$ dont la densit\'e s'exprime par la formule
$$\frac{d \thinspace D \nu_n(x)}{d \thinspace \nu_n (x)} = 
\int_{\Gamma} Jac_n(g^{-1}) (x) \thinspace d\mu(g) =
\int_{\Gamma} Jac_n(g)(x) \thinspace d\mu(g).$$
Par la concavit\'e de la fonction logarithme nous avons 
$$\frac{d \thinspace D \nu_{n_i}(x)}{d \thinspace \nu_{n_i} (x)} \geq 
\exp \left( \int_{\Gamma} \log \big( Jac_{n_i} (g) (x)\big) 
\thinspace d\mu(g) \right ) \geq \exp(-1/n_i),$$
c'est-\`a-dire \thinspace $D \nu_{n_i} \geq \exp(-1/n_i) \thinspace \nu_{n_i}$. \thinspace 
\`A la limite nous obtenons $D \bar{\nu} \geq \bar{\nu}$, et puisque $\bar{\nu}$ et $D \bar{\nu}$ 
sont des mesures de probabilit\'e, elles sont \'egales~; $\bar{\nu}$ est donc harmonique.

Nous sommes maintenant en mesure de d\'emontrer que $\bar{\nu}$ est invariante par n'importe quel 
\'el\'ement de $\Gamma$. Pour cela, donnons nous un intervalle $I$ tel que $\bar{\nu}(I) > 0$, 
et consid\'erons les fonctions $\psi_{n,I}\!: \Gamma\!\rightarrow ]0,1]$ d\'efinies par 
\thinspace $\psi_{n,I}(g) = \nu_n \big( g (I) \big)$. \thinspace 
Les in\'egalit\'es suivantes montrent que le laplacien de $\log (\psi_{n,I})$ est minor\'e 
par $-1/n$ en l'\'el\'ement neutre $e$ de $\Gamma$~:
\begin{eqnarray*}
\Delta \log (\psi_{n,I}) (e)\!\! &=&\!\! \int_{\Gamma} \log 
\Big( \frac{\nu_{n} (g(I))}{\nu_n(I)} \Big) \thinspace d\mu(g) 
= \int_{\Gamma} \log \Big( \int_I Jac_{n}(g)(x) \thinspace \frac{d \nu_n(x)}{\nu_n(I)}\Big) 
\thinspace d\mu(g)\\
&\geq&\!\!\int_{\Gamma} \Big( \int_I \log \big( Jac_{n}(g)(x) \big) \thinspace 
\frac{d\nu_n(x)}{\nu_n(I)} \Big) \thinspace d\mu(g) 
= \int_I \Big( \int_{\Gamma} \log \big( Jac_{n}(g)(x) \big) \thinspace d\mu(g) \Big) 
\thinspace \frac{d\nu_n(x)}{\nu_n(I)} \geq -\esp\frac{1}{n}.
\end{eqnarray*}
Ceci est valable pour tout intervalle $I$ v\'erifiant 
$\bar{\nu}(I) > 0$. Par cons\'equent, et \`a cause des relations 
$$\psi_{n,I}(gf)= \psi_{n,f(I)}(g),$$ 
le laplacien de $\log(\psi_{n,I})$ est partout minor\'e 
par $-1/n$. Donc, si l'on d\'esigne par $\psi_I$ la limite simple de la suite des 
$\psi_{n_i,I}$, \esp {\em i.e.}  \thinspace $\psi_{I}(g) = \bar{\nu} \big( g(I) \big)$, 
\esp alors la fonction $\log (\psi_{I})$ est surharmonique, dans le sens que son laplacien 
est positif. Par ailleurs, puisque $\mu$ est sym\'etrique, $\psi_{I}$ est 
harmonique. Par cons\'equent, pour tout $f \in \Gamma$ nous avons partout 
\'egalit\'e dans les in\'egalit\'es suivantes~:
\[ \log (\psi_{I}) (f) \leq \int_{\Gamma} \log (\psi_{I}) (gf) \thinspace d\mu(g) 
\leq \log \Big( \int_{\Gamma} \psi_{I}(gf)\thinspace d\mu(g)\Big)=\log (\psi_{I}) (f).\]
La fonction $\psi_I$ est donc constante. Or, comme l'intervalle $I$ v\'erifiant \esp $\bar{\nu}(I)>0$ 
\esp est arbitraire, nous en d\'eduisons que la mesure $\bar{\nu}$ est invariante par tous les 
\'el\'ements de $\Gamma$. La proposition \ref{neg} est d\'emontr\'ee.


\begin{footnotesize}

\vspace{0.1cm}

Bertrand Deroin\\

IH\'ES, 35 route de Chartres, 91440 Bures sur Yvette, France (bderoin@ihes.fr)\\

\vspace{0.2cm}

Victor Kleptsyn\\

UMPA, ENS-Lyon, 46 all\'ee d'Italie, 69007 Lyon, France (victor.kleptsyn@umpa.ens-lyon.fr)\\

\vspace{0.2cm}

Andr\'es Navas\\

IH\'ES, 35 route de Chartres, 91440 Bures sur Yvette, France (anavas@ihes.fr)\\

Univ. de Chile, Las Palmeras 3425, \~Nu\~noa, Santiago, Chile (andnavas@uchile.cl)\\

\end{footnotesize}

\end{document}